\documentclass[12pt]{article}
\usepackage{amssymb,bm}
\usepackage{amsmath}
\usepackage{amsthm}
\usepackage{graphicx}
\usepackage[active]{srcltx} % Includes line number info into dvi file
\usepackage{hyperref} % for xdvi viewing
\hypersetup{pdfborder=0 0 0}
\newtheorem{theorem}{{\sc Theorem}}[section]

\newtheorem{definition}[theorem]{Definition}

\newcommand{\bb}[1]{\mathbb{ #1}}

\bmdefine\Bone{1}
\newcommand{\thus}{\Rightarrow}
\newcommand{\End}{\mathrm{End}}
\newcommand{\Ann}{\mathrm{Ann}}
\newcommand{\dev}[1]{\mathrm{dev}(#1)}

\newcommand{\Sym}{\mathrm{Sym}}

\newcommand{\Skew}{\mathrm{Skew}}
\newcommand{\Span}{\mathrm{Span}}
\newcommand{\defeq}{{\buildrel\rm def\over=}}

\newcommand{\eqv}{\Longleftrightarrow}
\newcommand{\bra}[1]{\overline{#1}}

\newcommand{\Trc}{\mathrm{Tr}\,}

\newcommand{\cof}{\mathrm{cof}}

\newcommand{\tns}[1]{#1\otimes #1}
\newcommand{\hf}{\displaystyle\frac{1}{2}}
\newcommand{\nth}[1]{\displaystyle\frac{1}{#1}}

\newcommand{\dif}[2]{\displaystyle\frac{\partial #1}{\partial #2}}
\newcommand{\Grad}{\nabla}
\newcommand{\Div}{\nabla \cdot}
\newcommand{\Curl}{\nabla \times}

\renewcommand{\Hat}[1]{\widehat{#1}}
\newcommand{\Tld}[1]{\widetilde{#1}}

\newcommand{\av}[1]{\langle #1 \rangle}

\def\XXint#1#2#3{{\setbox0=\hbox{$#1{#2#3}{\int}$ }
\vcenter{\hbox{$#2#3$ }}\kern-.6\wd0}}

\newcommand{\im}{\mathfrak{Im}}
\newcommand{\re}{\Re\mathfrak{e}}
\newcommand{\jump}[1]{\lbrack\!\lbrack #1 \rbrack\!\rbrack}

\newcommand{\mat}[4]{\left[\begin{array}{cc}
\displaystyle{#1}&\displaystyle{#2}\\[1ex]
\displaystyle{#3}&\displaystyle{#4}\end{array}\right]}
\newcommand{\vect}[2]{\left[\begin{array}{c}
\displaystyle{#1}\\[1ex]\displaystyle{#2}\end{array}\right]}

\newcommand{\bc}{boundary condition}

\newcommand{\rhs}{right-hand side}
\newcommand{\lhs}{left-hand side}
\newcommand{\mc}{microstructure}

\newcommand{\WLOG}{without loss of generality}
\newcommand{\nbh}{neighborhood}
\newcommand{\IFF}{if and only if }

\newcommand{\Ga}{\alpha}
\newcommand{\Gb}{\beta}
\newcommand{\Gd}{\delta}
\newcommand{\Ge}{\epsilon}

\newcommand{\Gvf}{\varphi}
\newcommand{\Gg}{\gamma}

\newcommand{\Gk}{\kappa}

\newcommand{\Gl}{\lambda}

\newcommand{\Gth}{\theta}

\newcommand{\Gs}{\sigma}

\newcommand{\GD}{\Delta}

\newcommand{\GL}{\Lambda}

\newcommand{\GS}{\Sigma}

\newcommand{\GO}{\Omega}

\bmdefine\BGa{\alpha}
\bmdefine\BGb{\beta}
\bmdefine\BGd{\delta}
\bmdefine\BGe{\epsilon}
\bmdefine\BGve{\varepsilon}
\bmdefine\BGf{\phi}
\bmdefine\BGvf{\varphi}
\bmdefine\BGg{\gamma}
\bmdefine\BGc{\chi}
\bmdefine\BGi{\iota}
\bmdefine\BGk{\kappa}
\bmdefine\BGl{\lambda}
\bmdefine\BGn{\eta}
\bmdefine\BGm{\mu}
\bmdefine\BGv{\nu}
\bmdefine\BGp{\pi}
\bmdefine\BGth{\theta}
\bmdefine\BGvth{\vartheta}
\bmdefine\BGr{\rho}
\bmdefine\BGvr{\varrho}
\bmdefine\BGs{\sigma}
\bmdefine\BGvs{\varsigma}
\bmdefine\BGt{\tau}
\bmdefine\BGj{\tau}
\bmdefine\BGu{\upsilon}
\bmdefine\BGo{\omega}
\bmdefine\BGx{\xi}
\bmdefine\BGy{\psi}
\bmdefine\BGz{\zeta}
\bmdefine\BGD{\Delta}
\bmdefine\BGF{\Phi}
\bmdefine\BGG{\Gamma}
\bmdefine\BGL{\Lambda}
\bmdefine\BGP{\Pi}
\bmdefine\BGT{\Theta}
\bmdefine\BGS{\Sigma}
\bmdefine\BGU{\Upsilon}
\bmdefine\BGO{\Omega}
\bmdefine\BGX{\Xi}
\bmdefine\BGY{\Psi}

\bmdefine\BFM{\mathfrak{M}}
\bmdefine\BFb{\mathfrak{b}}
\bmdefine\BFk{\mathfrak{k}}
\bmdefine\BFm{\mathfrak{m}}
\bmdefine\BFu{\mathfrak{u}}
\bmdefine\BFv{\mathfrak{v}}

\newcommand{\CA}{{\mathcal A}}

\newcommand{\ClD}{{\mathcal D}}
\newcommand{\CE}{{\mathcal E}}

\newcommand{\CH}{{\mathcal H}}
\newcommand{\CI}{{\mathcal I}}
\newcommand{\CJ}{{\mathcal J}}
\newcommand{\CK}{{\mathcal K}}
\newcommand{\CL}{{\mathcal L}}

\newcommand{\CT}{{\mathcal T}}

\bmdefine\BCA{{\mathcal A}}
\bmdefine\BCB{{\mathcal B}}
\bmdefine\BCC{{\mathcal C}}
\bmdefine\BCD{{\mathcal D}}
\bmdefine\BCE{{\mathcal E}}
\bmdefine\BCF{{\mathcal F}}
\bmdefine\BCG{{\mathcal G}}
\bmdefine\BCH{{\mathcal H}}
\bmdefine\BCI{{\mathcal I}}
\bmdefine\BCJ{{\mathcal J}}
\bmdefine\BCK{{\mathcal K}}
\bmdefine\BCL{{\mathcal L}}
\bmdefine\BCM{{\mathcal M}}
\bmdefine\BCN{{\mathcal N}}
\bmdefine\BCO{{\mathcal O}}
\bmdefine\BCP{{\mathcal P}}
\bmdefine\BCQ{{\mathcal Q}}
\bmdefine\BCR{{\mathcal R}}
\bmdefine\BCS{{\mathcal S}}
\bmdefine\BCT{{\mathcal T}}
\bmdefine\BCU{{\mathcal U}}
\bmdefine\BCV{{\mathcal V}}
\bmdefine\BCW{{\mathcal W}}
\bmdefine\BCX{{\mathcal X}}
\bmdefine\BCY{{\mathcal Y}}
\bmdefine\BCZ{{\mathcal Z}}

\bmdefine\Bzr{ 0}
\bmdefine\Ba{ a}
\bmdefine\Bb{ b}
\bmdefine\Bc{ c}
\bmdefine\Bd{ d}
\bmdefine\Be{ e}
\bmdefine\Bf{ f}
\bmdefine\Bg{ g}
\bmdefine\Bh{ h}
\bmdefine\Bi{ i}
\bmdefine\Bj{ j}
\bmdefine\Bk{ k}
\bmdefine\Bl{ l}
\bmdefine\Bm{ m}
\bmdefine\Bn{ n}
\bmdefine\Bo{ o}
\bmdefine\Bp{ p}
\bmdefine\Bq{ q}
\bmdefine\Br{ r}
\bmdefine\Bs{ s}
\bmdefine\Bt{ t}
\bmdefine\Bu{ u}
\bmdefine\Bv{ v}
\bmdefine\Bw{ w}
\bmdefine\Bx{ x}
\bmdefine\By{ y}
\bmdefine\Bz{ z}
\bmdefine\BA{ A}
\bmdefine\BB{ B}
\bmdefine\BC{ C}
\bmdefine\BD{ D}
\bmdefine\BE{ E}
\bmdefine\BF{ F}
\bmdefine\BG{ G}
\bmdefine\BH{ H}
\bmdefine\BI{ I}
\bmdefine\BJ{ J}
\bmdefine\BK{ K}
\bmdefine\BL{ L}
\bmdefine\BM{ M}
\bmdefine\BN{ N}
\bmdefine\BO{ O}
\bmdefine\BP{ P}
\bmdefine\BQ{ Q}
\bmdefine\BR{ R}
\bmdefine\BS{ S}
\bmdefine\BT{ T}
\bmdefine\BU{ U}
\bmdefine\BV{ V}
\bmdefine\BW{ W}
\bmdefine\BX{ X}
\bmdefine\BY{ Y}
\bmdefine\BZ{ Z}

\newcommand{\SFA}{\mathsf{A}}
\newcommand{\SFB}{\mathsf{B}}
\newcommand{\SFC}{\mathsf{C}}

\newcommand{\SFF}{\mathsf{F}}

\newcommand{\SFI}{\mathsf{I}}
\newcommand{\SFJ}{\mathsf{J}}
\newcommand{\SFK}{\mathsf{K}}
\newcommand{\SFL}{\mathsf{L}}
\newcommand{\SFM}{\mathsf{M}}

\newcommand{\SFP}{\mathsf{P}}

\newcommand{\SFS}{\mathsf{S}}
\newcommand{\SFT}{\mathsf{T}}

\newcommand{\SFX}{\mathsf{X}}

\usepackage{color}
\author{Yury Grabovsky}
\title{Exact relations and links for two-dimensional thermoelectric composites}
\begin{document}
\maketitle 
\tableofcontents

\section{Nonintroduction}
This is a report of the massive multi-year effort by the author and two
graduate students Huilin Chen and Sarah Childs to compute all exact relations
and links for two-dimensional thermoelectric composites. The size of this
report is due to the inclusion of all technical details of calculations, which
are customarily omitted in journal articles. At the moment I have no time to
prepare a proper ``archival quality'' manuscript with a good introduction and
references. However, I believe that the results, concisely summarized in the
last three sections of this report, should be made available to the
research community even in this unfinished form.

\section{Equations of thermoelectricity}
Thermoelectric properties of a material are described by the relations between
the gradient $\Grad\mu$ of an electrochemical potential, temperature gradient
$\Grad T$, current density $\Bj_{E}$ and entropy flux $\Bj_{S}$.
The total energy $U=U(S,N)$ is a function of entropy and the number of charge
carriers $N$. Therefore, the energy flux $\dot{U}$ is given by
\[
\dot{U}=T\dot{S}+\mu\dot{N},\quad T=\dif{U}{S},\quad\mu=\dif{U}{N}.
\]
where $T$ is the absolute temperature and $\mu$ is the electrochemical potential.
Thus, in a general heterogeneous medium we have
\[
\Bj_{U}=T\Bj_{S}+\mu\Bj_{E},
\]
where $\Bj_{U}$ is the total energy flux, $\Bj_{S}$ is the entropy flux and
$\Bj_{E}$ is the electric current (charge crrier flux). 
The consevation of charge and energy laws are expressed by the equations
\begin{equation}
  \label{conslaws}
  \Div\Bj_{E}=0,\qquad \Div\Bj_{U}=0.
\end{equation}
In addition to conservation laws we also postulate linear constitutive laws
that relate the electric current and the entropy flux to the nonuniformity of
electrochemical potential and temperature. In a thermoelectric material these
two driving forces are coupled:
\begin{equation}
  \label{constitlaw}
  \begin{cases}
      \Bj_{E}=\BGs\Grad(-\mu)+\BGs\BS\Grad(-T),\\
      \Bj_{S}=\BS^{T}\BGs\Grad(-\mu)+\BGg\Grad(-T)/T,
  \end{cases}\qquad\BGs^{T}=\BGs,\quad\BGg^{T}=\BGg.
\end{equation}
The Onsager reciprocity relation is incorporated in the above constitutive laws.
The form of the cross-property coupling tensors is chosen in such a way as to
make the thermoelectric coupling laws more transparent.  We will now show how
the general equations (\ref{conslaws}), (\ref{constitlaw}) relate to the
well-known thermoelectric effects.

\subsection{Seebeck effect and the figure of merit}
The electrochemical potential $\mu$ is a sum of the electrostatic potential
and a chemical potential. The latter depends only on the temperature and is
therefore constant, when the temperature is constant. In this case
$\BE=\Grad(-\mu)$ is the electric field and the first equation in
(\ref{constitlaw}) reads $\Bj_{E}=\BGs\BE$. Therefore, $\BGs$ has the physical
meaning of the isothermal conductivity tensor. As such it must be represented
by a symmetric, positive definite $3\times 3$ matrix. In the absence of the
electrical current ($\Bj_{E}=\Bzr$) the gradient of $-\mu$ has the meaning of
the electromotive force generated by a temperature gradient. This is called
the \emph{Seebeck effect}. From the first equation in (\ref{constitlaw}) we obtain
\[
\Be_{\rm emf}=-\Grad\mu=\BS\Grad T.
\]
The $3\times 3$ matrix $\BS$ is called the Seebeck coefficient (tensor). In
the literature the Seebeck coefficient is often assumed to be a
scalar. However, we will see that in general, a composite made of such
materials will have an anisotropic Seeback coefficient. Another a priori assumption
is that $\BS$ is symmetric (see e.g. [lusi18]). We will again see that
symmetry of $\BS$ is not preserved under homogenization.

The heat flux at zero electric current is characterized by the heat
conductivity tensor $\Bj_{U}=-\BGk\Grad T$, which gives a formula for the
tensor $\BGg$ in the constitutive equtions in terms of the symmetric, positive
definite heat conductivity tensor $\BGk$:
\[
\BGk=\BGg-T\BS^{T}\BGs\BS.
\]

Thus, imposing a temperature gradient on a thermoelectic material creates stored
electrical energy with density
\[
e_{\rm el}=\BGs\Be_{\rm emf}\cdot\Be_{\rm emf}=(\BS^{T}\BGs\BS\Grad T)\cdot(\Grad T).
\]
This phenomenon can be used to make a ``Seebeck generator'', converting heat
flux (temperature differences) directly into electrical energy. The efficiency
of Seebeck generator is called the \emph{figure of merit}.

The body not in thermal equilibrium can be used to produce mechanical
work. However, not all thermal energy can be used. One of the physical
interpretations of entropy is that it is a measure of the inaccessible portion
of the total internal energy per degree of temperature. Thus,
the density of this non-extractable thermal energy is the product of
temperature and the entropy production density:
\[
e_{\rm th}=T\Div\Bj_{S}=\Div(T\Bj_{S})-\Bj_{S}\cdot\Grad T=
\Div\Bj_{U}-\frac{\Bj_{U}}{T}\cdot\Grad T=\nth{T}(\BGk\Grad T)\cdot\Grad T.
\]
In a thermoelectric device we want to maximize the stored electrical energy
while minimizing unusable thermal energy.
The ratio $E_{\rm electrical}/E_{\rm entropy}$ is therefore a measure of
efficiency of the
thermoelectric \emph{device}, since the values of energies depend on specific
\bc s. If we want a \emph{material property} that is independent of the \bc s
we may define the figure of merit as follows
\[
Z=\sup_{\Bh}\frac{\BS^{T}\BGs\BS\Bh\cdot\Bh}{\BGk\Bh\cdot\Bh}.
\]
Thus, $Z$ is the largest eigenvalue of $\BGk^{-1}\BS^{T}\BGs\BS$.

In the isotropic case, where $\BGs$, $\BS$, and $\BGk$ are all
constant multiples of the identity, we have
$Z=S^{2}\Gs/\Gk$. 

In summary, our assumptions on the
possible values of the tensors $\BGs$, $\BS$ and $\BGg$ are equivalent to the
symmetry and positive definiteness of the $6\times 6$ matrix
\begin{equation}
  \label{physL}
  \SFL'=\mat{\BGs}{\BGs\BS}{\BS^{T}\BGs}{\BGk/T+\BS^{T}\BGs\BS},
\end{equation}
that describes constitutive relation (\ref{constitlaw})
\[
\vect{\Bj_{E}}{\Bj_{S}}=\SFL'\vect{\Grad(-\mu)}{\Grad(-T)}.
\]

\subsection{The Thomson and Peltier effects}
In physically relevant variables we can write equations of thermoelectricity
in the form of the following system
\begin{equation}
  \label{physeq}
  \begin{cases}
    \Bj_{E}=\BGs\Grad(-\mu)+\BGs\BS\Grad(-T),\\
    \Bj_{Q}=T\BS^{T}\Bj_{E}+\BGk\Grad(-T),\\
    \Bj_{U}=\Bj_{Q}+\mu\Bj_{E},\\
    \Div\Bj_{E}=\Div\Bj_{U}=0,
  \end{cases}
\end{equation}
where $\Bj_{Q}=T\Bj_{S}$ is the heat flux. In this form it is immediately
apparent that adding a constant to the electrochemical potential $\mu$ does
not change the flux $\Bj_{E}$, while adding a constant multiple of $\Bj_{E}$
to $\Bj_{U}$. Since $\Div\Bj_{E}=0$ then adding a constant to the
electrochemical potential $\mu$ gives another solution of balance equations
(\ref{physeq}). This observation will be useful later.

Let us write the conservation of energy law:
\[
0=\Div\Bj_{U}=\Div(\BGk\Grad(-T))+\Div(T\BS^{T}\Bj_{E})+\Grad\mu\cdot\Bj_{E},
\]
where we have used the conservation of charge law $\Div\Bj_{E}=0$.
From the first equation in (\ref{physeq}) we have
\[
\Grad\mu=-\BGs^{-1}\Bj_{E}+\BS\Grad(-T),
\]
so that the conservation of energy has the form
\[
0=\Div(\BGk\Grad(-T))+\Div(T\BS^{T}\Bj_{E})-(\BGs^{-1}\Bj_{E})\cdot\Bj_{E}-(\BS\Grad T)\cdot\Bj_{E}.
\]
We can rewrite it as
\begin{equation}
  \label{enerbal}
  \Div(\BGk\Grad(-T))=(\BGs^{-1}\Bj_{E})\cdot\Bj_{E}-T\Div(\BS^{T}\Bj_{E}).
\end{equation}
On the left we have heat production density. On the right we have two heat
sources: Joule heating, represented by the first term and the thermoelectric
heating or cooling. This second term represents those thermoelectric effects
that occur when the current flows through the thermoelectric material. The
commonly encountered description of these effects assumes that the Seebeck
tensor is scalar: $\BS=S\BI_{3}$. In that case the conservation of charge law
allows us to simplify the second term on the \rhs\ of (\ref{enerbal}):
\[
\dot{Q}_{\rm thel}=-T\Div(\BS^{T}\Bj_{E})=-T(\Grad S)\cdot\Bj_{E}.
\]
The \emph{Thompson effect} is related to the dependence of the Seebeck
coefficient $S$ on $T$. In this case the additional thermoelectric heat
production density is
\[
\dot{Q}_{\rm thel}=-T\Grad S\cdot\Bj_{E}=-TS'(T)\Grad T\cdot\Bj_{E}=-\CK\Grad T\cdot\Bj_{E} .
\]
The coefficient $\CK=TS'(T)$ is called the Thompson coefficient. The
\emph{Peltier effect} occurs at an isothermal junction $\GS$ of two
different materials with different Seebeck coefficients. At every point $\Bs\in\GS$
\[
\dot{Q}_{\rm thel}=-T\Grad S\cdot\Bj_{E}=-T\jump{S}(\Bj_{E}\cdot\Bn)\Gd_{\Bs}(\Bx)=-\jump{\Pi}(\Bj_{E}\cdot\Bn)\Gd_{\Bs}(\Bx),
\]
where $\Pi=TS$ is called the Peltier coefficient and the normal charge flux
$\Bj_{E}\cdot\Bn$ is continuous across the junction $\GS$. 
In general, when $\BS$ is not scalar we can rewrite the thermoelectric heat
production term $\dot{Q}_{\rm thel}$ as
follows
\begin{equation}
  \label{anisiS}
\dot{Q}_{\rm thel}=-T(\Div\BS)\cdot\Bj_{E}-T\av{\dev{\BS},\Grad\Bj_{E}},
\end{equation}
where
\[
\dev{\BS}=\BS-\nth{3}(\Trc\BS)\BI_{3}
\]
is the deviatoric part of $\BS$. The second term in (\ref{anisiS}) represents
the thermoelectric effects of anisotropy, of which there seems to be no
evidence in the literature. 

In conclusion, aside from the completely undocumented anisotropic effects from
(\ref{anisiS}) the most commonly described effects are due to either
inhomogeneity (Peltier effect) or essential nonlinearity (Thomson's effect due
to the dependence of $\BS$ on $T$). In what follows we will focus on the case
of small perturbations 
\[
\mu=\mu_{0}+\Ge\Tld{\mu},\qquad T=T_{0}+\Ge\Tld{T}.
\]
As $\Ge\to 0$ the equations of theremoelectricity become linear with respect
to $\Tld{\mu}$ and $\Tld{T}$, with temperature-dependent coefficients set to
their values corresponding to $T=T_{0}$. In what follows we use notation $\mu$
and $T$ instead of $\Tld{\mu}$ and $\Tld{T}$.

\subsection{The canonical form of equations of thermoelectricity}
Many physical phenomena and processes are described by systems of linear PDE
(partial differential equations). A very large class of these have a common
structure that I would like to emphasise. These phenomena deal with various
properties of solid bodies (materials). For example, we may be interested in
how materials respond to electromagnetic fields, heat or mechanical forces. In
each of these cases we identify a pair of vector fields, defined at each point
inside the material and taking values in appropriate vector spaces (different
in each physical context). The first field in the pair describes what is being
done to a material: applied deformation, or an electric field, or a temperature
distribution, etc. The second describes how material responds to the applied
field, such as forces (stress) that arise in response to a deformation or an
electrical current that arises in response to an applied electric field, etc.
These physical vector fields obey fundamental laws of classical physics,
such as conservation of energy, for example. These laws can be expressed as
system of linear differential equations, which combined with the constitutive
laws give a full quantitative description of the respective phenomena.
 
The constitutive law is a linear relation between the two fields in a
pair. The linear operator effecting this relation describes material
properties in question. If one adds information of how the disturbance is
applied to the body (usually through a particular action on the boundary of
the body), then one obtains a unique solution. To summarize, we will be
looking at
\begin{itemize}
\item A pair of vector fields (we will call them $\BE(\Bx)$ and $\BJ(\Bx)$)
  defined at each point $\Bx$ inside the body $\GO$, with values
  in some finite dimensional vector space, equipped with a physically natural
  inner product;
\item Systems of constant coefficient PDEs obeyed by $\BE(\Bx)$ and $\BJ(\Bx)$
\item A linear relation between $\BE(\Bx)$ and $\BJ(\Bx)$, written in operator
  form $\BJ(\Bx)=\SFL(\Bx)\BE(\Bx)$, where the linear operator $\SFL(\Bx)$
  describes material properties (that can be different at different points
  $\Bx\in\GO$). This operator is almost always symmeteric and positive
  definite.
\end{itemize}
We do not include \bc s in the above list because answers to questions that we
are interested in do not depend on \bc s. We will now show how equations of
thermoelectricity (\ref{physeq}) can be rewritten as a linear relation between
a pair of curl-free fields $(\Be_{1},\Be_{2})$ and a pair of divergence-free
fields $(\Bj_{1},\Bj_{2})$
\begin{equation}
  \label{constrel}
  \begin{cases}
  \Bj_{1}=\BL_{11}\Be_{1}+\BL_{12}\Be_{2},\\
\Bj_{2}=\BL_{12}^{T}\Be_{1}+\BL_{22}\Be_{2}.
\end{cases}
\end{equation}
Following Callen's textbook we define new
potentials
\[
\psi_{1}=\frac{\mu}{T},\qquad\psi_{2}=\nth{T},
\]
denoting
\[
\Be_{1}=\Grad\psi_{1},\quad\Be_{2}=\Grad\psi_{2},\quad\Bj_{1}=-\Bj_{E},\quad\Bj_{2}=\Bj_{U},
\]
we obtain the form (\ref{constrel}), were
\begin{equation}
  \label{mathL}
  \BL_{11}=T\BGs,\quad\BL_{12}=-T(\mu\BGs+T\BGs\BS),\quad
\BL_{22}=T[\mu^{2}\BGs+T\BGg+T\mu(\BGs\BS+\BS^{T}\BGs)].
\end{equation}
In general the coefficients $\BL_{ij}$ depend on the values of $T$ and $\mu$,
and we are considering situations where these quantities change little and
equations (\ref{constrel}) represent the linearization around the fixed values
$T_{0}$ and $\mu_{0}$. We observe that the new material tensor
\[
\SFL=T\mat{\BGs}{-\BGs(\mu+T\BS)}{-(\mu+T\BS)^{T}\BGs}
{T\BGk+(\mu+T\BS)^{T}\BGs(\mu+T\BS)}
\]
is symmetric and positive definite \IFF $\SFL'$, given by (\ref{physL}), is
symmetric and positive definite, i.e. \IFF $\BGs$ and $\BGk$ are symmetric and
positive definite $3\times 3$ matrices. In full generality equations of
thermoelectricity are very nonlinear, especially in view of the fact that all
physical property tensors $\BGs$, $\BGk$ and $\BS$ depend on temperature
$T$. We will be working with linearized version of the equations where both
$\mu$ and $T$ do not vary a lot. Mathematically, we look at the leading order
asymptotics of solutions $(\mu,T)$ of the form $\mu=\mu_{0}+\Ge\Tld{\mu}$ and
$T=T_{0}+\Ge\Tld{T}$. We have already observed that the full thermoelectric
system is invariant with respect to addition of a constant to the
electrochemical potential $\mu$. Thus, modifying the potential $\psi_{1}$
\[
\psi_{1}\mapsto\frac{\mu-\mu_{0}}{T},
\] 
we can set $\mu_{0}=0$, \WLOG. Thus, for linearized problems we can write
\begin{equation}
  \label{Lcanon}
  \SFL=T_{0}^{2}\mat{\BGs/T_{0}}{-\BGs\BS}{-\BS^{T}\BGs}{\BGk+T_{0}\BS^{T}\BGs\BS},
\end{equation}
where all physical property tensors $\BGs$, $\BGk$ and $\BS$ are evaluated at
$T=T_{0}$---the working temperature. We note (for no particular reason other
than curiousity) that
\begin{equation}
  \label{Lcanoninv}
  \SFL^{-1}=\nth{T_{0}}\mat{\BGs^{-1}+T_{0}\BS\BGk^{-1}\BS^{T}}{-\BS\BGk^{-1}}{-\BGk^{-1}\BS^{T}}{\BGk^{-1}/T_{0}}.
\end{equation}

Now, the vector fields $\BE=(\Be_{1},\Be_{2})$ and $\BJ=(\Bj_{1},\Bj_{2})$
take their values in the $2d$-dimesnional vector space
$\CT=\bb{R}^{d}\oplus\bb{R}^{d}$, $d=2$ or 3. (It will be 2 in this paper.)
The natural inner product on $\CT$ is defined by
\[
(\BE,\BE')_{\CT}=\Be_{1}\cdot\Be'_{1}+\Be_{2}\cdot\Be'_{2}.
\]
The differential equations satisfied by $\BE$ and $\BJ$ are
\begin{equation}
  \label{PDE}
  \Curl\Be_{1}=\Curl\Be_{2}=0,\qquad\Div\Bj_{1}=\Div\Bj_{2}=0.
\end{equation}
The material properties tensor $\SFL(\Bx)$ can therefore be written as a
$2\times 2$ block matrix
\begin{equation}
  \label{BML}
  \SFL(\Bx)=\mat{\BL_{11}(\Bx)}{\BL_{12}(\Bx)}{\BL_{12}^{T}(\Bx)}{\BL_{22}(\Bx)},
\end{equation}
where $\BL_{11}$ and $\BL_{22}$ are symmetric (and positive definite) $3\times
3$ matrices. The constitutive relation $\BJ=\SFL\BE$ can also be written as 
$\BJ=\SFL\BE$. From the block-components of $\SFL$ we can recover the physical
tensors:
\begin{equation}
  \label{L2phys}
  \BGs=\Gb_{0}\BL_{11},\quad\BS=-\Gb_{0}\BL_{11}^{-1}\BL_{12},\quad
\BGk=\Gb_{0}^{2}(\BL_{22}-\BL_{12}^{T}\BL_{11}^{-1}\BL_{12}),\quad\Gb_{0}=\nth{T_{0}}.
\end{equation}
With these formulas the figure of merit form is
\begin{equation}
  \label{ZTL}
  ZT=\max_{\Bh}\frac{\BL_{12}^{T}\BL_{11}^{-1}\BL_{12}\Bh\cdot\Bh}{(\BL_{22}-\BL_{12}^{T}\BL_{11}^{-1}\BL_{12})\Bh\cdot\Bh}=\frac{\Gl}{1-\Gl},
\end{equation}
where $\Gl\in(0,1)$ is the largest eigenvalue of
$\BL_{22}^{-1}\BL_{12}^{T}\BL_{11}^{-1}\BL_{12}$. 

For isotropic materials $\SFL=\BL\otimes\BI_{3}$ and their figure of merit is
\[
ZT=\frac{L_{12}^{2}}{\det\BL},\qquad\BL=\mat{L_{11}}{L_{12}}{L_{12}}{L_{22}}.
\]

\section{Periodic composites}
Let $Q=[0,1]^{d}$. It is a unit square when $d=2$ and unit cube when
$d=3$. Let us suppose that $Q$ is divided into two complementary subsets $A$
and $B$. We place one thermoelectric material in $A$ and another in $B$. If
the corresponding tensors of material properties and denoted by $\SFL_{A}$ and
$\SFL_{B}$, then the function
\[
\SFL(\Bx)=\SFL_{A}\chi_{A}(\Bx)+\SFL_{B}\chi_{B}(\Bx)
\]
describes this situation mathematically, since $\SFL(\Bx)=\SFL_{A}$, \IFF
$\Bx\in A$ and $\SFL(\Bx)=\SFL_{B}$, \IFF $\Bx\in B$. Here $\chi_{S}(\Bx)$ is
the characteristic function of a subset $S$, taking value 1, when $\Bx\in S$
and value 0, otherwise.

Now we are going to tile the entire space $\bb{R}^{d}$ with the copies of the
``period cell'' $Q$, generating a $Q$-periodic function $\SFL_{\rm per}(\Bx)$,
$\Bx\in\bb{R}^{d}$. Specifically, in order to find the value of $\SFL_{\rm
  per}(\Bx)$ at a specific point $\Bx\in\bb{R}^{d}$ we first find a vector
$\Bz$ with integer components, such that $\Bx-\Bz\in Q$ and then define
$\SFL_{\rm per}(\Bx)=\SFL(\Bx-\Bz)$. In general $\SFL_{\rm
  per}(\Bx_{1})=\SFL_{\rm per}(\Bx_{2})$, whenever $\Bx_{1}-\Bx_{2}$ has
integer components.

A periodic composite material would have such a structure on a
\emph{microscopic level}. Mathematically, we choose $\Ge>0$, representing a
microscopic length scale and define $\SFL_{\Ge}(\Bx)=\SFL_{\rm per}(\Bx/\Ge)$,
restricting $\Bx$ to lie in a subset $\GO\subset\bb{R}^{d}$ occupied by our
composite. On a macroscopic level, such a composite will look as though it is
a homogeneous thermoelectric material. Its thermoelectric tensor $\SFL^{*}$,
called the effective tensor of the composite, is a complicated function not
only of the tensors $\SFL_{A}$ and $\SFL_{B}$ of its constituents, but also of
the set $A$ ($B=Q\setminus A$). Specifically, if we keep $\SFL_{A}$ and
$\SFL_{B}$ fixed and change only the shape of $A$, then the effective tensor
$\SFL^{*}$ will change as well. Understanding how $\SFL^{*}$ depends on the
shape of $A$ is an important (and difficult) problem, that could help design
thermoelectric composites with desired properties. Even though, there is a
mathematical description of $\SFL^{*}$ as a function of $A$, it is complicated
and we will not be needing or using this description.

\section{Exact relations}
Let us recall that the thermoelectric tensor $\SFL$ is a $2\times 2$
block-matrix
\begin{equation}
  \label{blckM}
  \SFL=\mat{\BL_{11}}{\BL_{12}}{\BL_{12}^{T}}{\BL_{22}},
\end{equation}
where $\BL_{11}$ and $\BL_{22}$ are symmetric $3\times 3$ matrices.
Therefore, we are going to think of each such tensor as a point in an
$N$-dimensional vector space, where $N=2d^{2}+d$.

Now, let us imagine that we have fixed two such points, representing tensors
$\SFL_{A}$ and $\SFL_{B}$ and we are making periodic composites with all
possible subsets $A\subset Q$. For each choice of the set $A$ we get a point 
$\SFL^{*}$ in our $N$-dimensional vector space. The set of all such points
corresponding to all possible subsets $A\subset Q$ is called the G-closure of
of a two-point set $\{\SFL_{A},\SFL_{B}\}$. Generically, this G-closure set
will have a non-empty interior is the $N$-dimensional vector space of material
tensors. However, there are special cases, all of which we want to describe,
where the G-closure set is a submanifold of positive codimension.
Equations describing such a submanifold are
called exact relations. In the language of composite materials, these relations will be
satisfied by \emph{all} composites, as long as they are made of materials that
satisfy these equations.

\section{Polycrystals}
A general thermoelectric tensor $\SFL$ is \emph{anisotropic}, i.e. its $N$
components will change when we rotate the material. Nevertheless, there are
\emph{isotropic} materials, whose tensors are given by
\begin{equation}
  \label{Liso3d}
  \SFL=\mat{\Gl_{11}\BI_{d}}{\Gl_{12}\BI_{d}}{\Gl_{12}\BI_{d}}{\Gl_{22}\BI_{d}}
  =\BGL\otimes\BI_{d},\qquad\BGL=\mat{\Gl_{11}}{\Gl_{12}}{\Gl_{12}}{\Gl_{22}},
\end{equation}
when $d=3$. When $d=2$ there is an additional isotropic tensor
\begin{equation}
  \label{Liso2d}
  \SFL=\BGL\otimes\BI_{2}+\nu\BR_{\perp}\otimes\BR_{\perp},\qquad\BR_{\perp}=\mat{0}{-1}{1}{0}.
\end{equation}
However, if we think that the 2D case is just a special case of 3D, where
fields do not change in one of the direction, then the isotropy (\ref{Liso2d})
can be exhibited by anisotropic thermoelectrics that are, for example, only
transversely isotropic.

Operators (\ref{Liso3d}) are positive definite \IFF $\Gl_{11}>0$ and
$\det\BGL>0$, while operators (\ref{Liso2d}) are positive definite \IFF $\Gl_{11}>0$ and
$\det\BGL>\nu^{2}$.

If tensors $\SFL_{A}$ and $\SFL_{B}$
are anisotropic, it means that in a composite described above we have to use
these materials in one fixed orientation. This is very often impractical, and
we will restrict our attention to polycrystals, where we are permitted to use
each anisotropic material in any orientation, so that at different points we
may have different orientation of the same material. There are a lot fewer
exact relations and links for polycrystals, and they will be easier (not easy)
to find.

\section{Exact relations for thermoelectricity}
Recall that in space dimension $d$ the space $\CT=\bb{R}^{d}\oplus\bb{R}^{d}$
is $2d$-dimensional and the space $\Sym(\CT)$ of all symmetric operators on
$\CT$ is $N=2d(2d+1)/2=2d^{2}+d$ dimensional. A positive definite operator
$\SFL\in\Sym(\CT)$ will be thought of as a description of thermoelectric
properties of a material via (\ref{constrel}), (\ref{BML}) and will be
referred to as a \emph{tensor of material properties} or a
\emph{thermoelectric tensor}.  The set of all thermoelectric tensors, i.e. the
set of all positive definite symmetric operators on $\CT$ will be denoted
$\Sym^{+}(\CT)$. Our first task is to identify all \emph{exact
  relations}---submanifolds $\bb{M}$ (think surfaces or curves in space) in
$\Sym^{+}(\CT)$, such that the thermoelectric tensor of any composite made with
materials from $\bb{M}$ must necessarily be in $\bb{M}$. To be precise we are
only interested in polycrystalline exact relations $\bb{M}$ that have the
additional property that $\BR\cdot\SFL\in\bb{M}$ for any $\SFL\in\bb{M}$ and
for any rotation $\BR\in SO(d)$. In fact, the complete list of them is known
for $d=3$. Our first goal will be to compute all polycrystalline exact
relations when $d=2$. This is done by applying the general theory of exact
relations that states that every exact relation $\bb{M}$ corresponds to a peculiar algebraic
object called \emph{Jordan multialgebra}. Jordan algebras are very-well
studied object in algebra. The particle ``multi'' comes from the fact that in
our case each Jordan algebra carries several Jordan multiplications,
parametrized by a particular subspace $\CA\subset\Sym(\CT)$.
\begin{definition}
  We say that a subspace $\Pi\subset\Sym(\CT)$ is a Jordan $\CA$-multialgebra
  if 
\[
\SFK_{1}*_{\SFA}\SFK_{2}=\hf(\SFK_{1}\SFA\SFK_{2}+\SFK_{2}\SFA\SFK_{1})\in\Pi,\quad
\forall\SFK\in\Pi,\ \SFA\in\CA.
\]
\end{definition}
The subspace $\CA$ of Jordan multiplications is defined by the formula
\begin{equation}
  \label{Adef}
  \CA=\Span\{\BGG_{0}(\Bn)-\BGG_{0}(\Bn_{0}):|\Bn|=1\},
\end{equation}
where $\BGG_{0}(\Bn)$ is accociated to an isotropic tensor $\SFL_{0}$, through
which the exact relations manifold $\bb{M}$ is passing and is determined by
the differential equations (\ref{PDE}) satified by the fields $\BE$ and $\BJ$, written in
Fourier space
\begin{equation}
  \label{FPDE}
  \BGx\times\Hat{\Be_{1}}=\BGx\times\Hat{\Be_{2}}=0,\qquad
\BGx\cdot\Hat{\Bj_{1}}=\BGx\cdot\Hat{\Bj_{2}}=0.
\end{equation}
We view these equations as definitions of two subspaces $\CE_{\BGx}$ and
$\CJ_{\BGx}$ of pairs $(\Hat{\Be_{1}},\Hat{\Be_{2}})$ and
$(\Hat{\Bj_{1}},\Hat{\Bj_{2}})$, respectively, regarding the Fourier wave
vector $\BGx$ as fixed. Specifically,
\[
\CE_{\BGx}=\{(\Gg_{1}\BGx,\Gg_{2}\BGx):\Gg_{1}\in\bb{R},\ \Gg_{1}\in\bb{R}\},\qquad
\CJ_{\BGx}=\{(\Bv_{1},\Bv_{2})\in\bb{R}^{d}\oplus\bb{R}^{d}:\BGx\cdot\Bv_{1}=\BGx\cdot\Bv_{2}=0\}.
\]
We observe that vectors $\BGx$ and $c\BGx$, where $c\in\bb{R}\setminus\{0\}$
produce the same subspaces $\CE_{c\BGx}-\CE_{\BGx}$ and
$\CJ_{c\BGx}=\CJ_{\BGx}$. Therefore, we only need to refer to subspaces
$\CE_{\Bn}$ and $\CJ_{\Bn}$ for unit vectors $\Bn$.

Now let $\SFL_{0}\in\Sym(\CT)$ be isotropic (and positive definite), then we define
\begin{equation}
  \label{G0def}
  \BGG_{0}(\Bn)=\SFL_{0}^{-1}\BGG'(\Bn),
\end{equation}
where $\BGG'(\Bn)$ is the projection onto $\SFL_{0}\CE_{\Bn}$ along
$\CJ_{\Bn}$. 

In order to compute $\BGG_{0}(\Bn)$ we take an arbitrary vector
$(\Bu_{1},\Bu_{2})\in\CT$ and decompose it into the sum
\[
(\Bu_{1},\Bu_{2})=\SFL_{0}\BE+\BJ,\qquad\BE\in\CE_{\Bn},\quad\BJ\in\CJ_{\Bn}.
\]
Then $\SFL_{0}\BE=\BGG'(\Bn)(\Bu_{1},\Bu_{2})$ and therefore, 
\[
\BE=\SFL_{0}^{-1}\BGG'(\Bn)(\Bu_{1},\Bu_{2})=\BGG_{0}(\Bn)(\Bu_{1},\Bu_{2}).
\]
The vector $\BE\in\CE_{\Bn}$ is uniquely determined by two scalars $\Gg_{1}$, $\Gg_{2}$:
$\BE=(\Gg_{1}\Bn,\Gg_{2}\Bn)$, while $\BJ=(\Bj_{1},\Bj_{2})$ must satisfy 
\begin{equation}
  \label{Geq}
  \Bj_{1}\cdot\Bn=0,\qquad\Bj_{2}\cdot\Bn=0.
\end{equation}
Finding expressions for $(\Bj_{1},\Bj_{2})$ from
\[
\BJ=(\Bu_{1},\Bu_{2})-\SFL_{0}(\Gg_{1}\Bn,\Gg_{2}\Bn),
\]
where $\SFL_{0}$ is given by (\ref{Liso3d}) or (\ref{Liso2d}), and substituting
into (\ref{Geq}) we will obtain two linear equations for the two unknowns
$\Gg_{1}$, $\Gg_{2}$. Solving this linear system we will obtain the explicit
expressions for $\Gg_{1}$, $\Gg_{2}$ in terms of $\Bu_{1}$, $\Bu_{2}$, $\Bn$
and $\SFL_{0}$. The obtained expressions will be linear in $\Bu_{1}$, $\Bu_{2}$,
permitting us to write the desired operator $\BGG_{0}(\Bn)$ in block-matrix
form (\ref{BML}).
\begin{equation}
  \label{G0}
  \BGG_{0}(\Bn)=\BGL^{-1}\otimes(\tns{\Bn}).
\end{equation}
Formula (\ref{G0}) is valid in both cases $d=2$ and $d=3$. We can now use
formula (\ref{G0}) in (\ref{Adef}) and obtain
the explicit formula for the subspace $\CA$:
\begin{equation}
  \label{Athel}
  \CA=\{\BGL^{-1}\otimes\BA:\BA^{T}=\BA,\Trc\BA=0\}.
\end{equation}
Our first task is to identify
(explicitly) all SO(d)-invariant Jordan $\CA$-multialgebras
$\Pi\subset\Sym(\CT)$. Once this is done, the theory of exact relations, gives
an explicit formula for the corresponding exact
relation $\bb{M}$
\begin{equation}
  \label{Mdef}
  \bb{M}=\{\SFL\in\Sym^{+}(\CT): W_{\Bn}(\SFL)\in\Pi\}
\end{equation}
for some unit vector $\Bn$, where
\[
W_{\Bn}(\SFL)=[(\SFL-\SFL_{0})^{-1}+\BGG_{0}(\Bn)]^{-1}.
\]
We ephasize that even though transformations $W_{\Bn}$ are all different for
different $\Bn$, the submanifold $\bb{M}$ in (\ref{Mdef}) does not depend on
the choice of $\Bn$. In fact, we can also compute $\bb{M}$ using the transformation
\[
W_{\SFM}(\SFL)=[(\SFL-\SFL_{0})^{-1}+\SFM]^{-1},
\]
where the ``inversion key'' $\SFM$ is found as the ``simplest'' isotropic
tensor satisfying
\begin{equation}
  \label{invkeydef}
  \SFK_{1}*_{\BGG_{0}(\Bn)-\SFM}\SFK_{2}\in\Pi,\quad\forall\SFK\in\Pi.
\end{equation}

At this point we note that the subspace $\CA$ is different for different
isotropic reference tensors $\SFL_{0}$ through which the exact relations
manifolds are passing. In many cases, and in ours in particular, this technical
complication can be eliminated by means of the  ``covariance
transformations''. The idea is to observe that for any invertible operators
$\SFB$ and $\SFC$ on $\CT$ we have
\[
\SFB(\SFK_{1}*_{\SFA}\SFK_{2})\SFC=(\SFB\SFK_{1}\SFC)*_{\SFC^{-1}\SFA\SFB^{-1}}(\SFB\SFK_{2}\SFC).
\]
It means that if $\Pi$ is an Jordan $\CA$-multialgebra, then $\SFB\Pi\SFC$ is
a Jordan $\SFC^{-1}\SFA\SFB^{-1}$-multialgebra. In order to preserve symmetry
of operators and $SO(d)$-invariance of subspaces we have to set $\SFB=\SFC^{T}$ and use
only isotropic operators $\SFC$. In the case of $\CA$, given by (\ref{Athel})
we can use $\SFC=\BGL^{-1/2}\otimes\BI_{d}$, so that
\[
\CA_{0}=\SFC^{-1}\SFA\SFC^{-1}=\{\BI_{2}\otimes\BA:\BA^{T}=\BA,\Trc\BA=0\}
\]
is independent of $\SFL_{0}$. Now, if $\Pi_{0}$ is a Jordan
$\CA_{0}$-multialgebra then we compute a corresponding inversion key
$\SFM_{0}$, which must be an isotropic tensor satisfying
\begin{equation}
  \label{invkeyeq}
  \SFK_{1}*_{\BI_{2}\otimes\BI_{d}-d\SFM_{0}}\SFK_{2}\in\Pi_{0},\quad\forall\SFK\in\Pi_{0}.
\end{equation}
In particular, the choice 
\begin{equation}
  \label{M0univ}
  \SFM_{0}=\nth{d}\BI_{2}\otimes\BI_{d}=\nth{d}\CI_{\CT}
\end{equation}
satisfies (\ref{invkeyeq}). When $d=2$ we will also try two other simpler
choices for $\SFM$: $\SFM=0$ and $\SFM=\hf\BI_{2}\otimes(\tns{\Be_{1}})$.
Once the inversion key $\SFM_{0}$ is determined, the corresponding exact
relation $\bb{M}$ will be computed using
\[
\bb{M}=\{\SFC^{-1}\SFL\SFC^{-1}:\SFL\in\bb{M}_{0}\},\qquad
\bb{M}_{0}=\{\SFL\in\Sym^{+}(\CT): W_{0}(\SFL)\in\Pi_{0}\},
\]
where $\SFC=\BGL^{-1/2}\otimes\BI_{d}$ and
\[
W_{0}(\SFL)=[(\SFL-\SFL_{0}^{0})^{-1}+\SFM_{0}]^{-1},
\]
where
\[
\SFL_{0}^{0}=\SFC\SFL_{0}\SFC=
\begin{cases}
  \BI_{2}\otimes\BI_{3},&d=3,\\
 \BI_{2}\otimes\BI_{2}+\frac{\nu}{\sqrt{\det\BGL}}\tns{\BR_{\perp}},&d=2.
\end{cases}
\]
In summary, our first task is to solve a (very nontrivial) problem of identifying
all SO(2)-invariant subspaced $\Pi_{0}\subset\Sym(\bb{R}^{2}\oplus\bb{R}^{2})$,
that are Jordan $\CA_{0}$-multialgebras. Very often a difficult problem can be
made easier by identifying its symmetries. In our case a symmetry is
an SO(2)-invariant linear operator $\Phi:\Sym(\CT)\to\Sym(\CT)$, such that
\begin{equation}
  \label{autodef}
  \Phi(\SFK\SFA\SFK)=\Phi(\SFK)\SFA\Phi(\SFK),\quad\forall\SFK\in\Sym(\CT),\ \SFA\in\CA_{0}.
\end{equation}
Such a transformation will be called a gobal SO(2)-invariant Jordan
$\CA_{0}$-multialgebra automorphism of $\Sym(\CT)$.

\section{SO(2)-invariant subspaces of $\Sym(\CT)$}
Our task of finding SO(2)-invariant Jordan $\CA_{0}$-multialgebras will be
significantly simplified by first identifying all SO(2)-invariant subspaces of
$\Sym(\CT)$, a standard task in the representation theory of compact Lie
groups, which is particularly easy for a commutative ``circle group'' $SO(2)$.
It is well know that all irreducible representations of $SO(2)$ are of complex
type. Therefore, it will be convenient to identify the physical space $\bb{R}^{2}$
with complex numbers, so that\footnote{The image in $\bb{C}$ of a vector in
  $\bb{R}^{2}$, denoted by a bold letter, is represented by the same letter in
  normal font.}  $\Bx=(x_{1},x_{2})\mapsto x=x_{1}+ix_{2}\in\bb{C}$. Then
\begin{equation}
  \label{TCn}
  \CT=\bb{R}^{2}\oplus\bb{R}^{2}\cong\bb{C}\oplus\bb{C}\cong\bb{C}^{2},
\end{equation}
With corresponding identification
\[
\CT\ni\vect{\Bu}{\Bv}\mapsto(u,v)\in\bb{C}^{2},\quad u=u_{1}+iu_{2},\ v=v_{1}+iv_{2}.
\]

The utility of this isomorphism of $4$-dimensional real vector spaces ($\CT$
and $\bb{C}^{2}$) comes from the fact that the set $\bb{C}^{2}$ also has a
structure of a complex vector space. In order to characterize all rotationally
invariant subspaces in $\Sym(\CT)$ we observe that rotations $\BR_{\Gth}$ of
$\bb{R}^{2}$ through the angle $\Gth$ counterclockwise act on vectors
$\Bu\in\CT\cong\bb{C}^{2}$ by $\BR_{\Gth}\cdot\Bu=e^{i\Gth}\Bu$. Every
real operator $\SFK$ on $\CT$ can be described by two complex $2\times 2$ matrices
$X$ and $Y$ via
\begin{equation}
  \label{Ku}
  \SFK\Bu=X\Bu+Y\bra{\Bu},\qquad\Bu\in\bb{C}^{2},
\end{equation}
where $\Bu$ on the \lhs\ is an element of $\CT$, while $\Bu$ on the \rhs\ is
its $\bb{C}^{2}$ representation. Henceforth, we will write $K(X,Y)$ to
indicate this parametrization of $\End(\CT)$. 

We compute
\[
\left(\vect{\Bu_{1}}{\Bu_{2}},\vect{\Bv_{1}}{\Bv_{2}}\right)_{\CT}=
\Bu_{1}\cdot\Bv_{1}+\Bu_{2}\cdot\Bv_{2}=\re(u_{1}\bra{v_{1}})+\re(u_{2}\bra{v_{2}})=
\re(\Bu,\Bv)_{\bb{C}^{2}}.
\]
We compute
\[
(K(X,Y)\Bu,\Bv)_{\bb{C}^{2}}=(X\Bu+Y\bra{\Bu},\Bv)_{\bb{C}^{2}}=
(\Bu,X^{H}\Bv)_{\bb{C}^{2}}+\bra{(\Bu,\bra{Y^{H}}\bra{\Bv})_{\bb{C}^{2}}},
\]
where $X^{H}=\bra{X}^{T}$ denotes Hermitian conjugation\footnote{We do not use
  the standard notation $X^{\ast}$ to avoid confusion with our notation for
  the effective tensor.}.  Hence
\[
(K(X,Y)\Bu,\Bv)_{\CT}=\re(\Bu,X^{H}\Bv)_{\bb{C}^{2}}+\re(\Bu,\bra{Y^{H}}\bra{\Bv})_{\bb{C}^{2}}=
\re(\Bu,X^{H}\Bv+Y^{T}\bra{\Bv})_{\bb{C}^{2}}=(\Bu,\SFK^{T}\Bv)_{\CT}.
\]
This shows that $K(X,Y)^{T}=K(X^{H},Y^{T})$. It follows that
$K(X,Y)\in\Sym(\CT)$ \IFF $X$ is a complex Hermitian $2\times 2$ matrix
($X^{H}=X$) and $Y$ is a complex symmetric $2\times 2$ matrix ($Y^{T}=Y$).

Let us find the characterization of positive definiteness of
$K(X,Y)\in\Sym(\CT)$ in terms of complex matrices $X$ and $Y$. The first
observation is that
\[
(K(X,Y)\Bu,\Bu)_{\CT}=\hf\left(\Hat{K}(X,Y)\vect{\Bu}{\bra{\Bu}},
\vect{\Bu}{\bra{\Bu}}\right)_{\bb{C}^{4}},\qquad\Hat{K}(X,Y)=\mat{X}{Y}{\bra{Y}}{\bra{X}}.
\]
We see that $\Hat{K}(X,Y)\in\mathfrak{H}(\bb{C}^{4})$. We now view
$\bb{C}^{4}=\bb{C}^{2}\oplus\bb{C}^{2}$ as a real (8-dimensional) vector
space with the standard inner product
$(\BGx,\BGn)=\re(\BGx,\BGn)_{\bb{C}^{4}}$. It is then easy to check that
$\bb{C}^{4}$ can be split into the orthogonal sum of subspaces
$\bb{C}^{4}=S_{+}\oplus S_{-}$,
\[
S_{\pm}=\left\{\vect{\Bu}{\pm\bra{\Bu}}:\Bu\in\bb{C}^{2}\right\}.
\]
Moreover, both $S_{+}$ and $S_{-}$ are invariant subspaces for
$\Hat{K}(X,Y)$. The final observation is that $S_{-}=iS_{+}$. Now, the
positive definiteness of $K(X,Y)$ is equivelent to the positive definiteness
of $\Hat{K}(X,Y)$ on $S_{+}$. But $i\BGx\in S_{+}$ for any $\BGx\in
S_{-}$, and therefore,
\[
(\Hat{K}(X,Y)\BGx,\BGx)_{\bb{C}^{4}}=(\Hat{K}(X,Y)i\BGx,i\BGx)_{\bb{C}^{4}}>0.
\]
This implies that the positive definiteness
of $\Hat{K}(X,Y)$ on $S_{+}$ is equivalent to the positive definiteness
of $\Hat{K}(X,Y)$ on $\bb{C}^{4}$. In turn, the positive definiteness
of $\Hat{K}(X,Y)$ on $\bb{C}^{4}$ is equivalent to
\[
X>0,\qquad S_{X}=X-Y\bra{X}^{-1}\bra{Y}>0.
\]
We will see later that
\[
K(X,Y)^{-1}=K(S_{X}^{-1},-S_{X}^{-1}Y\bra{X}^{-1}).
\]
In other words, positive definiteness of $K(X,Y)$ is equivalent to positive
definiteness of the ``X-components'' of both $K(X,Y)$ and $K(X,Y)^{-1}$.

We easily compute the action of rotations $\BR_{\Gth}$ on $K(X,Y)$ from the
``distributive law''
\[
\BR_{\Gth}\cdot(K(X,Y)\Bu)=(\BR_{\Gth}\cdot K(X,Y))(\BR_{\Gth}\cdot\Bu)
\]
and the formula $\BR_{\Gth}\cdot\Bu=e^{i\Gth}\Bu$:
\[
e^{i\Gth}(X\Bu+Y\bra{\Bu})=(\BR_{\Gth}\cdot K(X,Y))e^{i\Gth}\Bu.
\]
Denoting $e^{i\Gth}\Bu$ be $\Bv$ and substituting $\Bu=e^{-i\Gth}\Bv$ we
obtain
\[
(\BR_{\Gth}\cdot K(X,Y))\Bv=X\Bv+e^{2i\Gth}Y\bra{\Bv},
\]
which means that
\begin{equation}
  \label{Raction}
   \BR_{\Gth}\cdot K(X,Y)=K(X,e^{2i\Gth}Y).
\end{equation}
Therefore, if $\Pi$ is an SO(2)-invariant subspace of $\Sym(\CT)$ then
\[
\Pi=\CL(V,W)\defeq\{K(X,Y):X\in V\subset\CH(\bb{C}^{2}),\ Y\in W\subset\Sym(\bb{C}^{2})\},
\]
where $V$ can be any subspace of $\CH(\bb{C}^{2})$|the set of all complex
Hermitian $2\times 2$ matrices, regarded as a real vector space, and $W$ can
be any subspace of $\Sym(\bb{C}^{2})$|the set of all complex symmetric
$2\times 2$ matrices, regarded as a complex vector space. In this notation the
subspace $\CA_{0}$ corresponds to $V=\{0\}$ and $W=\{z\BI_{2}:z\in\bb{C}\}$:
\[
\CA_{0}=\{K(0,z\BI_{2}):z\in\bb{C}\}.
\]

\section{SO(2)-invariant Jordan $\CA_{0}$-multialgebras}
Using definition (\ref{Ku}) of the action of an operator $\SFK$ we compute
\begin{equation}
  \label{multrule}
  K(X_{1},Y_{1})K(X_{2},Y_{2})=K(X_{1}X_{2}+Y_{1}\bra{Y_{2}},X_{1}Y_{2}+Y_{1}\bra{X_{2}}),
\end{equation}
Using this  multiplication rule we compute
\[
K(X,Y)K(0,z\BI_{2})K(X,Y)=K(zX\bra{Y}+\bar{z}YX,zX\bra{X}+\bar{z}Y^{2})
\]
This formula implies that
a subspace $\Pi=\CL(V,W)$ is a Jordan $\CA_{0}$-multialgebra \IFF
\begin{equation}
    \label{2dJMA}
Y^{2}+XX^{T}\in W,\quad YX+XY^{H}\in V\text{ for all }X\in V,\ Y\in W.
  \end{equation}
  The goal is therefore find all solutions $\CL(V,W)$ of
  (\ref{2dJMA}). Equations (\ref{2dJMA}) suggest an obvious strategy. We first
  identify all 0, 1, 2 and 3-dimensional complex subspace
  $W\subset\Sym(\bb{C}^{2})$ satisfying $Y^{2}\in W$ for all $Y\in W$. Then,
  for each such $W$ we will look for 0, 1, 2, 3 and 4-dimensional subspaces
  $V\subset\mathfrak{H}(\bb{C}^{2})$, the space of complex Hermitian $2\times
  2$ matrices. However, before we begin it will be helpful to identify all
  symmetries of (\ref{2dJMA}), i.e. all global SO(2)-invariant Jordan
  $\CA_{0}$-multialgebra automorphisms.

\section{Global SO(2)-invariant Jordan $\CA_{0}$-multialgebra automorphisms}
Let $\Phi:\Sym(\CT)\to\Sym(\CT)$ be SO(2)-invariant. Any such linear map must
have the form
\[
\Phi(K(X,Y))=K(\Phi_{11}(X)+\Phi_{12}(Y),\Phi_{21}(X)+\Phi_{22}(Y)),
\]
where $\Phi_{ij}$ are real linear maps between the appropriate spaces
Then the ``distributive law'' for rotations says
\[
\BR_{\Gth}\cdot\Phi(K(X,Y))=\Phi(\BR_{\Gth}\cdot K(X,Y))
\]
Using formula (\ref{Raction}) we obtain
\[
K(\Phi_{11}(X)+\Phi_{12}(Y),e^{2i\Gth}\Phi_{21}(X)+e^{2i\Gth}\Phi_{22}(Y))=
K(\Phi_{11}(X)+\Phi_{12}(e^{2i\Gth}Y),\Phi_{21}(X)+\Phi_{22}(e^{2i\Gth}Y))
\]
It follows that
\[
\Phi_{12}(Y)=\Phi_{12}(e^{2i\Gth}Y),\quad
e^{2i\Gth}\Phi_{21}(X)=\Phi_{21}(X),\quad
e^{2i\Gth}\Phi_{22}(Y)=\Phi_{22}(e^{2i\Gth}Y).
\]
The first two equations imply that $\Phi_{12}=0$ and $\Phi_{21}=0$, while the
third equation implies that $\Phi_{22}$ is a complex-linear map on
$\Sym(\bb{C}^{2})$. Thus, any linear $SO(2)$ automorphism $\Phi$ of
$\Sym(\CT)$ can be written as 
\[
\Phi(K(X,Y))=K(\Phi_{0}(X),\Phi_{2}(Y)),
\]
where $\Phi_{0}$ is a real-linear automorphism  of
$\mathfrak{H}(\bb{C}^{2})$ and $\Phi_{2}$ is a complex-linear automorphism  of
$\Sym(\bb{C}^{2})$.  

Let us now assume that $\Phi$ is also a Jordan $\CA_{0}$-multialgebra
automorphism. In that case the maps $\Phi_{0}$ and $\Phi_{2}$ must satisfy,
additionally,
\begin{equation}
  \label{PhiXX}
\Phi_{2}(XX^{T})=\Phi_{0}(X)\Phi_{0}(X)^{T},\qquad\Phi_{2}(Y^{2})=\Phi_{2}(Y)^{2},  
\end{equation}
\begin{equation}
  \label{PhiXY}
\Phi_{0}(X\bra{Y}+YX)=\Phi_{0}(X)\bra{\Phi_{2}(Y)}+\Phi_{2}(Y)\Phi_{0}(X).  
\end{equation}
Our task is to determine all maps $\Phi_{0}$ and $\Phi_{2}$, satisfying
(\ref{PhiXX}), (\ref{PhiXY}).

Observe that $\Phi_{2}(Y^{2})=\Phi_{2}(Y)^{2}$ means $\Phi_{2}$ maps
projections (idempotents) into projections. Conversely, if $\Phi_{2}(Y)$ is a
projection for some
$Y\in\Sym(\bb{C}^{2})$, then $\Phi_{2}(Y^{2})=\Phi_{2}(Y)$, which implies that
$Y^{2}=Y$, since $\Phi_{2}$ is a bijection. Every non-zero idempotent in
$\Sym(\bb{C}^{2})$ is either $I_{2}$ or $\Ba\otimes\Ba$, where
$\Ba\cdot\Ba=1$. Since $\Phi_{2}$ is a bijection it must map all idempotents
of the form $\Ba\otimes\Ba$, except
possibly one, into idempotents of the same form. Then the map
\[
\Ba\mapsto\Trc\Phi_{2}(\Ba\otimes\Ba)
\]
is continuous and has constant value 1 on almost all $\Ba$. Hence, by
continuity it must have value 1 on all $\Ba$. This implies that
$\Phi_{2}(I_{2})=I_{2}$ and $\Phi_{2}(\Ba\otimes\Ba)=\tns{C\Ba}$ for
some complex linear map $C$ that has the property
$C\Ba\cdot C\Ba=\Ba\cdot\Ba=1$. Thus, 
$\Phi_{2}(Y)=CYC^{T}$, where $C\in
O(2,\bb{C})=\{C\in\End_{\bb{C}}(\bb{C}^{2}):CC^{T}=I_{2}\}$.

Now we need to compute the map $\Phi_{0}$. We start by determining all
possible values of $\Phi_{0}(I_{2})$. To this end we take $X=I_{2}$ in the
first equation in (\ref{PhiXX}). Then
$\Phi_{0}(I_{2})\Phi_{0}(I_{2})^{T}=I_{2}$. Hence, $\Phi_{0}(I_{2})\in
O(2,\bb{C})\cap\mathfrak{H}(\bb{C}^{2})$.  Let us take $X=I_{2}$, $Y=iS$,
$S\in\Sym(\bb{R}^{2})$ in (\ref{PhiXY}). Then
\[
0=-i\Phi_{0}(I)\bra{C}S\bra{C}^{T}+iCSC^{T}\Phi_{0}(I).
\]
Using the fact that $C\in O(2,\bb{C})$ we obtain that
$C^{T}\Phi_{0}(I)\bra{C}$ commutes with every $S\in\Sym(\bb{R}^{2})$. This
quickly leads, via $\Phi_{0}(I)\Phi_{0}(I)^{T}=I$, to $\Phi_{0}(I)=\pm
CC^{H}$. Notice that if $\Phi_{0}$ satisfies our equations, then so does
$-\Phi_{0}$. Hence, \WLOG\ we assume that $\Phi_{0}(I)=CC^{H}$. 

Next we determine $\Phi_{0}$ on real symmetric matrices.
Taking $X=I$, $Y=S\in\Sym(\bb{R}^{2})$ in (\ref{PhiXY}) we then obtain
$\Phi_{0}(S)=CSC^{H}$. It remains to figure out the value of $\Phi_{0}$ on $i\BR_{\perp}$.

Now we take $X=S_{1}$, $Y=iS_{2}$ in (\ref{PhiXY}),
where $\{S_{1},S_{2}\}\subset\Sym(\bb{R}^{2})$. Then
\[
\Phi_{0}(i[S_{2},S_{1}])=-iCS_{1}C^{H}\bra{C}S_{2}C^{H}+iCS_{2}C^{T}CS_{1}C^{H}=
iC[S_{2},S_{1}]C^{H}.
\]
Hence, $\Phi_{0}(X)=CXC^{H}$ for all $X\in\mathfrak{H}(\bb{C}^{2})$. Thus the
set of all $SO(2)$ Jordan $\CA_{0}$-multialgebra automorphisms is given by
\begin{equation}
  \label{Jautoex}
\Phi(K(X,Y))=K(\pm CXC^{H},CYC^{T}),\qquad C\in O(2,\bb{C}).  
\end{equation}
Finally every $C\in O(2,\bb{C})$ has a representation
\[
C=C_{+}=\mat{\cos z}{\sin z}{-\sin z}{\cos z},\text{ or }
C=C_{-}=\mat{\cos z}{\sin z}{\sin z}{-\cos z},\qquad z\in\bb{C}.
\]

In fact, there is a general theorem that guarantees that the set of all
$SO(d)$-invariant Jordan multialgebra automorphisms has the form
$\Phi(\SFX)=\SFC\SFX\SFC^{T}$ for any isotropic $\SFC$ preserving $\SFA$,
$\SFC\Tld{\BGG}_{0}\SFC^{T}=\Tld{\BGG}_{0}$. In addition to these there could be
additional automorphisms of the form $\Phi(\SFX)=-\SFC\SFX\SFC^{T}$ for those $\SFC$ for
which $\SFC\Tld{\BGG}_{0}\SFC^{T}=-\Tld{\BGG}_{0}$ (which can only happen when
$\Tld{\BGG}_{0}$ has the same number of positive and negative
eigenvalues). From this and (\ref{multrule}) it is easy to get the general
form obtained above (using the fact that any isotropic tensor $\SFC$ has the
form $K(C,0)$).

\section{Describing all Jordan $\CA_{0}$-multialgebras}
\label{sec:JMA}
Complex subspace $W\subset\Sym(\bb{C}^{2})$ can have dimension 0, 1, 2 or 3. 
\begin{itemize}
\item $\dim W=0$. Then  $W=\{0\}$ and $V$ may contain only those $X$ for
  which $X^{T}X=0$. If $X\not=0$, then $X$ is rank 1 and Hermitian. Thus,
  $X=\Ba\otimes\bra{\Ba}$ for some $\Ba\in\bb{C}^{2}$, satisfying
  $\Ba\cdot\Ba=0$, i.e. $a_{1}^{2}=-a_{2}^{2}$, which is equivalent to
  $a_{2}=\pm ia_{1}$. Thus, all such $X$ must be real
  multiples of one of the following 2 matrices
\[
X_{1}=\mat{1}{i}{-i}{1},\qquad X_{2}=\bra{X_{1}}.
\]
Hence, either $V=\{0\}$ or $V=\bb{R}X_{j}$ for some $j=1,2$.
These two are isomorphic by $C=\mat{0}{1}{1}{0}\in O(2,\bb{C})$ that map $W$ into
    itslef and maps $X_{1}$ into $X_{2}$ and vice versa.
\item  $\dim W=1$. If $W$ contains an invertible matrix $Y$, then the
  Cayley-Hamilton theorem implies that $W$ contains $\BI_{2}$, since
\[
\BI_{2}=\frac{\Trc(Y)Y-Y^{2}}{\det(Y)}\in W.
\]
Thus we have two possibilities
\begin{itemize}
\item $W=\bb{C}\BI_{2}$. In that case $V$ must contain only
  matrices $X$ such that $XX^{T}=\Gl\BI_{2}$ for some $\Gl\in\bb{C}$. We compute that all
  such matrices $X$ must have one of two possible forms
\[
F_{1}(\Ga,\Gb)=\mat{\Ga}{i\Gb}{-i\Gb}{\Ga},\quad F_{2}(\Ga,\Gb)=\mat{\Ga}{\Gb}{\Gb}{-\Ga},\qquad
\{\Ga,\Gb\}\subset\bb{R}.
\]
Then $V$ is either $\{0\}$ or $V_{1}=\{F_{1}(\Ga,\Gb):\{\Ga,\Gb\}\subset\bb{R}\}$, 
$V_{2}=\{F_{2}(\Ga,\Gb):\{\Ga,\Gb\}\subset\bb{R}\}$, or any 1D subspace of $V_{1}$
or $V_{2}$.
\item $W=\bb{C}\tns{\Ba}$ for some $\Ba\in\bb{C}^{2}$. Then $V$ must contain only
  matrices $X$ such that $XX^{T}=\Gl\tns{\Ba}$. In particular, $X$ must be
  rank-1 (if it is non-zero). Thus, $V$ is either $\{0\}$ or
  $\bb{R}\Ba\otimes\bra{\Ba}$. In the latter case $\Pi$ is the annihilator of
  the 1D complex subspace $U$ of $\bb{C}^{2}$, where $U=\bb{C}\Ba^{\perp}$,
  where
\[
\Ba^{\perp}=\BR_{\perp}\Ba=(-a_{2},a_{1}).
\]
\end{itemize}
\item $\dim W=2$. Then $W$ must contain an invertible matrix. Indeed, the
  set of non-zero complex singular $2\times 2$ matrices is a 2D complex
  manifold and is not a subspace. Hence any 2D subspace cannot be contained in
  it. By Cayley-Hamilton $\BI_{2}\in W$. Let $W=\Span_{\bb{C}}\{\BI_{2},A\}$ for some
  $A\in\Sym(\bb{C}^{2})\setminus\{\Gl\BI_{2}\}$. Observe that \WLOG, we may assume that
  $A=\tns{\Ba}$ for some $\Ba\in\bb{C}^{2}$ (if $\Gl\in\bb{C}$ is an
  eigenvalue of $A$, then $A-\Gl\BI_{2}\in W$ has rank 1.) So,
\[
W=W_{\Ba}=\Span\{\BI_{2},\tns{\Ba}\},
\]
which obviously satisfies $Y^{2}\in W_{\Ba}$ for every $Y\in W_{\Ba}$. We can
apply the global automorphism $\Phi$ to the Jordan multialgebra $\Pi$ and
reduce $W$ to the algebra of complex $2\times 2$ diagonal matrices, if
$\Ba\cdot\Ba\not=0$ or to
\[
W=\left\{\mat{\Ga-\Gb}{\pm i\Gb}{\pm i\Gb}{\Ga+\Gb}:\{\Ga,\Gb\}\subset\bb{C}\right\}.
\]
 In order to compute all possible subspace $V$, we need to consider the two cases above separately
 \begin{itemize}
 \item $W$ is the algebra of complex $2\times 2$ diagonal matrices.
  \begin{itemize}
  \item $\dim V=0$, $V=\{0\}$
  \item $\dim V=1$. Then $V=\bb{R}H_{0}$ for some
    $H_{0}\in\mathfrak{H}(\bb{C}^{2})$. Condition $H_{0}H_{0}^{T}\in W$
    results in 4 possibilities for $H_{0}$:
    \begin{equation}
      \label{Hsqdiag}
      \mat{h_{1}}{0}{0}{h_{2}},\quad\mat{h}{i\Ga}{-i\Ga}{h},\quad
\mat{h}{\Ga}{\Ga}{-h},\quad\mat{0}{a}{\bra{a}}{0}.
    \end{equation}
If both components $H_{11}$ and $H_{22}$ are nonzero, or $H_{12}\not=0$,
then $\{YH_{0}+H_{0}\bra{Y}:Y\in W\}\subset V$ will be at least
two-dimensional. Hence, there are two
    possibilities: $H_{0}=\tns{\Be_{1}}$ and $H_{0}=\tns{\Be_{2}}$. These two
    are isomorphic by $C=\mat{0}{1}{1}{0}\in O(2,\bb{C})$ that map $W$ into
    itslef and maps $\tns{\Be_{1}}$ into $\tns{\Be_{2}}$.
  \item $\dim V\ge 2$. We notice that the set of all
    $H\in\mathfrak{H}(\bb{C}^{2})$, such that $HH^{T}$ is diagonal is the
    union of 4 two-dimensional vector spaces (\ref{Hsqdiag}). Thus, there are
    no solutions $V$ with dimension greater than 2, while solutions $V$ with
    $\dim V=2$ must be one of the 4 spaces in (\ref{Hsqdiag}). We only need to
    check which of the 4 subspaces $V$ in in (\ref{Hsqdiag}) have the property
    $\{YX+X\bra{Y}:Y\in W\}\subset V$ for any $X\in V$. It is easy to verify
    that only the first and the fourth ones have that property. Hence,
For $\dim V=2$ we have the following choices:
\begin{itemize}
\item[(a)] $V=\left\{\mat{\Ga}{0}{0}{\Gb}:\{\Ga,\Gb\}\subset\bb{R}\right\}$
\item[(b)] $V=\left\{\mat{0}{\bra{a}}{a}{0}:a\in\bb{C}\right\}$
\end{itemize}
\end{itemize}
\item $\displaystyle
W=\left\{\mat{\Ga-\Gb}{\pm i\Gb}{\pm i\Gb}{\Ga+\Gb}:\{\Ga,\Gb\}\subset\bb{C}\right\}.
$\\
We note that using $C=\mat{1}{0}{0}{-1}\in O(2,\bb{C})$ we can transform
``$-$'' sign above into the ``$+$'' sign. So, that \WLOG
\[
W=\left\{\mat{\Ga-\Gb}{i\Gb}{i\Gb}{\Ga+\Gb}:\{\Ga,\Gb\}\subset\bb{C}\right\}.
\]
  In this case Maple worksheet %\texttt{dimW2.mw} 
shows that
  \begin{itemize}
  \item $\dim V=0$, $V=\{0\}$
\item $\dim V=1$, $V=\bb{R}\mat{1}{i}{-i}{1}$
\item $\dim V=2$. There is a 1-parameter family of 2D subspaces permuted by the
  automorphism $\Phi$:
\[
V_{t}=\left\{\mat{x}{t(x-y)+i\frac{x+y}{2}}{t(x-y)-i\frac{x+y}{2}}{y}:\{x,y\}\subset\bb{R}\right\},
\]
together with 
\[
V_{\infty}=\left\{\mat{y}{x+iy}{x-iy}{y}:\{x,y\}\subset\bb{R}\right\},
\]
which we select to be the representative of the entire family.
\item $\dim V=3$,
\[
V=\left\{X\in\mathfrak{H}(\bb{C}^{2}):
\Trc X=2\im(X_{12})\right\}.
 \]
Is the only 3D solution.
  \end{itemize}
\end{itemize}
\item $\dim W=3$. Then $W=\Sym(\bb{C}^{2})$.\\ 
Let us assume that there is a non-zero $X\in V$. Then for any
$\{Y_{1},Y_{2}\}\subset W$ $X'=Y_{1}X+XY_{1}^{H}\in V$ and therefore,
$Y_{2}X'+X'Y_{2}^{H}\in V$. We compute
\[
Y_{2}X'+X'Y_{2}^{H}=
Y_{2}Y_{1}X+X(Y_{2}Y_{1})^{H}+Y_{2}XY_{1}^{H}+Y_{1}XY_{2}^{H}\in V.
\]
Switching $Y_{1}$ and $Y_{2}$ we also obtain
\[
Y_{1}Y_{2}X+X(Y_{1}Y_{2})^{H}+Y_{1}XY_{2}^{H}+Y_{2}XY_{1}^{H}\in V.
\]
Adding the two expressions we obtain
\[
(Y_{1}Y_{2}+Y_{2}Y_{1})X+X(Y_{1}Y_{2}+Y_{2}Y_{1})^{H}+2Y_{1}XY_{2}^{H}+2Y_{2}XY_{1}^{H}\in
V.
\]
But $Y_{1}Y_{2}+Y_{2}Y_{1}\in W$ for all $\{Y_{1},Y_{2}\}\subset W$ and
therefore,
\[
Y_{1}XY_{2}^{H}+Y_{2}XY_{1}^{H}\in V.
\]
for every $\{Y_{1},Y_{2}\}\subset W$.

Now let $Z\in \mathfrak{H}(\bb{C}^{2})$ be orthogonal to $V$. Then for every
$\{Y_{1},Y_{2}\}\subset W$ we must have
$\av{Y_{1}XY_{2}^{H}+Y_{2}XY_{1}^{H},Z}=0$. We compute
\[
0=\av{Y_{1}XY_{2}^{H}+Y_{2}XY_{1}^{H},Z}=\av{Y_{1},ZY_{2}X}+\av{Y_{1}^{H},XY_{2}^{H}Z}
\]
Restricting to $Y_{1}\in \Sym(\bb{R}^{2})$ we obtain
\[
\av{Y_{1},ZY_{2}X+XY_{2}^{H}Z}=0
\]
The matrix $M=ZY_{2}X+XY_{2}^{H}Z$ is self-adjoint and therefore has the form
$M_{1}+iM_{2}$, where $M_{1}\in\Sym(\bb{R}^{2})$,
$M_{2}\in\Skew(\bb{R}^{2})$. Thus, $\av{Y_{1},M_{1}}=0$ for every $Y_{1}\in
\Sym(\bb{R}^{2})$. Hence, $M_{1}=0$ and we conclide that for every $Y_{2}\in
\Sym(\bb{C}^{2})$
\[
ZY_{2}X+XY_{2}^{H}Z=\mat{0}{-i\Gb}{i\Gb}{0}
\]
for some $\Gb=\Gb(Y_{2})\in\bb{R}$. We repeat the same argument, now
restricting $Y_{1}$ to be of the form
$Y_{1}=iY_{0}$, $Y_{0}\in \Sym(\bb{R}^{2})$. In that case we obtain
\[
\av{Y_{0},ZY_{2}X-XY_{2}^{H}Z}=0,\quad\forall Y_{0}\in \Sym(\bb{R}^{2}).
\]
In this case the matrix $M=ZY_{2}X-XY_{2}^{H}Z$ is skew-adjoint and therefore
has the form $M=M_{1}+iM_{2}$, where $M_{1}\in\Skew(\bb{R}^{2})$,
$M_{2}\in\Sym(\bb{R}^{2})$ resulting in $\av{Y_{0},M_{2}}=0$ for all $Y_{0}\in
\Sym(\bb{R}^{2})$. It follows that $M_{2}=0$ and
\[
ZY_{2}X-XY_{2}^{H}Z=\mat{0}{-\Ga}{\Ga}{0}
\]
for some $\Ga=\Ga(Y_{2})\in\bb{R}$. Adding the two results we obtain that for
every $Y_{2}\in\Sym(\bb{C}^{2})$ there exists $b=b(Y_{2})\in\bb{C}$, such that
\[
ZY_{2}X=\mat{0}{-b}{b}{0}.
\]
We now choose $Y_{2}=\tns{\Ba}$ obtaining
\[
Z\Ba\otimes X^{T}\Ba=\mat{0}{-b}{b}{0}.
\]
The matrix on the \lhs\ has rank at most 1, while the matrix on the \rhs\ has
rank 2, unless $b=0$. Therefore, $Z\Ba\otimes X^{T}\Ba=0$ for every $\Ba\in\bb{C}^{2}$.
If $X\not=0$ then either $X^{T}\Be_{1}\not=0$ or $X^{T}\Be_{2}\not=0$. To fix
ideas suppose that $X^{T}\Be_{1}\not=0$. But then, for $\Ba=\Be_{1}$ we must
have $Z\Be_{1}=0$. If $Z\Be_{2}\not=0$ then, for $\Ba=\Be_{2}$, we must have that
$X^{T}\Be_{2}=0$. Now writing $\Ba=a_{1}\Be_{1}+a_{2}\Be_{2}$ we obtain
\[
0=Z(a_{1}\Be_{1}+a_{2}\Be_{2})\otimes X^{T}(a_{1}\Be_{1}+a_{2}\Be_{2})=a_{1}a_{2}Z\Be_{2}\otimes
X^{T}\Be_{1}.
\]
We can just choose $a_{1}=a_{2}=1$ and conclude, recalling that
$X^{T}\Be_{1}\not=0$, that $Z\Be_{2}=0$ in contradiction to our
assumption. Thus, we must have $Z=0$, implying that $V=\mathfrak{H}(\bb{C}^{2})$.
\end{itemize}
\newpage
\begin{center}
  \textbf{Summary}
\end{center}
It will be convenient to introduce the ``square-free'' vector
$\Bz_{0}=(1,-i)$. We order the 23 solutions by dimension of $(W,V)$ in
lexicographic order. We also give them short names for easy reference and identify families of
equivalent solutions.

~

  \begin{itemize}
  \item $W=\{0\}$
    \begin{itemize}
    \item $V=\{0\}$ \textcolor{red}{$(0,0)$}
    \item $V=\bb{R}\Bz_{0}\otimes\bra{\Bz_{0}}$ and the equivalent 
$V=\bb{R}\bra{\Bz_{0}}\otimes\Bz_{0}$ \textcolor{red}{$(0,\bb{R}\BZ_{0})$}$\sim(0,\bb{R}\bra{\BZ}_{0})$
    \end{itemize}
\item $W=\bb{C}\BI$. 
  \begin{itemize}
  \item $V=\{0\}$ \textcolor{red}{$(\bb{C}\BI,0)$}
  \item $V=V_{1}=\left\{\mat{\Ga}{i\Gb}{-i\Gb}{\Ga}:\{\Ga,\Gb\}\subset\bb{R}\right\}$,
\textcolor{red}{$(\bb{C}\BI,\BGF)$}
  \item $V=V_{2}=\left\{\mat{\Ga}{\Gb}{\Gb}{-\Ga}:\{\Ga,\Gb\}\subset\bb{R}\right\}$,
  \textcolor{red}{$(\bb{C}\BI,\BGY)$}
\item $V$ is any 1D subspace of $V_{1}$ or $V_{2}$, which can be ``rotated''
    by $C\in O(2,\bb{C})$ into one of the following subspaces
    \begin{itemize}
    \item $V=\bb{R}\BI_{2}$ \textcolor{red}{$(\bb{C}\BI,\bb{R}\BI)$}$\sim(\bb{C}\BI,\bb{R}\phi_{t})$,
$
\phi_{t}=\bb{R}\mat{\cosh t}{i\sinh t}{-i\sinh t}{\cosh t},
$
 $t\in\bb{R}$
   \item $V=\bb{R}\mat{0}{1}{1}{0}$ \textcolor{red}{$(\bb{C}\BI,\bb{R}\psi(i))$}$\sim(\bb{C}\BI,\bb{R}\psi(e^{i\Ga}))$, 
$
\psi(e^{i\Ga})=\bb{R}\mat{\cos\Ga}{\sin\Ga}{\sin\Ga}{-\cos\Ga},
$
$\Ga\in[0,\pi)$,
    \item $V=\bb{R}\mat{0}{i}{-i}{0}$ \textcolor{red}{$(\bb{C}\BI,i\BR_{\perp})$}$\sim(\bb{C}\BI,\bb{R}\phi'_{t})$,
$
\phi'_{t}=\bb{R}\mat{\sinh t}{i\cosh t}{-i\cosh t}{\sinh t},
$
 $t\in\bb{R}$
    \item $V=\bb{R}\Bz_{0}\otimes\bra{\Bz_{0}}$ \textcolor{red}{$(\bb{C}\BI,\bb{R}\BZ_{0})$}$\sim(\bb{C}\BI,\bb{R}\bra{\BZ}_{0})$
    \end{itemize}
  \end{itemize}
\item $W=\bb{C}\tns{\Be_{1}}\sim\bb{C}\tns{\Ba}$, if $\Ba\cdot\Ba=1$
  \begin{itemize}
  \item $V=\{0\}$ \textcolor{red}{$(\tns{\Be_{1}},0)$}$\sim(\tns{\Ba},0)$,
      $\Ba\cdot\Ba=1$, ($\pm\Ba$ defining the same subspaces)
  \item $V=\bb{R}\tns{\Be_{1}}\sim\bb{R}\Ba\otimes\bra{\Ba}$ (any vector
    $\Ba$, satisfying $\Ba\cdot\Ba=1$ can be rotated by $C\in O(2,\bb{C})$ into either
$\Be_{1}$. \textcolor{red}{${\rm Ann}(\bb{C}\Be_{2})$}$\sim{\rm
  Ann}(\bb{C}\bra{\Ba}^{\perp})$, where
$
\Ba^{\perp}=\BR_{\perp}\Ba=(-a_{2},a_{1}).
$
  \end{itemize}
\item $W=\bb{C}\tns{\Bz_{0}}\sim\bb{C}\tns{\bra{\Bz}_{0}}$
  \begin{itemize}
  \item $V=\{0\}$ \textcolor{red}{$(\tns{\Bz_{0}},0)$}$\sim(\tns{\bra{\Bz}_{0}},0)$,
  \item $V=\bb{R}\Bz_{0}\otimes\bra{\Bz}_{0}$ 
    \textcolor{red}{${\rm Ann}(\bb{C}\bra{\Bz_{0}})$}$\sim{\rm Ann}(\bb{C}\Bz_{0})$
  \end{itemize}
\item $W=\Span_{\bb{C}}\{\tns{\Be_{1}},\tns{\Be_{2}}\}$, (representing an
  infinite $O(2,\bb{C})$-orbit)
  \begin{itemize}
  \item $V=\{0\}$, \textcolor{red}{$(\ClD,0)$}$\sim(W_{\Ba},0)$,
$W_{\Ba}=\{Y\in\Sym(\bb{C}^{2}):Y\Ba\cdot\Ba^{\perp}=0\}$.
  \item $V=\bb{R}\Be_{1}\otimes\Be_{1}$,
\textcolor{red}{$(\ClD,\Be_{1}\otimes\Be_{1})$}$\sim(W_{\Ba},\bb{R}\Ba\otimes\bra{\Ba})$
\item $V=\Span_{\bb{R}}\{\tns{\Be_{1}},\tns{\Be_{2}}\}$
\textcolor{red}{$(\ClD,\ClD)$}$\sim(W_{\Ba},V_{\Ba})$, where we define
\[
V_{\Ba}=\{X\in\mathfrak{H}(\bb{C}^{2}):(X\Ba,\Ba^{\perp})_{\bb{C}^{2}}=0\}.
\]
  \item $V=\left\{\mat{0}{\bra{c}}{c}{0}:c\in\bb{C}\right\}$
\textcolor{red}{$(\ClD,\ClD')$}$\sim(W_{\Ba},V'_{\Ba})$, where we define
\[
V'_{\Ba}=\{X\in\mathfrak{H}(\bb{C}^{2}):(X\Ba,\Ba)_{\bb{C}^{2}}=(X\Ba^{\perp},\Ba^{\perp})_{\bb{C}^{2}}=0\}.
\]
In all cases above $\Ba\cdot\Ba=1$, (where $\pm\Ba$ define the same $\Pi$). I
all cases, except $(\ClD,\Be_{1}\otimes\Be_{1})$, vectors $\pm\Ba^{\perp}$
also define the same subspace as $\Ba$.
\end{itemize}
\item $\displaystyle
W=\left\{\mat{\Ga-\Gb}{i\Gb}{i\Gb}{\Ga+\Gb}:\{\Ga,\Gb\}\subset\bb{C}\right\}
$, (equivalent set coming from $\bra{W}$, $\bra{V}$)
\begin{itemize}
\item $\dim V=0$, $V=\{0\}$ \textcolor{red}{$(W,0)$}$\sim(\bra{W},0)$
\item $\dim V=1$, $V=\bb{R}\Bz_{0}\otimes\bra{\Bz_{0}}$
  \textcolor{red}{$(W,\bb{R}\BZ_{0})$}$\sim(\bra{W},\bb{R}\bra{\BZ}_{0})$
\item $\dim V=2$. There is a 1-parameter family of 2D subspaces permuted by the
  automorphism $\Phi$. This $\Phi$-orbit can be represented by
\[
V=\left\{\mat{y}{x+iy}{x-iy}{y}:\{x,y\}\subset\bb{R}\right\}
\]
We denote this set by 
\textcolor{red}{$(W,V_{\infty})$}$\sim(W,V_{t})\sim(\bra{W},\bra{V_{t}})$, 
\[
V_{t}=\left\{\mat{x}{t(x-y)+i\frac{x+y}{2}}{t(x-y)-i\frac{x+y}{2}}{y}:\{x,y\}\subset\bb{R}\right\},
\]
\item $\dim V=3$,
$
V=\left\{X\in\mathfrak{H}(\bb{C}^{2}):
\Trc X=2\im(X_{12})\right\}
$ \textcolor{red}{$(W,V)$}$\sim(\bra{W},\bra{V})$
\end{itemize}

\item $W=\Sym(\bb{C}^{2})$
  \begin{itemize}
  \item $V=\{0\}$ \textcolor{red}{$(\Sym(\bb{C}^{2}),0)$}
  \item $V=\mathfrak{H}(\bb{C}^{2})$ \textcolor{red}{$\Sym(\CT)$}
  \end{itemize}
\end{itemize}
The Summary table below also lists subalgebras, squares and ideals of each of
the algebras $\Pi(V,W)$. They have been computed with Maple computer algebra
package by Huilin Chen.

\newpage For the purposes of writing Maple code we will refer to each Jordan
multialgebra (labeled in red) by its item number in the list below. The orbit
of each equivalence class is also indicated, but will not be used in Maple
directly, unless explicitly stated.\\
\begin{tabular}[h]{|c|c|c|c|c|}
  \hline
item \# & representative & orbit & dimensions & subalgebras\\
\hline
1 & \textcolor{red}{$(0,0)$} & * & (0,0) & []\\
\hline
2 & \textcolor{red}{$(0,\bb{R}\BZ_{0})$} & $(0,\bb{R}\bra{\BZ}_{0})$ & (0,1) & [\textcolor{red}{1}]\\
\hline
3 & \textcolor{red}{$(\bb{C}\BI,0)$} & * & (1,0) & [\textcolor{blue}{1}]\\
\hline
4&\textcolor{red}{$(\bb{C}\BI,\bb{R}\BI)$}&$(\bb{C}\BI,\bb{R}\phi_{t})$, $t\in\bb{R}$&
(1,1)&[\textcolor{blue}{1},3]\\
\hline
5 &\textcolor{red}{$(\bb{C}\BI,\bb{R}\psi(i))$}&
$(\bb{C}\BI,\bb{R}\psi(e^{i\Ga}))$, $\Ga\in[0,\pi)$& (1,1)&[\textcolor{blue}{1},3]\\
\hline
6& \textcolor{red}{$(\bb{C}\BI,i\BR_{\perp})$}&$(\bb{C}\BI,\bb{R}\phi'_{t})$, $t\in\bb{R}$&
(1,1)&[\textcolor{blue}{1},3]\\
\hline
7& \textcolor{red}{$(\bb{C}\BI,\bb{R}\BZ_{0})$}&$(\bb{C}\BI,\bb{R}\bra{\BZ}_{0})$&
(1,1)&[\textcolor{blue}{1},\textcolor{blue}{2},3]\\
\hline
8& \textcolor{red}{$(\bb{C}\BI,\BGF)$}&*&(1,2)&[\textcolor{blue}{1},2,3,4,6,7]\\
\hline
9& \textcolor{red}{$(\bb{C}\BI,\BGY)$}&*&(1,2)&[\textcolor{blue}{1},3,5]\\
\hline
10&\textcolor{red}{$(\tns{\Be_{1}},0)$}&$(\tns{\Ba},0)$,  $\Ba\sim-\Ba$&
(1,0)&[\textcolor{blue}{1}]\\
\hline
11& \textcolor{red}{${\rm Ann}(\bb{C}\Be_{2})$}&${\rm Ann}(\bb{C}\bra{\Ba}^{\perp})$, 
 $\Ba\sim-\Ba$&(1,1)&[\textcolor{blue}{1},10]\\
\hline
12& \textcolor{red}{$(\tns{\Bz_{0}},0)$}&$(\tns{\bra{\Bz}_{0}},0)$&(1,0)&[\textcolor{red}{1}]\\
\hline
13& \textcolor{red}{${\rm Ann}(\bb{C}\bra{\Bz_{0}})$}&${\rm Ann}(\bb{C}\Bz_{0})$&(1,1)&
[\textcolor{red}{1},\textcolor{blue}{2},\textcolor{blue}{12}]\\
\hline
14&\textcolor{red}{$(\ClD,0)$}&$(W_{\Ba},0)$,  $\pm\Ba\sim\pm\Ba^{\perp}$&
(2,0)&[\textcolor{blue}{1},3,\textcolor{blue}{10}]\\
\hline
15&\textcolor{red}{$(\ClD,\Be_{1}\otimes\Be_{1})$}&$(W_{\Ba},\bb{R}\Ba\otimes\bra{\Ba})$,
 $\Ba\sim-\Ba$& (2,1)&[\textcolor{blue}{1},3,10,\textcolor{blue}{-10},\textcolor{blue}{11},14]\\
\hline
16&\textcolor{red}{$(\ClD,\ClD)$}&$(W_{\Ba},V_{\Ba})$, 
$\pm\Ba\sim\pm\Ba^{\perp}$& (2,2)&[\textcolor{blue}{1},3,4,-5,10,\textcolor{blue}{11},14,15]\\
\hline
17& \textcolor{red}{$(\ClD,\ClD')$}&$(W_{\Ba},V'_{\Ba})$, 
$\pm\Ba\sim\pm\Ba^{\perp}$& (2,2)&[\textcolor{blue}{1},3,5,6,10,14]\\
\hline
18&\textcolor{red}{$(W,0)$}&$(\bra{W},0)$&(2,0)&[\textcolor{blue}{1},3,\textcolor{blue}{12}]\\
\hline
19&\textcolor{red}{$(W,\bb{R}\BZ_{0})$}&$(\bra{W},\bb{R}\bra{\BZ}_{0})$&(2,1)&
[\textcolor{blue}{1},\textcolor{blue}{2},3,7,\textcolor{blue}{12},\textcolor{blue}{13},18]\\
\hline
20&\textcolor{red}{$(W,V_{\infty})$}&$(W,V_{t})\sim(\bra{W},\bra{V_{t}})$&(2,2)&
[\textcolor{blue}{1},2,3,5,7,12,\textcolor{blue}{13},18,19]\\
\hline
21&\textcolor{red}{$(W,V)$}&$(\bra{W},\bra{V})$&(2,3)&[\textcolor{blue}{1},2,3,5,7,9,12,\textcolor{blue}{13},18,19,20]\\
\hline
22&\textcolor{red}{$(\Sym(\bb{C}^{2}),0)$}&*&(3,0)&[\textcolor{blue}{1},3,10,12,14,18]\\
\hline
23&\textcolor{red}{$\Sym(\CT)$}&*&(3,4)&$[\textcolor{blue}{1},2,\ldots,21,22]$\\
\hline
\end{tabular}
\begin{itemize}
\item Symbol * in the ``orbit'' column means that the orbit consists of a single
  algebra listed in the ``representative'' column.
\item Vectors $\Ba$ always lie on
  the ``complex circle'' 
$
\bb{S}_{\bb{C}}^{1}=\{\Ba\in\bb{C}^{2}:\Ba\cdot\Ba=1\}.
$
\item Algebra -10 in item 15, refers to an algebra from the orbit of item 10, corresponding
  to $\Ba=\Be_{2}$: $(\tns{\Be_{2}},0)$. It is there because among all global
  automorphisms mapping item 15 into itself none map 10 into -10. Therefore,
  within item 15 algebras 10 and -10 are not equivalent. There are no other occurrances
  of such a situation.
\item Algebra -5 in item 16 refers to an algebra from the orbit of item 5, corresponding
  to $\Ga=0$: $(\bb{C}\BI,\bb{R}\psi(1))$.
\item If an algebra is not listed as a subalgebra of a particular algebra it
  means that no algebra from its orbit is a subalgebra of that particular algebra.
\item The subalgebras listed in red are squares, the subalgebras listed in
  blue are ideals.
\end{itemize}
\newpage

\section{Theory of links}
Another related and important part of the project is to discover all possible
links. In order to describe a link, consider the opposite exercise. Instead of
fixing tensors $\SFL_{A}$ and $\SFL_{B}$ we are fixing the set $A$ and varying
$\SFL_{A}$ and $\SFL_{B}$. If we know $\SFL^{*}$ for one pair $\SFL_{A}$,
$\SFL_{B}$, does it give us any information about $\SFL^{*}$ for other pairs?
If the answer is yes for any subset $A$, then we say that we have discovered a
link. Links are much harder to characterise. Since, in general they contain a
lot more information than exact relations.

In the framework of the theory, links are described by Jordan
$\Hat{\CA}$-multialgebras in $\Sym(\CT)\oplus\Sym(\CT)$, where 
\[
\Hat{\CA}=\Span\left\{
\mat{\BGG_{0}^{(1)}(\Bn)-\BGG_{0}^{(1)}(\Bn_{0})}{0}{0}{\BGG_{0}^{(2)}(\Bn)-\BGG_{0}^{(2)}(\Bn_{0})}:
|\Bn|=1\right\},
\]
where $\BGG_{0}^{(1)}(\Bn)$ and $\BGG_{0}^{(2)}(\Bn)$ are constructed using
different reference media $\SFL_{0}^{(1)}$ and $\SFL_{0}^{(2)}$, respectively.
In our case
\[
\Hat{\CA}=\Span\left\{
\mat{\BGL_{1}^{-1}\otimes\BA}{0}{0}{\BGL_{2}^{-1}\otimes\BA}:\BA^{T}=\BA,\Trc\BA=0\right\}.
\]
We say that $\Hat\Pi\subset\Sym(\CT)\oplus\Sym(\CT)$ describes a link if
\[
\mat{\SFK_{1}}{0}{0}{\SFK_{2}}\mat{\SFA_{1}}{0}{0}{\SFA_{2}}\mat{\SFK_{1}}{0}{0}{\SFK_{2}}
\in\Pi,\quad\forall\mat{\SFK_{1}}{0}{0}{\SFK_{2}}\in\Hat{\Pi},\ 
\forall\mat{\SFA_{1}}{0}{0}{\SFA_{2}}\in\Hat{\CA}.
\]
From now on we will be using a more compact notation 
\[
[\SFK_{1},\SFK_{2}]=\mat{\SFK_{1}}{0}{0}{\SFK_{2}},\qquad
[\SFA_{1},\SFA_{2}]=\mat{\SFA_{1}}{0}{0}{\SFA_{2}}.
\]

As in the case of Jordan $\CA$-multialgebras we will first apply the
convariance transformation 
\[
\Hat{\CA}_{0}=\Hat{\SFC}\Hat{\CA}\Hat{\SFC}^{T}
\]
where 
$\Hat{\SFC}=[\BC_{1}\otimes\BI_{d},\BC_{2}\otimes\BI_{d}]$, and $\BC_{1}$,
$\BC_{2}$ are as before: $\BC_{1}=\BGL_{1}^{1/2}$, $\BC_{2}=\BGL_{2}^{1/2}$, so that 
\[
\Hat{\CA}_{0}=\{[\BA,\BA]:\BA\in\CA_{0}\},\qquad\CA_{0}=\{\BI_{2}\otimes\BA:\BA^{T}=\BA,\Trc\BA=0\}.
\]
All Jordan $\Hat{\CA}$-multialgebras can be described entirely in terms of the
algebraic structure of Jordan $\CA$-multialgebras. In order to describe an
$\Hat{\CA}$-multialgebra $\Hat{\Pi}$ we need the following algebraic data:
Jordan $\CA$-ideals $\CI_{1}\subset\Pi_{1}$, $\CI_{2}\subset\Pi_{2}$, such
that the factor-algebras $\Pi_{1}/\CI_{1}$ and $\Pi_{2}/\CI_{2}$ are
isomorphic, since we will also require a Jordan $\CA$-factoralgebra
isomorphism $\Phi:\Pi_{1}/\CI_{1}\to\Pi_{2}/\CI_{2}$. In that case 
\begin{equation}
  \label{Pihatstr}
  \Hat{\Pi}=\{[\SFK_{1},\SFK_{2}]\in\Pi_{1}\times\Pi_{2}:\Phi([\SFK_{1}])=[\SFK_{2}],\},
\end{equation}
where $[\SFK_{j}]$ denotes the equivalence class of $\SFK_{j}$ in
$\Pi_{j}/\CI_{j}$, $j=1,2$.

The most common occurrence is the situation, where $\Pi_{1}=\Pi_{2}=\Pi$ and
$\CI_{1}=\CI_{2}=\{0\}$, in which case
\begin{equation}
  \label{Philink}
  \Hat{\Pi}=\{[\SFK,\Phi(\SFK)]:\SFK\in\Pi\}.
\end{equation}
Another common occurrence happens when there exists a Jordan
$\CA$-multialgebra $\Pi'\subset\Pi$, such that $\Pi=\CI\oplus\Pi'$, where
$\CI$ is an ideal in $\Pi$. That means that every $\SFK\in\Pi$ can be written
uniquely as $\SFK=\SFK'+\SFJ$, where $\SFK'\in\Pi'$ and $\SFJ\in\CI$. The map 
$\Phi([\SFK])=\SFK'$ is obviously a factor-algebra isomorphism
$\Phi:\Pi/\CI\to\Pi'/\{0\}$. In that case
\begin{equation}
  \label{idlink}
  \Hat{\Pi}=\{[\SFK'+\SFJ,\SFK']:\SFK'\in\Pi',\ \SFJ\in\CI\}.
\end{equation}
We will see later that in our case only links of the above two types are present.

\section{Factor algebra isomorphism classes}
A Maple ideal checker has found 11 nontrivial ideals. In each case we have a situation where
$\Pi=\Pi'\oplus J$, where $\Pi'$ is a subalgebra and $J$ is an ideal. 
In that
case, every $\SFK\in\Pi$ can be written uniquely as $\SFK=\SFK'+\SFJ$ and
hence $\SFK'$ becomes a natural choice of the representative of the
equivalence class of $\SFK$ in $\Pi/J$. This identification is obviously an
algebra isomorphism.
Then
natural projection $\pi:\Pi\to\Pi/J\cong\Pi'$ is an algebra isomorphism:
\[
[\SFK*_{\SFA}\SFK]=[\SFK'*_{\SFA}\SFK'+2\SFK'*_{\SFA}\SFJ+\SFJ*_{\SFA}\SFJ]=\SFK'*_{\SFA}\SFK'.
\]
\begin{enumerate}
\item $(0,\bb{R}\BZ_{0})\subset(\bb{C}\BI,\bb{R}\BZ_{0})$. In this case
  $\Pi=\Pi'\oplus J$, where $\Pi'=(\bb{C}\BI,0)$, and therefore, the factor
  algebra is naturally isomorphic to $\Pi'$.
\item $(\tns{\Be_{1}},0)\subset(\ClD,0)$.  In this case
  $\Pi=\Pi'\oplus J$, where $\Pi'=(\tns{\Be_{2}},0)$, and therefore, the factor
  algebra is naturally isomorphic to $\Pi'$.
\item $(\tns{\Be_{2}},0)\subset(\ClD,\Be_{1}\otimes\Be_{1})$. In this case
  $\Pi=\Pi'\oplus J$, where $\Pi'={\rm Ann}(\bb{C}\Be_{2})$, and therefore, the factor
  algebra is naturally isomorphic to $\Pi'$.
\item ${\rm Ann}(\bb{C}\Be_{2})\subset(\ClD,\Be_{1}\otimes\Be_{1})$. In this case
  $\Pi=\Pi'\oplus J$, where $\Pi'=(\tns{\Be_{2}},0)$, and therefore, the factor
  algebra is naturally isomorphic to $\Pi'$.
\item ${\rm Ann}(\bb{C}\Be_{2})\subset(\ClD,\ClD)$. In this case
  $\Pi=\Pi'\oplus J$, where $\Pi'={\rm Ann}(\bb{C}\Be_{1})$, and therefore, the factor
  algebra is naturally isomorphic to $\Pi'$. It remains to recall that the
  algebras ${\rm Ann}(\bb{C}\Be_{1})$ and ${\rm Ann}(\bb{C}\Be_{2})$ are isomorphic my
  means of the global isomorphism. Thus, $(\ClD,\ClD)/{\rm
    Ann}(\bb{C}\Be_{2})\cong{\rm Ann}(\bb{C}\Be_{2})$. 
\item $(\tns{\Bz_{0}},0)\subset(W,0)$. In this case
  $\Pi=\Pi'\oplus J$, where $\Pi'=(\bb{C}\BI,0)$, and therefore, the factor
  algebra is naturally isomorphic to $\Pi'$.
\item $(0,\bb{R}\BZ_{0})\subset(W,\bb{R}\BZ_{0})$.  In this case
  $\Pi=\Pi'\oplus J$, where $\Pi'=(W,0)$, and therefore, the factor
  algebra is naturally isomorphic to $\Pi'$.
\item $(\tns{\Bz_{0}},0)\subset(W,\bb{R}\BZ_{0})$.  In this case
  $\Pi=\Pi'\oplus J$, where $\Pi'=(\bb{C}\BI,\bb{R}\BZ_{0})$, and therefore, the factor
  algebra is naturally isomorphic to $\Pi'$.
\item ${\rm Ann}(\bb{C}\bra{\Bz_{0}})\subset(W,\bb{R}\BZ_{0})$.  In this case
  $\Pi=\Pi'\oplus J$, where $\Pi'=(\bb{C}\BI,0)$, and therefore, the factor
  algebra is naturally isomorphic to $\Pi'$.
\item ${\rm Ann}(\bb{C}\bra{\Bz_{0}})\subset(W,V_{\infty})$. In this case
  $\Pi=\Pi'\oplus J$, where $\Pi'=(\bb{C}\BI,\bb{R}\psi(i))$, and therefore, the factor
  algebra is naturally isomorphic to $\Pi'$.
\item ${\rm Ann}(\bb{C}\bra{\Bz_{0}})\subset(W,V)$. In this case
  $\Pi=\Pi'\oplus J$, where $\Pi'=(\bb{C}\BI,\BGY)$, and therefore, the factor
  algebra is naturally isomorphic to $\Pi'$.
\end{enumerate}
Since there are no new factor-algebras in addition to the 23 algebras above,
we only need to know the algebra-ideal pairs. This information can kept in a
more economical list, observing that if $J\subset\Pi$ is an ideal, then for
any other algebra $\Pi'$ $J\cap\Pi'$ is an ideal in $\Pi\cap\Pi'$.
In this way, we have a reduced set of algebra-ideal pairs.
\begin{enumerate}
\item $(\ClD,\Be_{1}\otimes\Be_{1})/(\tns{\Be_{2}},0)\cong{\rm Ann}(\bb{C}\Be_{2})$
\item $(\ClD,\ClD)/{\rm Ann}(\bb{C}\Be_{2})\cong{\rm Ann}(\bb{C}\Be_{2})$
\item $(W,\bb{R}\BZ_{0})/(0,\bb{R}\BZ_{0})\cong(W,0)$
\item $(W,\bb{R}\BZ_{0})/(\tns{\Bz_{0}},0)\cong(\bb{C}\BI,\bb{R}\BZ_{0})$
\item $(W,V)/{\rm Ann}(\bb{C}\bra{\Bz_{0}})\cong(\bb{C}\BI,\BGY)$
\end{enumerate}
We remark that in the list above algebra/ideal pairs 1 and 2 represent
links in the absense of
thermoelectric coupling. Item 1 corresponds to the KDM link for 2D
conductivity, while item 2 just corresponds to a pair of uncoupled conducting
composites and it says that the effective tensor for the pair is a pair of
effective tensors of each of the composites.

\section{SO(2)-invariant Jordan multialgebra automorphisms}
We can partially determine all possible automorphisms $\Phi$ of each of the
algebras $\Pi$ by describing all transformations $\Phi_{2}:W\to W$ satisfying
(\ref{PhiXX})$_{2}$. There are only 7 possibilities for $W$
\begin{enumerate}
\item $W=\{0\}$. Then the only choice is the ``identity map'' $\Phi_{2}(Y)=Y$.
\item $W=\bb{C}\BI_{2}$. Then the only choice is the ``identity map''
  $\Phi_{2}(Y)=Y$.
\item $W=\bb{C}\tns{\Be_{1}}$. Then the only choice is the ``identity map''
  $\Phi_{2}(Y)=Y$.
\item $W=\bb{C}\tns{\Bz_{0}}$. In that case every nonzero linear map satisfies
  (\ref{PhiXX})$_{2}$: $\Phi(Y)=aY$ for some $a\in\bb{C}\setminus\{0\}$.
\item $W=\ClD$. Then in addition to the ``identity map'' $\Phi_{2}(Y)=Y$ there
  is one more possibility:
\[
\Phi_{2}\left(\mat{x}{0}{0}{y}\right)=\mat{y}{0}{0}{x}.
\]
\item $W=\Span_{\bb{C}}\{\BI_{2},\tns{\Bz_{0}}\}$. In that case $\Phi_{2}$ is
determined by its values on basis vectors:
\[
\Phi_{2}(\BI_{2})=\BI_{2},\qquad\Phi_{2}(\tns{\Bz_{0}})=a\tns{\Bz_{0}},\qquad
a\in\bb{C}\setminus\{0\}.
\]
\item $W=\Sym(\bb{C}^{2})$. This case has already been examined (in Section
  10). The set of all maps $\Phi_{2}$ is described by
\[
\Phi_{2}(Y)=CYC^{T},\qquad C\in O(2,\bb{C}).
\]
\end{enumerate}
The problem of determination of $\Phi_{0}$ is trivial when $V=\{0\}$, which is
true in 7 cases. It has also been solved by $\Pi=\Sym(\CT)$. In another 9
cases $\dim V=1$. Then in order to determine $\Phi_{0}$ we only need to find
all real non-zero numbers $\Ga$ for which $\Phi_{0}(X)=\Ga X$. If $X_{0}\in
V\setminus\{0\}$, then equations (\ref{PhiXX}), (\ref{PhiXY}) imply
\[
\Phi_{2}(Y)X_{0}=YX_{0},\qquad\Phi_{2}(X_{0}X_{0}^{T})=\Ga^{2}X_{0}X_{0}^{T}.
\]
In particular, if $X_{0}X_{0}^{T}X_{0}\not=0$, then $\Ga=\pm 1$ are the only
choices that can work. If $X_{0}=\BZ_{0}$, then any $\Ga\not=0$
works. Finally, we need to note that in the case
$\Pi=(\ClD,\Be_{1}\otimes\Be_{1})$ the nontrivial map $\Phi_{2}$ is ruled out,
since
\[
\mat{y}{0}{0}{x}\Be_{1}\otimes\Be_{1}\not=\mat{x}{0}{0}{y}\Be_{1}\otimes\Be_{1}.
\]
Thus, for this algebra, the only nontrivial automorphism is defined by
$\Phi_{0}(X)=-X$ and $\Phi_{2}(Y)=Y$. There are only 6 cases (except
$\Pi=\Sym(\CT)$), where $\dim V>1$. In 2 of these 6 cases $W=\bb{C}\BI_{2}$
and therefore $\Phi_{2}(Y)=Y$, while $\Phi_{0}$ satisfies
\[
\Phi_{0}(X)\Phi_{0}(X)^{T}=XX^{T},\qquad X\in V.
\]

\newpage
\begin{center}
  List of all SO(2)-invariant Jordan multialgebra automorphisms
\end{center}
\begin{tabular}[h]{|c|c|c|}
  \hline
item \# & representative & automorphisms\\
\hline
1 & \textcolor{red}{$(0,0)$} & *\\
\hline
2 & \textcolor{red}{$(0,\bb{R}\BZ_{0})$} & $\Phi_{0}(\BZ_{0})=\Ga\BZ_{0}$\\
\hline
3 & \textcolor{red}{$(\bb{C}\BI,0)$} & *\\
\hline
4&\textcolor{red}{$(\bb{C}\BI,\bb{R}\BI)$}&$\Phi_{0}(\BI)=-\BI$\\
\hline
5 &\textcolor{red}{$(\bb{C}\BI,\bb{R}\psi(i))$}&$\Phi_{0}(\psi(i))=-\psi(i)$\\
\hline
6& \textcolor{red}{$(\bb{C}\BI,i\BR_{\perp})$}&$\Phi_{0}(i\BR_{\perp})=-i\BR_{\perp}$\\
\hline
7& \textcolor{red}{$(\bb{C}\BI,\bb{R}\BZ_{0})$}&$\Phi_{0}(\BZ_{0})=\Ga\BZ_{0}$\\
\hline
8& \textcolor{red}{$(\bb{C}\BI,\BGF)$}&see below\\
\hline
9& \textcolor{red}{$(\bb{C}\BI,\BGY)$}&see below\\
\hline
10&\textcolor{red}{$(\tns{\Be_{1}},0)$}&*\\
\hline
11& \textcolor{red}{${\rm Ann}(\bb{C}\Be_{2})$}&$\Phi_{0}(\tns{\Be_{1}})=-\tns{\Be_{1}}$\\
\hline
12& \textcolor{red}{$(\tns{\Bz_{0}},0)$}&$\Phi(\SFK)=a\SFK$\\
\hline
13& \textcolor{red}{${\rm Ann}(\bb{C}\bra{\Bz_{0}})$}&$\Phi_{2}(\tns{\Bz_{0}})=a\tns{\Bz_{0}}$, 
$\Phi_{0}(\BZ_{0})=\Ga\BZ_{0}$\\
\hline
14&\textcolor{red}{$(\ClD,0)$}&$\Phi_{2}(Y)=\psi(i)Y\psi(i)$\\
\hline
15&\textcolor{red}{$(\ClD,\Be_{1}\otimes\Be_{1})$}&$\Phi_{0}(X)=-X$\\
\hline
16&\textcolor{red}{$(\ClD,\ClD)$}&$\Phi_{0}(X)=\pm\psi(i)X\psi(i)$, $\Phi_{2}(Y)=\psi(i)Y\psi(i)$\\
\hline
17& \textcolor{red}{$(\ClD,\ClD')$}&$\Phi_{0}(X)=-X$, $\Phi_{2}(Y)=Y$ or 
$\Phi_{0}(X)=\pm X^{T}$, $\Phi_{2}(Y)=\psi(i)Y\psi(i)$\\
\hline
18&\textcolor{red}{$(W,0)$}&$\Phi_{2}(\BI_{2})=\BI_{2}$, $\Phi_{2}(\tns{\Bz_{0}})=a\tns{\Bz_{0}}$\\
\hline
19&\textcolor{red}{$(W,\bb{R}\BZ_{0})$}&$\Phi_{2}(\BI_{2})=\BI_{2}$,
$\Phi_{2}(\tns{\Bz_{0}})=a\tns{\Bz_{0}}$, $\Phi_{0}(\BZ_{0})=\Ga\BZ_{0}$\\
\hline
20&\textcolor{red}{$(W,V_{\infty})$}&see below\\
\hline
21&\textcolor{red}{$(W,V)$}&see below\\
\hline
22&\textcolor{red}{$(\Sym(\bb{C}^{2}),0)$}&$\Phi_{2}(Y)=CYC^{T}$, $C\in O(2,\bb{C})$\\
\hline
23&\textcolor{red}{$\Sym(\CT)$}&$\Phi(K(X,Y))=K(\pm CXC^{H},CYC^{T})$, $C\in O(2,\bb{C})$\\
\hline
\end{tabular}

Item 8: $\Phi_{2}(\BI_{2})=\BI_{2}$ and $\Phi_{0}(\BGF(x,y))=\BGF(x',y')$, where
\[
\BGF(x,y)=\mat{x}{iy}{-iy}{x},\qquad\vect{x'}{y'}=\BF\vect{x}{y},\qquad\BF^{T}\psi(1)\BF=\BI_{2}.
\]
Then,
\[
\BF=\pm\mat{\cosh t}{\sinh t}{\sinh t}{\cosh t}\text{ or }
\BF=\pm\mat{\cosh t}{\sinh t}{-\sinh t}{-\cosh t}.
\]
Item 9: $\Phi_{2}(\BI_{2})=\BI_{2}$ and $\Phi_{0}(\BGY(x,y))=\BGY(x',y')$, where
\[
\BGY(x,y)=\mat{x}{y}{y}{-x},\qquad\vect{x'}{y'}=\BF\vect{x}{y},\qquad\BF^{T}\BF=\BI_{2}.
\]
Then,
\[
\BF=\mat{\cos t}{\sin t}{-\sin t}{\cos t}\text{ or }
\BF=\mat{\cos t}{\sin t}{\sin t}{-\cos t}.
\]
If we denote $z=x+iy$, then $\BGY(x,y)=\psi(z)$ and
$\Phi_{0}(\psi(z))=\psi(e^{-it}z)$ or $\psi(e^{it}\bra{z})$.

Item 20: Every $X\in V_{\infty}$ has the general form
$X=\xi\psi(i)+\eta\BZ_{0}$, $\{\xi,\eta\}\subset\bb{R}$, while every $Y\in W$
has the general form $Y=x\BI_{2}+y\tns{\Bz_{0}}$, $\{x,y\}\subset\bb{C}$. Then
\[
\Phi_{0}(\xi\psi(i)+\eta\BZ_{0})=\pm(\xi\psi(i)+\Ga\eta\BZ_{0}),\qquad
\Phi_{2}(x\BI_{2}+y\tns{\Bz_{0}})=x\BI_{2}+\Ga y\tns{\Bz_{0}},\qquad\Ga\in\bb{R}\setminus\{0\}.
\]
Item 21: Every $X\in V$ has the general form
$X=\psi(z)+\eta\BZ_{0}$, $z\in\bb{C}$, $\eta\in\bb{R}$, while every $Y\in W$
has the general form $Y=x\BI_{2}+y\tns{\Bz_{0}}$, $\{x,y\}\subset\bb{C}$. Then
\[
\Phi_{0}(\psi(z)+\eta\BZ_{0})=\pm(\psi(e^{i\Gth}z)+\rho\eta\BZ_{0}),\quad
\Phi_{2}(x\BI_{2}+y\tns{\Bz_{0}})=x\BI_{2}+\rho e^{i\Gth}y\tns{\Bz_{0}},\quad \rho e^{i\Gth}\in\bb{C}\setminus\{0\}.
\]

Our next task is to determine which of these automorphisms are not
restrictions of the global one to the multialgebra in question.
For this purpose we define
\begin{equation}
  \label{Cpmdef}
  \BC_{+}(c)=\mat{\cos c}{\sin c}{-\sin c}{\cos c},\qquad
\BC_{-}(c)=\mat{\cos c}{\sin c}{\sin c}{-\cos c},\qquad c\in\bb{C}.
\end{equation}
We compute
\[
\BC_{+}(c)\Bz_{0}=e^{-ic}\Bz_{0},\qquad\BC_{-}(c)\Bz_{0}=e^{-ic}\bra{\Bz}_{0}
\]
Hence,
\[
\BC_{+}(c)\BZ_{0}\BC_{+}(c)^{H}=e^{2\im(c)}\BZ_{0},\qquad\BC_{-}(c)\BZ_{0}\BC_{-}(c)^{H}=e^{2\im(c)}\bra{\BZ}_{0}.
\]
\[
\BC_{+}(c)\psi(z)\BC_{+}(c)^{H}=\psi(e^{-2i\re(c)}z),\qquad\BC_{-}(c)\psi(z)\BC_{-}(c)^{H}=\psi(e^{2i\re(c)}\bra{z}),\quad z\in\bb{C}
\]
These formulas show that the Automorphisms of algebras 20 and 21, as well as 9 are all
restrictions of the global automorphism.

In Item 19 there are new automorphisms. We can use the global one to set
$a=1$. The remaining ones $\Phi_{2}(Y)=Y$, $\Phi_{0}(\BZ_{0})=\Ga\BZ_{0}$, are
all new (except when $\Ga=\pm 1$). The same remark hold for item 13. However,
all automorphisms are now restrictions of automorphisms of algebra 19.

In order to decide on item 8. We similarly denote $\BGF(x,y)$ by $\BGF(z)$,
where $z=x+iy$. In that case the transformation $\Phi_{0}$ acts by
\[
\Phi_{0}^{+}(\BGF(z))=\pm\BGF(z\cosh t+i\bra{z}\sinh t),\text{ or }
\Phi_{0}^{-}(\BGF(z))=\pm\BGF(\bra{z}\cosh t-iz\sinh t).
\]
It remains to compute (via Maple) that
\[
\BC_{+}(c)\BGF(z)\BC_{+}(c)^{H}=\Phi_{0}^{+}(\BGF(z)),\qquad
\BC_{-}(c)\BGF(z)\BC_{-}(c)^{H}=\Phi_{0}^{-}(\BGF(z)),
\]
where $t=2\im(c)$. Hence, all automorphisms of algebra 8 are generated by the
global ones.

The automorphisms of the remaining algebras are easily seen to come from the
global ones. This leaves a single family of non-global automorphisms for
algebra 19:
\[
\Phi_{2}(Y)=Y,\qquad\Phi_{0}(\BZ_{0})=\Ga\BZ_{0},\quad Y\in W,\ \Ga\in\bb{R}\setminus\{0\}.
\]

\section{Eliminating redundancies}
We can now eliminate some of the Jordan multialgebras from our list either
because they are physically trivial or because they can be obtained as
intersections of other multialgebras.
\begin{enumerate}
\item %1
is physically trivial
\item %2
$2=7\cap 13$ and $2^{2}=7^{2}\cap 13^{2}$; 
$\mat{\Gl\BI_{2}}{\pm(\Gl-1)\BR_{\perp}}{\mp(\Gl-1)\BR_{\perp}}{\Gl\BI_{2}}$, $\Gl>\hf$,
$(\Gl^{*})^{-1}=\av{\Gl^{-1}}$.
\item %3
$3=4\cap 7$
\item %4
$4=8\cap 16$;  $\mat{\BGs}{0}{0}{\BGs}$, $\BGs>0$.
\item %5
$5=9\cap 17$; $\mat{\BL}{t\BL}{t\BL}{\BL}$, $\det\BL=\nth{1-t^{2}}$, $|t|<1$, $\BL>0$;\\[2ex] 
$\BGs=(t+1)\BL$, $t^{*}=\frac{\det\BGs^{*}-1}{\det\BGs^{*}+1}$, 
$\BL^{*}=\dfrac{\det\BGs^{*}+1}{2}\cdot\dfrac{\BGs^{*}}{\det\BGs^{*}}$.\\[2ex]
$5'=(\bb{C}\BI_{2},\bb{R}\psi(1))=9\cap 16$, $\mat{\BGs}{0}{0}{\frac{\BGs}{\det\BGs}}$, $\BGs>0$.
\item %6
$6=8\cap 17$; $\mat{\BL}{-t\BR_{\perp}}{t\BR_{\perp}}{\BL}$, $\det\BL=1+t^{2}$, $t\in\bb{R}$, $\BL>0$.
\item %7
$7=8\cap 19$; (see (\ref{ER7}))
\item %8
is essential (see (\ref{ER8}))
\item %9 
$9=21\cap\bra{21}$ (see (\ref{ER9}))
\item %10
$10=11\cap 14$
\item %11
is physically trivial because it is an ER in the absense of thermoelectric coupling
\item %12
$12=13\cap 18$ and $12^{2}=13^{2}\cap 18^{2}$ (see
(\ref{Annz0}) and (\ref{FR}))
\item %13
is essential because of the volume fraction relation that accompanies it
\item %14
is physicaly trivial because it is an ER in the absense of thermoelectric coupling
\item %15
is physically trivial because it is an ER in the absense of thermoelectric coupling
\item %16
is physically trivial because it is an ER in the absense of thermoelectric coupling
\item %17
is essential (see (\ref{ER17}) or (\ref{ER17fin}))
\item %18
$18=19\cap 22$ (see (\ref{ER19fin}) and (\ref{ER18ad}))
\item %19
$19=20\cap 20_{t}$: $(W,\bb{R}\BZ_{0})=(W,V_{\infty})\cap(W,V_{t})$ for any $t\in\bb{R}$.
(see (\ref{ER19fin}))
\item %20 
is essential (see (\ref{ER20})--(\ref{ER20fin}))
\item %21
is essential (see (\ref{ER21}) or (\ref{ER21fin}))
\item %22
is essential (see (\ref{ER22}) or (\ref{ER22fin}))
\item %23
is physically trivial
\end{enumerate}

\section{Verifying 3 and 4-chain properties}
Recall that every exact relation corresponds to a Jordan multialgebra. However,
theoretically, not every Jordan multialgebra may correspond to an exact
relation. Validity of 3 and 4-chain properties for a Jordan multialgebra ensures that it 
corresponds to an exact relation. Specifically we need to verify
\begin{align}
  \label{3chain}
  &\SFK_1\SFA_1\SFK_2\SFA_2\SFK_3+\SFK_3\SFA_2\SFK_2\SFA_1\SFK_1\in\Pi,\\
&\SFK_1\SFA_1\SFK_2\SFA_2\SFK_3\SFA_3\SFK_4+\SFK_4\SFA_3\SFK_3\SFA_2\SFK_2\SFA_1\SFK_1\in\Pi
\label{4chain}
\end{align}
for every $\SFK_{j}\in\Pi$ and every $\SFA_{j}\in\CA$.  The algebraic meaning
of 3 and 4 chain properties is the existence of an associative
$\CA$-multialgebra $\Pi'$ (closed under the associative set of multiplications
$\SFK_{1}\circ_{\SFA}\SFK_{2}=\SFK_{1}\SFA\SFK_{2}$), such that $\SFK\in\Pi'$
implies $\SFK^{T}\in\Pi'$ and $\Pi'\cap\Sym(\CT)=\Pi$. If $\Pi$ is
rotationally invariant, then $\Pi$ must necessarily be rotationally invariant,
as well. In our setting an $SO(2)$-invariant associative $\CA$-multialgebra
$\Pi'$ is characterized by a real subspace $V'$ and a complex subspace $W'$ of
$2\times 2$ complex matrices, such that
\[
X\bra{Y},\ YX\in V',\quad X_{1}\bra{X}_{2},\ Y_{1}Y_{2}\in W'
\]
for all $X,X_{1},X_{2}\in V'$ and all $Y,Y_{1},Y_{2}\in W'$. If
\[
X^{H}\in V',\quad Y^{T}\in W',\quad\forall X\in V',\ Y\in W'
\]
and
\[
V=V'\cap\mathfrak{H}(\bb{C}^{2}),\qquad W=W'\cap\Sym(\bb{C}^{2}).
\]
Then $\Pi(V,W)$ satisfies the 3 and 4-chain properties.
For example for the algebra \#19 $(W,\bb{R}\BZ_{0})$ we can set $W'=W$ and $V'=\bb{C}\BZ_{0}$ and
verify that all the relations above hold. 

There is also a version of 3 and 4-chain properties for ideals and
automorphisms. We say that an ideal $\CI\subset\Pi$ satisifies the 3 and 4-chain
properties if 
\begin{align}
  \label{id3chain}
  &\SFJ\SFA_1\SFK_2\SFA_2\SFK_3+\SFK_3\SFA_2\SFK_2\SFA_1\SFJ\in\CI,\\
&\SFK_1\SFA_1\SFK_2\SFA_2\SFK_3\SFA_3\SFK_4+\SFK_4\SFA_3\SFK_3\SFA_2\SFK_2\SFA_1\SFK_1\in\CI
\label{id4chain}
\end{align}
for every $\SFK_{j}\in\Pi$, every $\SFA_{j}\in\CA$ and every $\SFJ\in\CI$.
Equivalently, if we happen to know the associative $\CA$-multialgebra $\Pi'$
that establishes the 3 and 4-chain properties of $\Pi$, we can look for an
associative ideal $\CI'\subset\Pi'$, such that $\CI'\cap\Sym(\CT)=\CI$.
For example, two of the 3 nontrivial ideals belong to the algebra \#19. It it
is easy to verify that $J'=(\bb{C}\tns{\Bz_{0}},0)$ and $J'=(0,\bb{C}\BZ_{0})$ are
ideals in $(V',W')$, establishing the 
3 and 4-chain properties for $J=(\bb{R}\tns{\Bz_{0}},0)$ and $J=(0,\bb{C}\BZ_{0})$.

The 3 and 4-chain properties for
automorphisms are
\[
\Phi(\SFK_{0}\SFA\SFK_{1}\SFA'\SFK_{2}+\SFK_{2}\SFA'\SFK_{1}\SFA\SFK_{0})=
\Phi(\SFK_{0})\SFA\Phi(\SFK_{1})\SFA'\Phi(\SFK_{2})+
\Phi(\SFK_{2})\SFA'\Phi(\SFK_{1})\SFA\Phi(\SFK_{0})
\] 
\begin{multline*}
\Phi(\SFK_{0}\SFA\SFK_{1}\SFA'\SFK_{2}\SFA''\SFK_{3}+
\SFK_{3}\SFA''\SFK_{2}\SFA'\SFK_{1}\SFA\SFK_{0})=\\
\Phi(\SFK_{0})\SFA\Phi(\SFK_{1})\SFA'\Phi(\SFK_{2})
\SFA''\Phi(\SFK_{3})+\Phi(\SFK_{3})\SFA''\Phi(\SFK_{2})
\SFA'\Phi(\SFK_{1})\SFA\Phi(\SFK_{0})
\end{multline*}
for all
$
\{\SFK_{0},\SFK_{1},\SFK_{2},\SFK_{3}\}\subset\Pi$, and all $
\{\SFA,\SFA',\SFA'',\SFA''\}\subset\CA.
$
In our setting the automorphisms $\Phi$ has the 3 and 4-chain properties, \IFF
it is a restriction to $(V,W)$ of the automorphism $\Phi'$ of $\Pi'$,
generated by real and complex linear maps $\Phi'_{0}$ and
$\Phi'_{2}$ on $V'$ and $W'$, respectively, satisfying
\[
\Phi'_{0}(X\bra{Y})=\Phi'_{0}(X)\bra{\Phi'_{2}(Y)},\qquad
\Phi'_{0}(YX)=\Phi'_{2}(Y)\Phi'_{0}(X),
\]
\[
\Phi'_{2}(X_{1}\bra{X}_{2})=\Phi'_{0}(X_{1})\bra{\Phi'_{0}(X_{2})},\qquad
\Phi'_{2}(Y_{1}Y_{2})=\Phi'_{2}(Y_{1})\Phi'_{2}(Y_{2}).
\]
It is easy to verify that the map defined by $\Phi'_{2}(Y)=Y$ and $\Phi'_{0}(X)=\Ga X$
satisfies all the relations above. Hence, the ER, corresponding to algebra
\#19 $(W,\bb{R}\BZ_{0})$
and the links corresponding to this family of automorphisms, as well as the
links corresponding to the two ideals in this algebra hold for all
thermoelectric composites.

Finally, there is also a version of 3 and 4-chain properties for volume
fraction relations corresponding to situations where $\Pi^{2}\not=\Pi$. In
this case we need to verify conditions (\ref{3chain}) and (\ref{4chain}),
except the chains must belong to $\Pi^{2}$, instead of $\Pi$.

Thus, we only need to check the 3 and 4-chain relations for the remaining 7
essential Jordan multialgebras (8,9,13,17,20,21,22), for the volume
fraction relation that accompanies algebra \#13, as well as for the
remaining algebra/ideal pair $(W,V)/{\rm Ann}(\bb{C}\bra{\Bz_{0}})$. All these
checks have been done with Maple by Huilin Chen and confirmed that the 3 and 4-chain
relations were satisfied in all cases.

\section{Computing inversion keys}
We have already mentioned the algorithm for computing the inversion keys for
exact relations. Let us restate it in the $K(X,Y)$-language. The inversion key
$\SFM_{0}$ is always sought is the form $\SFM_{0}=K(M_{0},0)$, where $M_{0}$
is one of the 4 choices: 0, $\tns{\Be_{1}}/2$, $\tns{\Be_{2}}/2$ or
$\BI_{2}/2$. It is found from the rule that
\[
K(X,Y)K\left(\hf\BI_{2}-M_{0},0\right)K(X,Y)\in\Pi,\quad\forall K(X,Y)\in\Pi.
\]
This property holds trivially for the choice $M_{0}=\BI_{2}/2$. Thus, only if
we want to use one of the three remaining choices there is something to
verify. Obviously, $M_{0}=0$ is the most desirable choice. We can use it only
if the Jordan multialgebra $\Pi$ satisifes
\begin{equation}
  \label{Meq0}
  K(X,Y)^{2}\in\Pi,\quad\forall K(X,Y)\in\Pi.
\end{equation}
If (\ref{Meq0}) fails we will try
\begin{equation}
  \label{Me1e1}
  K(X,Y)K\left(\mat{0}{0}{0}{1},0\right)K(X,Y)\in\Pi,\quad\forall K(X,Y)\in\Pi.
\end{equation}
If (\ref{Me1e1}) holds, then we will be able to use
$M_{0}=\tns{\Be_{1}}/2$. If this condition fails as well, then we will try
\begin{equation}
  \label{Me2e2}
  K(X,Y)K\left(\mat{1}{0}{0}{0},0\right)K(X,Y)\in\Pi,\quad\forall K(X,Y)\in\Pi.
\end{equation}
If (\ref{Me2e2}) holds, then we will be able to use
$M_{0}=\tns{\Be_{2}}/2$. If this condition also fails, then we choose $M_{0}=\BI_{2}/2$.

Now let us describe the algorithm for finding the inversion key for the links. In
our case there are 5 of them: 3 are of type (\ref{idlink}) and 2 are of type (\ref{Philink}).
Links $\Hat{\Pi}$ have two components and each can use its own inversion key,
so that $\Hat{M}_{0}=[M_{1},M_{2}]$.
For simplicity of notation we will denote 
\[
\GD_{j}=\hf\BI_{2}-M_{j},\quad j=1,2.
\]
The inversion key $\Hat{M}_{0}$ for $\Hat{\Pi}$, given by the algebra ideal pair
$\Pi=\CI\oplus\Pi'$ via (\ref{idlink}) is identified by checking the following 3
properties
\begin{enumerate}
\item $K(X',Y')K(\GD_{2},0)K(X',Y')\in\Pi'$ for all $K(X',Y')\in\Pi'$ ($M_{2}$
  must be an inversion key for $\Pi'$.)
\item $K(J_{X},J_{Y})K(\GD_{1},0)K(X,Y)+K(X,Y)K(\GD_{1},0)K(J_{X},J_{Y})\in\CI$
 for all $K(X,Y)\in\Pi$ and $K(J_{X},J_{Y})\in\CI$ ($M_{1}$
  must be an inversion key for $\CI$.)
\item $K(X',Y')K(M_{1}-M_{2},0)K(X',Y')\in\CI$ for all $K(X',Y')\in\Pi'$
\end{enumerate}
In the case of $\Hat{\Pi}$ corresponding to an automorphism of $\Pi$ the
inversion key $\Hat{M}_{0}$ is sought in the form $\Hat{M}_{0}=[M_{0},M_{0}]$,
where $M_{0}$ is an inversion key for $\Pi$, satisfying additionally the
relation
\[
\Phi(K(X,Y)K(\GD_{0},0)K(X,Y))=\Phi(K(X,Y))K(\GD_{0},0)\Phi(K(X,Y)).
\]
Let us show that $\Hat{M}_{0}=[0,0]$ is \emph{not} an inversion key for the global
automorphism. $\Hat{M}_{0}=[0,0]$ is equivalent to the property that all global automorphisms
satisfy 
\begin{equation}
  \label{glinvkey}
  \Phi(\SFK^{2})=\Phi(\SFK)^{2}\qquad\forall\SFK\in\Sym(\CT).
\end{equation}
Let us verify that this is not always the case.
There are two branches of the global
automotphisms:
\[
\Phi_{+}(\SFK)=\SFC\SFK\SFC^{T},\qquad\SFC=K(C,0),\quad C\in O(2,\bb{C}).
\]
and 
\[
\Phi_{-}(\SFK)=-\SFC\SFK\SFC^{T},\qquad\SFC=K(iC,0),\quad C\in O(2,\bb{C}).
\]
For $\Phi_{\pm}$ equation (\ref{glinvkey}) is equivalent to
$\SFC^{T}\SFC=\pm\SFI$. We compute for $C\in O(2,\bb{C})$
\[
K(C,0)^{T}K(C,0)=K(C^{H}C,0),\qquad K(iC,0)^{T}K(iC,0)=K(C^{H}C,0).
\]
We see that $\SFC^{T}\SFC=-\SFI$ is never satisfied, while $\SFC^{T}\SFC=\SFI$
holds \IFF $\SFC\in O(2,\bb{R})$. In fact the inversion key for the global
automorphism must be 
\begin{equation}
  \label{globinvkey}
  \Hat{M}_{\rm glob}=\left[\hf\BI_{2},\hf\BI_{2}\right].
\end{equation}
Nevertheless we can write an arbitrary transformation $\Phi_{+}$ as a
superposition of a transformation in $O(2,\bb{R})$ and a transformation
corresponding to
\begin{equation}
  \label{nontriv}
  \BC_{+}(t)=\mat{\cosh t}{i\sinh t}{-i\sinh t}{\cosh t},\quad t\in\bb{R}.
\end{equation}
To obtain transformations $\Phi_{-}$ we only need to compose a transformation
$\Phi_{+}$ with $\Phi_{*}(K(X,Y))=K(-X,Y)$.

\section{Summary of non-redundant and nontrivial ERs and Links}
\begin{center}
  \textbf{Exact relations}
\end{center}
\begin{center}
\begin{tabular}[h]{|c|c|c|}
  \hline
item \# & algebra & inversion key\\
\hline
8& $(\bb{C}\BI,\BGF)$ &$M_{0}=0$\\
\hline
13& ${\rm Ann}(\bb{C}\bra{\Bz_{0}})$ &$M_{0}=0$\\
\hline
17& $(\ClD,\ClD')$ &$M_{0}=\BI_{2}/2$\\
\hline
20&$(W,V_{\infty})$ &$M_{0}=\BI_{2}/2$\\
\hline
21&$(W,V)$ &$M_{0}=\BI_{2}/2$\\
\hline
22&$(\Sym(\bb{C}^{2}),0)$ &$M_{0}=\BI_{2}/2$\\
\hline
\end{tabular}
\end{center}
\begin{center}
  \textbf{Links}
\end{center}
\begin{tabular}[h]{|c|c|c|c|}
  \hline
item \# & algebra & link & inversion key\\
\hline
13& ${\rm Ann}(\bb{C}\bra{\Bz_{0}})$ &$\Pi^{2}=\{0\}$&$M_{0}=\BI_{2}/2$\\
\hline
19&$(W,\bb{R}\BZ_{0})$ &$\Phi_{2}(Y)=Y$,
$\Phi_{0}(\BZ_{0})=\Ga\BZ_{0}$&$\Hat{M}_{0}=[\BI_{2}/2,\BI_{2}/2]$\\
\hline
19&$(W,\bb{R}\BZ_{0})$ &
$\CI=(0,\bb{R}\BZ_{0})$&$\Hat{M}_{0}=[\BI_{2}/2,\BI_{2}/2]$\\
\hline
19&$(W,\bb{R}\BZ_{0})$ &
$\CI=(\tns{\Bz_{0}},0)$&$\Hat{M}_{0}=[\BI_{2}/2,\BI_{2}/2]$\\
\hline
21&$(W,V)$ &$\CI={\rm Ann}(\bb{C}\bra{\Bz_{0}})$&$\Hat{M}_{0}=[\BI_{2}/2,\BI_{2}/2]$\\
\hline
23&$\Sym(\CT)$ &$\Phi(K(X,Y))=K(\pm CXC^{H},CYC^{T})$, $C\in O(2,\bb{C})$&
$\Hat{M}_{0}=[\BI_{2}/2,\BI_{2}/2]$\\
\hline
\end{tabular}

\section{Global Link}
The global link will actually consist of 3 families of links:
\begin{enumerate}
\item  One comes from
choosing $\BC\in O(2,\bb{R})$ (it is obtained using the $\Hat{M}_{0}=[0,0]$
inversion key);
\item The second family comes from the single global automorphism
$\Phi(K(X,Y))=K(-X,Y)$;
\item The third family of links corresponds to the family (\ref{nontriv}) of
  global automorphisms.
\end{enumerate}
 The last two families require using
inversion key (\ref{globinvkey}).

\subsection{$O(2,\bb{R})$ family}
Let us begin with the simplest case, where $\Hat{M}_{0}=[0,0]$, i.e. with
finding links corresponding to global automorphisms defined by $\BC\in
O(2,\bb{R})$. The simplified version is
\[
\SFL_{2}=K(\BI_{2}+\Gb_{2}i\BR_{\perp},0)-K(\BC,0)(\SFL_{1}-K(\BI_{2}+\Gb_{1}i\BR_{\perp},0))K(\BC^{T},0).
\]
For $\BC\in O(2,\bb{R})$ we obtain
\[
\SFL_{2}=(\Gb_{2}\mp(\det\BC)\Gb_{1})\tns{\BR_{\perp}}+(\BC\otimes\BI_{2})\SFL_{1}(\BC^{T}\otimes\BI_{2})
\]
To obtain the general link we replace $\SFL_{j}$ in the above formula with
$(\BGL_{j}^{-1/2}\otimes\BI_{2})\SFL_{j}(\BGL_{j}^{-1/2}\otimes\BI_{2})$.
We obtain, solving for $\SFL_{2}$
\[
\SFL_{2}=(\Gb_{2}\mp(\det\BC)\Gb_{1})\sqrt{\det\BGL_{2}}\tns{\BR_{\perp}}+
(\BGL_{2}^{1/2}\BC\BGL_{1}^{-1/2}\otimes\BI_{2})\SFL_{1}(\BGL_{1}^{-1/2}\BC^{T}\BGL_{2}^{1/2}\otimes\BI_{2})
\]
We observe that a polar decomosition of a matrix implies that every
non-singular $2\times 2$ real matrix can be written as $\BGL^{1/2}\BC$ for
some symmetric positive definite matrix $\BGL$ and $\BC\in
O(2,\bb{R})$. Therefore, we obtain the link
\begin{equation}
  \label{lingloblink}
  \SFL_{2}=\Gb_{0}\SFT+(\BB_{0}\otimes\BI_{2})\SFL_{1}(\BB_{0}^{T}\otimes\BI_{2}),\qquad\SFT=\tns{\BR_{\perp}}.
\end{equation}
where $\Gb_{0}\in\bb{R}$ and $\BB_{0}\in GL(2,\bb{R})$ are parameters of the
family of links, restricted by the requirement that $\SFL_{1,2}$ be
positive definite. We note that by construction for any pair
$\Hat{\SFL}_{0}=[\SFL_{0}^{(1)},\SFL_{0}^{(2)}]$ of isotropic materials, there is a link of
the form (\ref{lingloblink}) passing through $\Hat{\SFL}_{0}$. That means that
\emph{any} link $\SFL_{2}=\mathfrak{L}(\SFL_{1})$ passing through $\Hat{\SFL}_{0}$ can
be obtained from a link passing through
$\Hat{\SFL}_{0}^{0}=[\SFI,\SFI]$. Indeed, let $\mathfrak{L}_{1}$ and
$\mathfrak{L}_{2}$ be the links of the form (\ref{lingloblink}), passing
through $[\SFL_{0}^{(1)},\SFI]$ and $[\SFL_{0}^{(2)},\SFI]$,
respectively. Then, the function
\[
\mathfrak{L}_{0}(\SFL_{1})=\mathfrak{L}_{2}(\mathfrak{L}(\mathfrak{L}_{1}^{-1}(\SFL_{1})))
\]
is a link passing through $\Hat{\SFL}_{0}^{0}$. Hence,
\[
\mathfrak{L}(\SFL_{1})=\mathfrak{L}_{2}^{-1}(\mathfrak{L}_{0}(\mathfrak{L}_{1}(\SFL_{1}))).
\]
We therefore, have the option of deriving only the links passing through
$\Hat{\SFL}_{0}^{0}=[\SFI,\SFI]$, which, combined with (\ref{lingloblink}),
will generate all global links. The same goes for exact relations: we only
need to compute the ones passing though $\SFL_{0}=\SFI$.

\subsection{$\Phi(K(X,Y))=K(-X,Y)$ family}
Next, let us compute the link corresponding to the  map
$\Phi(\SFK(X,Y))=\SFK(-X,Y)$. This has the form 
\begin{equation}
  \label{initlink0}
  [(\SFL_{2}-\SFI)^{-1}+\hf\SFI]^{-1}=K(i\BI_{2},0)[(\SFL_{1}-\SFI)^{-1}+\hf\SFI]^{-1}K(i\BI_{2},0),
\end{equation}
The idea is to identify two
fixed points of this transformation $\SFF_{+}$ and $\SFF_{-}$ and then rewrite
(\ref{initlink0}) as
\[
(\SFL_{2}-\SFF_{-})^{-1}(\SFL_{2}-\SFF_{+})=\SFS_{-}^{-1}(\SFL_{1}-\SFF_{-})^{-1}(\SFL_{1}-\SFF_{+})\SFS_{+}
\]
We can solve (\ref{initlink0}) for $\SFL_{2}$ and then express
$(\SFL_{2}-\SFF_{-})^{-1}(\SFL_{2}-\SFF_{+})$ as
$\SFS_{-}^{-1}(\SFL_{1}-\SFX_{-})^{-1}(\SFL_{1}-\SFX_{+})\SFS_{+}$, which
leads to the formulas for $\SFS_{\pm}$:
\[
\SFS_{\pm}=K(i\BI_{2},0)\SFF_{\pm}.
\]
It is reasonable to look for fixed points in the form 
\begin{equation}
  \label{fixedans}
  \SFF=K\left(\mat{x}{y}{\bra{y}}{x},0\right),
\end{equation}
since, all tensors in (\ref{initlink0}), except $\SFL_{1,2}$, have that form.
We compute (using Maple) that there exists a pair of fixed points of the form
$\SFF_{\pm}=\pm\SFT$. Therefore,
$\SFS_{\pm}=\mp K(\BR_{\perp},0)$. So we have the link
$\Hat{\bb{M}}_{0}$ of the form
\begin{equation}
  \label{gnonlin0}
(\SFL_{2}+\SFT)^{-1}(\SFL_{2}-\SFT)=K(\BR_{\perp},0)(\SFL_{1}+\SFT)^{-1}(\SFL_{1}-\SFT)K(\BR_{\perp},0).
\end{equation}
We can rewrite (\ref{gnonlin0}) by observing that
\begin{equation}
  \label{convform}
  (\SFL+\SFT)^{-1}(\SFL-\SFT)-\SFI=-2(\SFL+\SFT)^{-1}\SFT.
\end{equation}
Thus, (\ref{gnonlin0}) becomes
\[
(\SFL_{2}+\SFT)^{-1}=\SFT-(\BR_{\perp}\otimes\BI_{2})(\SFL_{1}+\SFT)^{-1}(\BR_{\perp}^{T}\otimes\BI_{2})
\]
Applying the link (\ref{lingloblink}) ($\SFL_{2}'=\SFL_{2}+\SFT$ and
$\SFL_{1}'=(\BR_{\perp}\otimes\BI_{2})\SFL_{1}(\BR_{\perp}^{T}\otimes\BI_{2})+\SFT$)
we can simplify our link to 
\begin{equation}
  \label{gnonlinalt}
\SFL_{2}^{-1}=\SFT-\SFL_{1}^{-1}.
\end{equation}

\subsection{(\ref{nontriv}) family}
Finally, we need to compute the link, corresponding to $\BC\in O(2,\bb{C})$,
given by (\ref{nontriv}). We first compute the link $\Hat{\bb{M}}_{0}$ defined
by
\begin{equation}
  \label{initlink}
  [(\SFL_{2}-\SFI)^{-1}+\hf\SFI]^{-1}=\SFC[(\SFL_{1}-\SFI)^{-1}+\hf\SFI]^{-1}\SFC,
\end{equation}
where $\SFC=K(\BC,0)$, where $\BC$ is given by (\ref{nontriv}).

We rewrite (\ref{initlink}) in a symmetrical way with respect to
$\SFL_{1}$ and $\SFL_{2}$ using the fixed points $\SFF_{\pm}$ of the form 
\begin{equation}
  \label{fixedans0}
  \SFF=K\left(\mat{0}{y}{\bra{y}}{0},0\right),
\end{equation}
as we discovered before. We also find that
\[
\SFS_{\pm}=\hf(\SFC^{-1}-\SFC)(\SFI-\SFF_{\pm})+\SFC.
\]

Using Maple, we find that $\SFF_{\pm}=\pm\SFT$. Then
$\SFS_{\pm}=e^{\mp t}\SFI$, resulting in the formula
\[
(\SFL_{2}+\SFT)^{-1}(\SFL_{2}-\SFT)=e^{-2t}(\SFL_{1}+\SFT)^{-1}(\SFL_{1}-\SFT),
\]
Equivalently, using (\ref{convform}),
\begin{equation}
  \label{gnonlin}
(\SFL_{2}+\SFT)^{-1}=e^{-2t}(\SFL_{1}+\SFT)^{-1}+e^{-t}\sinh(t)\SFT.
\end{equation}
Applying the link (\ref{lingloblink}) we can simplify our link to 
\begin{equation}
  \label{gnonlin1}
\SFL_{2}^{-1}=\Ga_{0}\SFT+\SFL_{1}^{-1},
\end{equation}
which forms a group of transformations $\mathfrak{L}_{\Ga_{0}}$, such that
$\mathfrak{L}_{\Ga_{0}}\circ\mathfrak{L}_{\Gb_{0}}=\mathfrak{L}_{\Ga_{0}+\Gb_{0}}$

\subsection{General form and properties of the global link}
Combining the results obtained so far we conclude that any global link can be
obtained as a superposition of the following subgroups of links
\begin{equation}
  \label{globlinkcomb}
  \begin{cases}
\SFL_{2}^{-1}=\SFT-\SFL_{1}^{-1},\\
\SFL_{2}^{-1}=\Ga_{0}\SFT+\SFL_{1}^{-1},\\
\SFL_{2}=\Gb_{0}\SFT+\SFL_{1},\\
\SFL_{2}=(\BB_{0}\otimes\BI_{2})\SFL_{1}(\BB_{0}^{T}\otimes\BI_{2}).
  \end{cases}
\end{equation}
The most general transformation that can be made by compositing
transformations (\ref{globlinkcomb}) with each other is
\begin{equation}
  \label{globauto}
  \Psi(\SFL)=(\BB_{0}\otimes\BI_{2})\SFT(\Ga_{1}\SFL+\Gb_{1}\SFT)^{-1}(\Ga_{0}\SFL+\Gb_{0}\SFT)(\BB_{0}^{T}\otimes\BI_{2}).
\end{equation}
We remark that
\[
 \Psi(\SFL)=(\BB_{0}\otimes\BI_{2})(\Ga_{0}\SFL+\Gb_{0}\SFT)(\Ga_{1}\SFL+\Gb_{1}\SFT)^{-1}\SFT(\BB_{0}^{T}\otimes\BI_{2}).
\]
If $\Psi(\SFL)$ is given by (\ref{globauto}) we will write
$\Psi_{\BA_{0},\BB_{0}}(\SFL)$ to refer to it, where
\[
\BA_{0}=\mat{\Ga_{0}}{\Gb_{0}}{\Ga_{1}}{\Gb_{1}}.
\]
Different pairs of matrices $\{\BA,\BB\}\subset GL(2,\bb{R})$ can define
the same transformation $\Psi_{\BA,\BB}$. Specifically, 
\begin{equation}
  \label{projinv}
  \Psi_{\Gl\BA,\BB}=\Psi_{\BA,\BB},\quad\Psi_{\BA,\Gl\BB}=\Psi_{\BA_{\Gl},\BB},\qquad\BA_{\Gl}=\mat{\Gl^{2}}{0}{0}{1}\BA
\end{equation}
for any nonzero real number $\Gl$. Thus, \WLOG, we may assume that
$|\det\BA|=|\det\BB|=1$. Even with this assumption we still have symmetries 
\[
\Psi_{-\BA,\BB}=\Psi_{\BA,\BB},\qquad\Psi_{\BA,-\BB}=\Psi_{\BA,\BB}.
\]

Let us derive the formula for superposition of two transformations
(\ref{globauto}). We note that $\Psi_{\BA,\BI_{2}}(\SFL)=M_{\BA}(\SFL\SFT)\SFT$, where
$M_{\BA}(z)$ is a fractional-linear M\"obius transformation with real matrix
$\BA$. The composition law for M\"obius transformations then imply that
\[
\Psi_{\BA_{1},\BI_{2}}\circ\Psi_{\BA_{2},\BI_{2}}=\Psi_{\BA_{1}\BA_{2},\BI_{2}}.
\]
The composition formula
$\Psi_{\BI_{2},\BB_{1}}\circ\Psi_{\BA,\BB_{2}}=\Psi_{\BA,\BB_{1}\BB_{2}}$
is completely evident. Finally, a direct calculation shows that
\[
\Psi_{\BA,\BI_{2}}\circ\Psi_{\BI_{2},\BB}=\Psi_{\BA^{\BB},\BB},\qquad
\BA^{\BB}=\mat{\det\BB}{0}{0}{1}^{-1}\BA\mat{\det\BB}{0}{0}{1}.
\]
where we have used the projective invariance property (\ref{projinv}).
These formulas allow us to derive the full composion formula
\begin{multline*}
  \Psi_{\BA_{1},\BB_{1}}\circ\Psi_{\BA_{2},\BB_{2}}=\Psi_{\BI_{2},\BB_{1}}\circ\Psi_{\BA_{1},\BI_{2}}\circ
\Psi_{\BA_{2}^{\BB_{2}^{-1}},\BI_{2}}\circ\Psi_{\BI_{2},\BB_{2}}=
\Psi_{\BI_{2},\BB_{1}}\circ\Psi_{\BA_{1}\BA_{2}^{\BB_{2}^{-1}},\BI_{2}}\circ\Psi_{\BI_{2},\BB_{2}}=\\
\Psi_{\BI_{2},\BB_{1}}\circ\Psi_{(\BA_{1}\BA_{2}^{\BB_{2}^{-1}})^{\BB_{2}},\BB_{2}}=
\Psi_{\BI_{2},\BB_{1}}\circ\Psi_{\BA_{1}^{\BB_{2}}\BA_{2},\BB_{2}}=\Psi_{\BA_{1}^{\BB_{2}}\BA_{2},\BB_{1}\BB_{2}}.
\end{multline*}
We note that if $\det\BB_{2}=1$, then
$\Psi_{\BA_{1},\BB_{1}}\circ\Psi_{\BA_{2},\BB_{2}}=\Psi_{\BA_{1}\BA_{2},\BB_{1}\BB_{2}}$. 
The most general transformation $\Psi_{\BA,\BB}$, such that $\Psi_{\BA,\BB}(\SFI)=\SFI$ has the
form
\begin{equation}
  \label{Idfixed}
  \BA=\mat{\Ga_{0}}{\Gb_{0}}{\Gb_{0}}{\Ga_{0}},\quad\BB\in O(2,\bb{R}),\qquad
|\Ga_{0}^{2}-\Gb_{0}^{2}|=1.
\end{equation}

\section{Formulas for computing ERs}
The relation bewteen $K(X,Y)$ and the $2\times 2$ block-matrix representation
(\ref{blckM}) is
\begin{equation}
  \label{K2block}
  K(X,Y)=\mat{\Gvf(X_{11})+\psi(Y_{11})}{\Gvf(X_{12})+\psi(Y_{12})}{\Gvf(X_{21})+\psi(Y_{21})}{\Gvf(X_{22})+\psi(Y_{22})},
\end{equation}
where 
\[
\Gvf(\Ga+i\Gb)=\mat{\Ga}{-i\Gb}{i\Gb}{\Ga},\qquad\psi(\Ga+i\Gb)=\mat{\Ga}{\Gb}{\Gb}{-\Ga}.
\]

When $M_{0}=0$ we have
\begin{equation}
  \label{M0red}
  \bb{M}_{0}=\{\SFI+\SFK\in\Sym^{+}(\CT):\SFK\in\Pi_{0}\}.
\end{equation}

For $M_{0}=\BI_{2}/2$ we have
\begin{equation}
  \label{L2K}
  \SFK=\left[(\SFL-\SFI)^{-1}+\hf\BI_{2}\right]^{-1}=2(\SFL+\SFI)^{-1}(\SFL-\SFI)=
2\SFI-4(\SFL+\SFI)^{-1}
\end{equation}
Equivalently $4(\SFL+\SFI)^{-1}=2\SFI-\SFK$. Solving for $\SFL$ we obtain
\[
\SFL=\SFI+2(2\SFI-\SFK)^{-1}\SFK=2(\SFI-\SFK/2)^{-1}-\SFI.
\]
Since our goal is to compute the image of the \emph{subspace} $\Pi$ under the
above transformation, we may just as well use the formula
\begin{equation}
  \label{invform0}
 \bb{M}_{0}=\{\SFL\in\Sym^{+}(\CT):\SFL=2(\SFI+\SFK)^{-1}-\SFI,\ \SFK\in\Pi_{0}\}.
\end{equation}
Sometimes it might be easier to characterize $\bb{M}_{0}$ by computing $\SFK$
in terms of $\SFL$ and writing equations satisfied by $\SFK$ in terms of
$\SFL$. Then 
\begin{equation}
  \label{invform1}
 \bb{M}_{0}=\left\{\SFL\in\Sym^{+}(\CT):(\SFL+\SFI)^{-1}-\hf\SFI\in\Pi_{0}\right\}.
\end{equation}
Both formulas require inverting the $2\times 2$ block-matrices.
One may choose to compute block-matrix inverse in two ways: in the $2\times 2$
block-matrix notation or in the $K(X,Y)$ notation. The $2\times 2$
block-matrix formalism is standard. Let us assume that $\BF_{11}$ in
\[
\SFF=\mat{\BF_{11}}{\BF_{12}}{\BF_{21}}{\BF_{22}},
\]
is invertible. Then, in order to invert $\SFF$ we need to solve the system of
equations:
\[
\begin{cases}
  \BF_{11}\Bu_{1}+\BF_{12}\Bu_{2}=\Bv_{1},\\
  \BF_{21}\Bu_{1}+\BF_{22}\Bu_{2}=\Bv_{2}
\end{cases}
\]
We solve it using the method of elimination. We solve the first equation for $\Bu_{1}$:
\[
\Bu_{1}=\BF_{11}^{-1}\Bv_{1}-\BF_{11}^{-1}\BF_{12}\Bu_{2},
\]
and substitute the result into the second equation:
 \[
\BF_{21}\BF_{11}^{-1}\Bv_{1}+(\BF_{22}-\BF_{21}\BF_{11}^{-1}\BF_{12})\Bu_{2}=\Bv_{2}.
\]
We then solve this for $\Bu_{2}$:
\[
\Bu_{2}=-(\BF_{22}-\BF_{21}\BF_{11}^{-1}\BF_{12})^{-1}\BF_{21}\BF_{11}^{-1}\Bv_{1}+
(\BF_{22}-\BF_{21}\BF_{11}^{-1}\BF_{12})^{-1}\Bv_{2},
\]
and substitute this into the formula for $\Bu_{1}$:
\[
\Bu_{1}=(\BF_{11}^{-1}+\BF_{11}^{-1}\BF_{12}(\BF_{22}-\BF_{21}\BF_{11}^{-1}\BF_{12})^{-1}\BF_{21}\BF_{11}^{-1})\Bv_{1}-
\BF_{11}^{-1}\BF_{12}(\BF_{22}-\BF_{21}\BF_{11}^{-1}\BF_{12})^{-1}\Bv_{2}
\]
A little matrix algebra shows that
\[
\Bu_{1}=(\BF_{11}-\BF_{12}\BF_{22}^{-1}\BF_{21})^{-1}\Bv_{1}-
\BF_{11}^{-1}\BF_{12}(\BF_{22}-\BF_{21}\BF_{11}^{-1}\BF_{12})^{-1}\Bv_{2}
\]
This gives us a formula for $\SFF^{-1}$:
\[
\SFF^{-1}=\mat{\BS_{11}^{-1}}{-\BF_{11}^{-1}\BF_{12}\BS_{22}^{-1}}{-\BS_{22}^{-1}\BF_{21}\BF_{11}^{-1}}{\BS_{22}^{-1}},
\]
where matrices
\[
\BS_{11}=\BF_{11}-\BF_{12}\BF_{22}^{-1}\BF_{21},\qquad
\BS_{22}=\BF_{22}-\BF_{21}\BF_{11}^{-1}\BF_{12}
\]
are called Schur complements of $\BF_{11}$ and $\BF_{22}$, respectively. Of
course,
\[
\BF_{11}^{-1}\BF_{12}\BS_{22}^{-1}=\BS_{11}^{-1}\BF_{12}\BF_{22}^{-1},\qquad
\BS_{22}^{-1}\BF_{21}\BF_{11}^{-1}=\BF_{22}^{-1}\BF_{21}\BS_{11}^{-1}.
\]
Therefore, we can write $\SFF^{-1}$ in two equivalent more symmetrical forms
\begin{equation}
  \label{Finv}
  \SFF^{-1}=\mat{\BS_{11}^{-1}}{-\BS_{11}^{-1}\BF_{12}\BF_{22}^{-1}}{-\BS_{22}^{-1}\BF_{21}\BF_{11}^{-1}}{\BS_{22}^{-1}}=\mat{\BS_{11}^{-1}}{-\BF_{11}^{-1}\BF_{12}\BS_{22}^{-1}}{-\BF_{22}^{-1}\BF_{21}\BS_{11}^{-1}}{\BS_{22}^{-1}}.
\end{equation}
Equivalently,
\[
\SFF^{-1}=\mat{\BS_{11}}{0}{0}{\BS_{22}}^{-1}\mat{\BF_{11}}{-\BF_{12}}{-\BF_{21}}{\BF_{22}}
\mat{\BF_{11}}{0}{0}{\BF_{22}}^{-1}
\]
\[
\SFF^{-1}=
\mat{\BF_{11}}{0}{0}{\BF_{22}}^{-1}\mat{\BF_{11}}{-\BF_{12}}{-\BF_{21}}{\BF_{22}}
\mat{\BS_{11}}{0}{0}{\BS_{22}}^{-1}
\]

A necessary and sufficient condition for $\SFF$ to be in $\Sym^{+}(\CT)$ is
$\BF_{11}>0$ and $\BS_{22}>0$ (or $\BF_{22}>0$ and $\BS_{11}>0$).

We can also derive formulas for $K(X,Y)^{-1}$. In this case we solve the equation
\[
X\Bu+Y\bra{\Bu}=\Bv
\]
for $\Bu$ or $\bra{\Bu}$:
\[
\Bu=X^{-1}\Bv-X^{-1}Y\bra{\Bu},
\]
or
\[
\bra{\Bu}=Y^{-1}\Bv-Y^{-1}X\Bu.
\]
Taking complex conjugates we get
\[
\bra{\Bu}=\bra{X}^{-1}\bra{\Bv}-\bra{X}^{-1}\bra{Y}\Bu,
\]
or
\[
\Bu=\bra{Y}^{-1}\bra{\Bv}-\bra{Y}^{-1}\bra{X}\bra{\Bu}.
\]
We then substitute this into the original equation:
\[
(X-Y\bra{X}^{-1}\bra{Y})\Bu=\Bv-Y\bra{X}^{-1}\bra{\Bv},
\]
or
\[
(Y-X\bra{Y}^{-1}\bra{X})\bra{\Bu}=\Bv-X\bra{Y}^{-1}\bra{\Bv}.
\]
we now solve for $\Bu$ (or $\bra{\Bu}$):
\[
\Bu=(X-Y\bra{X}^{-1}\bra{Y})^{-1}\Bv-(X-Y\bra{X}^{-1}\bra{Y})^{-1}Y\bra{X}^{-1}\bra{\Bv},
\]
or
\[
\bra{\Bu}=(Y-X\bra{Y}^{-1}\bra{X})^{-1}\Bv-(Y-X\bra{Y}^{-1}\bra{X})^{-1}X\bra{Y}^{-1}\bra{\Bv}.
\]
Hence, we obtain two formulas for $K(X,Y)^{-1}$:
\begin{equation}
  \label{Kinv}
  K(X,Y)^{-1}=K(S_{X}^{-1},-S_{X}^{-1}Y\bra{X}^{-1})=
K(-S_{Y}^{-1}\bra{X}Y^{-1},S_{Y}^{-1})=K(S_{X}^{-1},S_{Y}^{-1})
\end{equation}
where
\[
S_{X}=X-Y\bra{X}^{-1}\bra{Y},\qquad S_{Y}=\bra{Y}-\bra{X}Y^{-1}X
\]
play the role of Schur complements of $X$ and $Y$, respectively.

\section{$\Pi_{0}=(\bb{C}\BI_{2},\BGF)$ and $\Pi_{0}={\rm Ann}(\bb{C}\bra{z_{0}})$} 
For $\Pi_{0}=(\bb{C}\BI_{2},\BGF)$ we have
\[
  \SFK=\mat{\psi(z)}{0}{0}{\psi(z)}+\mat{\Gvf(\Ga)}{\Gvf(-i\Gb)}{\Gvf(i\Gb)}{\Gvf(\Ga)}=
\mat{\psi(z)+\Gvf(\Ga)}{\Gvf(-i\Gb)}{\Gvf(i\Gb)}{\psi(z)+\Gvf(\Ga)}.
\]
Let $\BK=\psi(z)+\Gvf(\Ga)\in\Sym(\bb{R}^{2})$. This is our change of variables. There is a
1-1 correspondence between $\bb{C}\times\bb{R}$ and $\Sym(\bb{R}^{2})$, given
by $\BK=\psi(z)+\Gvf(\Ga)\in\Sym(\bb{R}^{2})$. Thus, we obtain
\[
\SFK=\mat{\BK}{-\Gb\BR_{\perp}}{\Gb\BR_{\perp}}{\BK}=\BI_{2}\otimes\BK+\Gb\BR_{\perp}\otimes\BR_{\perp}.
\]
Recall that
$\SFL_{0}=\tns{\BI_{2}}$. Then
\[
\SFL_{0}^{0}+\SFK=\BI_{2}\otimes(\BK+\BI_{2})+\Gb\BR_{\perp}\otimes\BR_{\perp}.
\]
We conclude that (denoting $\BL=\BK+\BI_{2}$)
\[
\bb{M}_{0}=\{\BI_{2}\otimes\BL+\Gb\BR_{\perp}\otimes\BR_{\perp}\in\Sym^{+}(\CT):
\BL\in\Sym(\bb{R}^{2}),\ \Gb\in\bb{R}\}.
\]
Finally (and this is optional, since the global link can map $\bb{M}_{0}$ into
$\bb{M}$), $\bb{M}=\{\SFC_{0}\SFL\SFC_{0}:\SFL\in\bb{M}_{0}\}$, where
$\SFC_{0}=\BGL_{0}^{1/2}\otimes\BI_{2}$. We compute
\[
(\BGL_{0}^{1/2}\otimes\BI_{2})(\BI_{2}\otimes\BL+\nu\BR_{\perp}\otimes\BR_{\perp})
(\BGL_{0}^{1/2}\otimes\BI_{2})=\BGL_{0}\otimes\BL+\nu\sqrt{\det\BGL_{0}}\BR_{\perp}\otimes\BR_{\perp}
\]
Hence, (introducing a new variable $t=\nu\sqrt{\det\BGL_{0}}$)
\begin{equation}
    \label{ER8}
\bb{M}=\{\BGL_{0}\otimes\BL+t\BR_{\perp}\otimes\BR_{\perp}:
\BL\in\Sym^{+}(\bb{R}^{2}),\ |t|<\sqrt{\det(\BGL_{0}\BL)}\}.
 \end{equation}
This is a whole family of exact relation manifolds (one for each choice of
$\BGL_{0}$) corresponding to $\Pi_{0}=(\bb{C}\BI_{2},\BGF)$. We have computed
$\bb{M}$ for the sole reason that it just as beautiful as $\bb{M}_{0}$. In our
next example, this does not seem to be the case, so we leave the exact
relations in the $\bb{M}_{0}$ form.

For $\Pi_{0}={\rm Ann}(\bb{C}\bra{z_{0}})$ we have
\[
  \SFK=\mat{\psi(z)}{\psi(-iz)}{\psi(-iz)}{-\psi(z)}+\mat{\Gvf(\Ga)}{\Gvf(i\Ga)}{\Gvf(-i\Ga)}{\Gvf(\Ga)}=
\mat{\psi(z)+\Gvf(\Ga)}{\psi(-iz)+\Gvf(i\Ga)}{\psi(-iz)+\Gvf(-i\Ga)}{-\psi(z)+\Gvf(\Ga)}.
\]
Let $\BK=\psi(z)+\Gvf(\Ga)\in\Sym(\bb{R}^{2})$. This is our change of
variables. Then
\[
\SFK=\mat{\BK}{\BK\BR_{\perp}}{-\BR_{\perp}\BK}{\cof(\BK)}.
\]
\[
\SFL=\SFI+\SFK=\mat{\BK+\BI_{2}}{\BK\BR_{\perp}}
{-\BR_{\perp}\BK}{\cof(\BK+\BI_{2})}.
\]
Denoting $\BL=\BK+\BI_{2}$ we obtain
\begin{equation}
  \label{Annz0}
  \SFL=\mat{\BL}{(\BL-\BI_{2})\BR_{\perp}}{\BR_{\perp}^{T}(\BL-\BI_{2})}{\cof(\BL)}=
\mat{\BL}{\BL\BR_{\perp}}{\BR_{\perp}^{T}\BL}{\cof(\BL)}+\SFT.
\end{equation}
One can check that $\SFL>0$ \IFF $\BL>\BI_{2}/2$ in
the sense of quadratic forms. The attempts to compute $\bb{M}$ have not lead
to a very beatiful representation of this exact relation, so we leave it in
the $\bb{M}_{0}$ form. Application of the inversion formula with
$M_{0}=\BI_{2}/2$ will lead the volume fraction relation in the form
\begin{equation}
  \label{FR}
\BL_{*}^{-1}=\av{\BL^{-1}}.
\end{equation}

\section{$\Pi_{0}=(\Sym(\bb{C}^{2}),0)$ and $\Pi_{0}=(\ClD,\ClD')$}
\[
\SFK=\mat{\psi(x)}{\psi(c)}{\psi(c)}{\psi(y)},\qquad
\{x,y,c\}\subset\bb{C}.
\]
We can compute the ER using the complex inversion formula
(\ref{Kinv}). According to this formula we have
\[
\SFL+\SFI=\left(\SFK+\hf\SFI\right)^{-1}=K\left(\hf\BI,\BY\right)^{-1}=K\left(2(\BI-4\BY\bra{\BY})^{-1},4(4\bra{\BY}-\BY^{-1})^{-1}\right).
\]
Thus, if we write $\SFL+\SFI=K(\BU,\BZ)$, then
\[
\BU=2(\BI-4\BY\bra{\BY})^{-1},\qquad\BZ=-4(\BI-4\BY\bra{\BY})^{-1}\BY
\]
Thus, $2\BY=-\BU^{-1}\BZ$. The symmetry of $\BY$ is equivalent to the equation
$\BU^{-1}\BZ=\BZ\bra{\BU^{-1}}$, or equivalently, to
\[
\BZ\bra{\BU}=\BU\BZ.
\]
Again, due to the symmetry of $\BY$ we have
\[
2\bra{\BY}=2\BY^{*}=-\bra{\BZ}\BU^{-1}.
\]
Using this in the formula for $\BU$ to eliminate $\BY$ and $\bra{\BY}$ we have
\[
2\BU^{-1}=\BI-\BU^{-1}\BZ\bra{\BZ}\BU^{-1}\eqv\BZ\bra{\BZ}=\BU^{2}-2\BU=(\BU-\BI)^{2}-\BI.
\]
Noting that $\SFL=K(\BU-\BI,\BZ)=K(\BV,\BZ)$ we obtain the description of this
exact relation as the system of equations:
\[
\BZ\bra{\BV}=\BV\BZ,\qquad\BV^{2}-\BZ\bra{\BZ}=\BI.
\]
This suggests that these equations can be rewritten in terms of the $4\times
4$ matrix multiplication. Let $\mathfrak{I}:\Sym(\CT)\to\Sym(\CT)$ be given by
its action $\mathfrak{I}(K(X,Y))=K(X,-Y)$. Then it is easy to see that our
exact relation says
\[
\SFL\mathfrak{I}(\SFL)=\SFI.
\]
We observe that
\[
\mathfrak{I}(\SFL)=(\BI\otimes\BR_{\perp})\SFL(\BI\otimes\BR_{\perp})^{T}.
\]
Hence, we have an alternative representation of the ER $(\Sym(\bb{C}^{2}),0)$:
\begin{equation}
  \label{ER22}
  \bb{M}=\{\SFL>0:\SFL(\BI\otimes\BR_{\perp})\SFL(\BI\otimes\BR_{\perp})^{T}=\SFI\}=
\{\SFL>0:\SFL(\BI\otimes\BR_{\perp})\SFL=\BI\otimes\BR_{\perp}\}.
\end{equation}
In block-components we can rewrite this as a system of equations
\begin{equation}
  \label{ER22comp}
  \begin{cases}
    \frac{\BL_{11}}{\det\BL_{11}}=\BL_{11}-\BL_{12}\BL_{22}^{-1}\BL_{12}^{T},\\
    \det\BL_{11}+\det\BL_{12}=1,\\
    \det\BL_{22}+\det\BL_{12}=1,
  \end{cases}
\end{equation}
where the last equation is redundant and is added to the system for the sake
of the symmetry. For the sake of reference
\begin{equation}
  \label{ER22fin}
  \bb{M}=\left\{\SFL>0:\BL_{11}=-\frac{\BL_{12}\cof(\BL_{22})\BL_{12}^{T}}{\det\BL_{12}},\
\det\BL_{22}+\det\BL_{12}=1\right\}.
\end{equation}

The form (\ref{ER22}) of the ER corresponding to $\Pi_{0}=(\Sym(\bb{C}^{2}),0)$ suggests
looking for other isotropic tensors $\SFA$, such that
$\SFL\SFA\SFL=\SFB=$const is an ER. Hence, we are looking for an isotropic
tensor $\SFA=K(\BA,0)$, such that $\Pi=\{\SFK:\SFK K(\BA,0)+K(\BA,0)\SFK=0\}$
is one of the algebras in our list. It is not hard to compute that the only
other choice of $\BA$ besides $\mat{i}{0}{0}{i}$ is $\mat{i}{0}{0}{-i}$, which
corresponds to $\Pi=(\ClD,\ClD')$ and
\[
\SFA=K\left(\mat{i}{0}{0}{-i},0\right)=\mat{\BR_{\perp}}{0}{0}{-\BR_{\perp}}=
\BJ\otimes\BR_{\perp},\quad\BJ=\psi(1)=\mat{1}{0}{0}{-1}.
\]
The Maple calculation confirms that
\begin{equation}
  \label{ER17}
  \bb{M}=\{:\SFL(\BJ\otimes\BR_{\perp})\SFL=\BJ\otimes\BR_{\perp}\},
\end{equation}
since the manifold has the same dimension as $\Pi_{0}$ and every matrix $\SFL$
of the formb $\SFL=2\SFF^{-1}-\SFI$, $\SFK\in\Pi_{0}$ satisfies (\ref{ER17}).
In block-components we can rewrite this as a system of equations:
\begin{equation}
  \label{ER17comp}
  \begin{cases}
    \frac{\BL_{11}}{\det\BL_{11}}=\BL_{11}-\BL_{12}\BL_{22}^{-1}\BL_{12}^{T},\\
    \det\BL_{11}-\det\BL_{12}=1,\\
    \det\BL_{22}-\det\BL_{12}=1,
  \end{cases}
\end{equation}
where the last equation is redundant and is added to the system for the sake
of the symmetry.
\begin{equation}
  \label{ER17fin}
  \bb{M}=\left\{\SFL>0:\BL_{11}=\frac{\BL_{12}\cof(\BL_{22})\BL_{12}^{T}}{\det\BL_{12}},\
\det\BL_{22}-\det\BL_{12}=1\right\}.
\end{equation}
The next idea comes from examining an application of the theory.

\section{$\Pi_{0}=(W,V_{\infty})$}
Figure~\ref{fig:bincomp} shows that $(W,V_{\infty})$, corresponding to
$\BY\not=0$ and $\BF\in\bb{R}\BZ_{0}$ is a limiting case of the generic
situation, as $\BF\to\BF_{0}\in\bb{R}\BZ_{0}$. The limiting position of a
family of ERs is also an ER. In other words, the set of ERs is closed in the
Grassmannian of $\Sym(\CT)$. Our study of the ERs applicable to binary
composites made of two isotropic phases shows that the $(W,V_{\infty})$ ER is
the limiting position of images of both $(\ClD,\ClD)$ and $(\ClD,\ClD')$ under
the action of global outomorphisms
\[
K(X,Y)\mapsto K(C(c)XC(c)^{H},C(c)YC(c)^{T}),\qquad 
C(c)=\mat{\cos(c)}{\sin(c)}{-\sin(c)}{\cos(c)},\quad c\in\bb{C}.
\]
We have understood the action of the automorphism as follows:
\begin{itemize}
\item \textbf{Action on X}. We can decompose every
  $X\in\mathfrak{H}(\bb{C}^{2})$ as
\[
X=\psi(x)+\xi\BZ_{0}+\eta\bra{\BZ_{0}},\qquad\BZ_{0}=\Bz_{0}\otimes\bra{\Bz_{0}},\
\Bz_{0}=[1,-i],\ x\in\bb{C},\ \{\xi,\eta\}\subset\bb{R}.
\]
Then
\[
C(c)XC(c)^{H}=\psi(e^{-2i\re(c)}x)+\xi e^{2\im(c)}\BZ_{0}+\eta e^{-2\im(c)}\bra{\BZ_{0}}.
\]
\item \textbf{Action on Y}. We can decompose every
  $Y\in\Sym(\bb{C}^{2})$ as
\[
Y=a\BI_{2}+y\tns{\Bz_{0}}+z\tns{\bra{\Bz_{0}}},\quad\{a,y,z\}\subset\bb{C}.
\]
Then,
\[
C(c)YC(c)^{T}=a\BI_{2}+ye^{-2ic}\tns{\Bz_{0}}+ze^{2ic}\tns{\bra{\Bz_{0}}}.
\]
\end{itemize}
We then describe the $(W,V_{\infty})$ ER in the appropriate basis:
\[
W=\{a\BI+y\tns{\Bz_{0}}:\{a,y\}\subset\bb{C}\},\quad
V_{\infty}=\{\psi(it)+\xi\BZ_{0}:\{t,\xi\}\subset\bb{R}\}.
\]
We also describe the $(\ClD,\ClD')$ ER in the same basis:
\[
\ClD_{Y}=\{b\BI_{2}+w(\tns{\Bz_{0}}+\tns{\bra{\Bz_{0}}}):\{b,w\}\subset\bb{C}\},\quad
\ClD'_{X}=\{\psi(is)+\eta(\BZ_{0}-\bra{\BZ_{0}}):\{s,\eta\}\subset\bb{R}\}.
\]
Now we set $c=iM$, where $M>0$ is large. We then reparametrize
$(\ClD,\ClD')$ as follows:
\[
b=a,\quad w=ye^{-2M},\quad s=t,\quad\eta=e^{-2M}\xi.
\]
Then
\[
C(c)\cdot(\ClD,\ClD')=(a\BI+y\tns{\Bz_{0}}+ye^{-4M}\tns{\bra{\Bz_{0}}},
\psi(it)+\xi\BZ_{0}-\xi e^{-4M}\bra{\BZ_{0}}).
\]
This shows that
\[
\lim_{M\to+\infty}C(iM)\cdot(\ClD,\ClD')=(W,V_{\infty}).
\]
We take as a starting point formula (\ref{ER17}) for $(\ClD,\ClD')$ and
formula (\ref{gnonlin}) as the action of the subgroup $C(it)$ on material
tensors. We then try to apply this transformation with $t=M$ to formula
(\ref{ER17}), discarding exponentially small terms\footnote{Terms like
  $e^{-2M}\SFL$ are not necessarily exponentially small, since some components
  of $\SFL$ can be exponentially large.} along the way. At the end we obtain
\begin{equation}
  \label{ER20}
  (\SFL-\SFT)(\BJ\otimes\BR_{\perp})(\SFL-\SFT)=0,\qquad\BJ=\psi(1)=\mat{1}{0}{0}{-1}.
\end{equation}
Maple verification confirms the correctness of (\ref{ER20}). We can also
rewrite (\ref{ER20}) in terms of the block-components of $\SFL$ in the form of
4 independent equations. If we write
\[
\SFL=\mat{\BL_{11}}{\BL_{12}}{\BL_{12}^{T}}{\BL_{22}},
\]
then (\ref{ER20}) is equivalent to
\begin{equation}
  \label{ER20comp}
  \begin{cases}
      \BL_{11}=(\BL_{12}+\BR_{\perp})\BL_{22}^{-1}(\BL_{12}+\BR_{\perp})^{T},\\
      \det\BL_{22}=\det\BS+(\Gb+1)^{2}=\det(\BL_{12}+\BR_{\perp}).
  \end{cases}
\end{equation}
Equivalently,
\begin{equation}
  \label{ER20compa}
  \begin{cases}
      \BL_{22}=(\BL_{12}+\BR_{\perp})^{T}\BL_{11}^{-1}(\BL_{12}+\BR_{\perp}),\\
      \det\BL_{11}=\det(\BL_{12}+\BR_{\perp}).
  \end{cases}
\end{equation}
We can also write this ER in terms of $\BM=\BL_{11}^{-1}(\BL_{12}+\BR_{\perp})$:
\[
\BL_{12}=\BL_{11}\BM-\BR_{\perp},\quad\BL_{22}=\BM^{T}\BL_{11}\BM,\quad\det\BM=1.
\]
\begin{equation}
  \label{ER20fin}
  \bb{M}=\left\{\mat{\BL}{\BL\BM-\BR_{\perp}}{\BM^{T}\BL+\BR_{\perp}}{\BM^{T}\BL\BM}:\det\BM=1,\
\BL>0,\ \BL+2\BR_{\perp}\BM\det\BL<0\right\}.
\end{equation}
The choice of the parameter $t=\infty$ in the family $V_{t}$ was probably not
the best. A better choice is $t=0$, giving
\begin{equation}
  \label{ER20t0}
  (\SFL-\SFT)(\BJ'\otimes\BR_{\perp})(\SFL-\SFT)=0,\qquad\BJ'=\psi(i)=\mat{0}{1}{1}{0}.
\end{equation}
Equivalently (obtained by substituting (\ref{ER21fin}) into (\ref{ER20t0})),
\begin{equation}
  \label{ER20t0fin}
  \bb{M}=\left\{\mat{\BL}{\BL\BM-\BR_{\perp}}{\BM^{T}\BL+\BR_{\perp}}{\BM^{T}\BL\BM}:\Trc\BM=0,\
\BL>0,\ \BL+2\BR_{\perp}\BM\det\BL<0\right\}.
\end{equation}

\section{$\Pi_{0}=(W,V)$}
Here are several different ways to describe $(W,V)$.
\[
W=\{Y\in\Sym(\bb{C}^{2}):Y_{22}-Y_{11}+2iY_{12}=0\},\qquad
V=\{X\in\CH(\bb{C}^{2}):2\im(X_{12})=\Trc X\}.
\]
or
\[
\Pi_{0}=\left\{\SFK=\mat{\BK_{11}}{\BK_{12}}{\BK_{12}^{T}}{\BK_{22}}:
\BK_{11}+\BR_{\perp}\BK_{22}\BR_{\perp}^{T}+\BK_{12}\BR_{\perp}-\BR_{\perp}\BK_{12}^{T}=0.
\right\}.
\]
$\dim\Pi_{0}=7$, co$\dim\Pi_{0}=3$. 

We observe that $(W,V)$ contains $(W,V_{\infty})$ as a codimension 1
subspace. It means that $(W,V)$ requires one less equation for its description
than $(W,V_{\infty})$. Formula (\ref{ER20comp}) describes $(W,V_{\infty})$ as
a set of one matrix and one scalar equation. It is natural to try to see if
eliminating the scalar equation results in the correct description of
$(W,V)$. Maple check confirms this hypothesis. So, the ER $(W,V)$ can be
described as
\begin{equation}
  \label{ER21}
  \BL_{11}=(\BL_{12}+\BR_{\perp})\BL_{22}^{-1}(\BL_{12}+\BR_{\perp})^{T}.
\end{equation}
The equation says that the Schur complement of $\BL_{22}$ in $\SFL-\SFT$ vanishes.
Equivalently,
\[
  \BL_{22}=(\BL_{12}+\BR_{\perp})^{T}\BL_{11}^{-1}(\BL_{12}+\BR_{\perp}).
\]
Hence, the Schur complement of $\BL_{11}$ in $\SFL-\SFT$ also vanishes.
\begin{equation}
  \label{ER21fin}
  \bb{M}=\left\{\mat{\BL}{\BL\BM-\BR_{\perp}}{\BM^{T}\BL+\BR_{\perp}}{\BM^{T}\BL\BM}:\BL>0,\
\BL+2\BR_{\perp}\BM\det\BL<0\right\}.
\end{equation}

\section{Link $(W,V)=(\bb{C}\BI,\BGY)\oplus{\rm Ann}(\bb{C}\bra{\Bz_{0}})$}
\subsection{$\Pi_{0}=(\bb{C}\BI,\BGY)$}
In this case $\SFK\in\Pi$ has the form
\[
\SFK=\mat{\psi(z)}{0}{0}{\psi(z)}+\mat{\phi(\Ga)}{\phi(\Gb)}{\phi(\Gb)}{\phi(-\Ga)}=
\mat{\BK}{\Gb\BI}{\Gb\BI}{-\cof(\BK)},\qquad\BK\in\Sym(\bb{R}^{2})
\]
\[
\SFI+\SFK=\mat{\BI_{2}+\BK}{\Gb\BI}{\Gb\BI}{\BI_{2}-\cof(\BK)},
\]
changing parametrization to $\BK'=\BI_{2}+\BK$ and dropping
primes we see that we need to compute 
\[
\SFL=2\mat{\BK}{\Gb\BI}{\Gb\BI}{\cof(2\BI_{2}-\BK)}^{-1}-\SFI
\]
Applying formula for the inverse of the block matrix we conclude that
\[
\SFL=\mat{\BL_{11}}{\Gl\BL_{11}}{\Gl\BL_{11}}{\eta(\Gl,\BL_{11})\BL_{11}}.
\]
Hence, one only needs to determine a scalar function $\eta(\Gl,\BL_{11})$. 
Applying formula for the inverse of the block matrix to $(\SFL+\SFI)^{-1}$ 
we find that the 12-block of $(\SFL+\SFI)^{-1}$ is 
\[
\BL_{12}=-\Gl(\BL_{11}+\BI_{2})^{-1}\BL_{11}(\eta\BL_{11}+\BI_{2}-\Gl^{2}\BL_{11}(\BL_{11}+\BI_{2})^{-1}\BL_{11})^{-1}.
\]
It is a multiple of the identity \IFF
\[
\eta\BL_{11}+\BL_{11}^{-1}-\Gl^{2}\BL_{11}
\]
is a multiple of the identity. In other words, we need to choose $\eta$, such
that the two eigenvalues of the above matrix are the same. If $\Gs_{1}$ and
$\Gs_{2}$ are the eigenvalues of $\BL_{11}$, then we must have
\[
\eta\Gs_{1}+\nth{\Gs_{1}}-\Gl^{2}\Gs_{1}=\eta\Gs_{2}+\nth{\Gs_{2}}-\Gl^{2}\Gs_{2},
\]
from which we find that
\[
\eta=\Gl^{2}+\nth{\Gs_{1}\Gs_{2}}=\Gl^{2}+\nth{\det\BL_{11}}.
\]
Hence,
\begin{equation}
  \label{ER9}
\bb{M}_{0}=\left\{\mat{1}{\Gl}{\Gl}{\Gl^{2}+\dfrac{1}{\det\BL}}\otimes\BL:
\Gl\in\bb{R},\ \BL>0\right\}=\{\SFL=\BGL\otimes\BL>0:
\det\SFL=\det\BGL\det\BL=1\}.
\end{equation}
The exact relation $\bb{M}_{(\bb{C}\BI,\BGY)}$ says that if the Seebeck tensor
is scalar and heat conductivity is $\BGk$ is a constant scalar multiple of
$\BGs/\det\BGs$, then the effective tensor will also have the same form. 

\subsection{Link calculation}
The strategy is use Maple to compute 
\begin{enumerate}
\item $\SFL=\SFL(\BL_{22},\BL_{12})$,
\item $\SFK=2(\SFL+\SFI)^{-1}-\SFI$,
\item Projection $\SFP$ of $\SFK$ onto $(\bb{C}\BI,\BGY)$
\item $\SFL'=2(\SFP+\SFI)^{-1}-\SFI$,
\item $\Gl$: $\Gl\BI_{2}=\BL'_{12}(\BL'_{11})^{-1}$.
\item $\eta$: $\eta\BI_{2}=\BL'_{22}(\BL'_{11})^{-1}$.
\item $\BL'_{11}$
\end{enumerate}
The link can be written as the map from 
\[
\bb{M}_{(W,V)}=\{\SFL>0:\BL_{22}=(\BL_{12}+\BR_{\perp})^{T}\BL_{11}^{-1}(\BL_{12}+\BR_{\perp})\}
\]
to
\[
\bb{M}_{(\bb{C}\BI,\BGY)}=\{\BGL\otimes\BL:\BL>0,\ \BGL>0,\ \det(\BL\BGL)=1\}.
\]
We obtain (setting $\GL_{11}=1$)
\[
\BGL=\mat{1}{\Gl}{\Gl}{\eta},\quad\Gl=\frac{\Trc\BM}{2},\quad
\eta=\det\BM,\quad\BM=\BL_{11}^{-1}(\BL_{12}+\BR_{\perp}),\quad
\BL=\BR_{\perp}\frac{\Gl\BI_{2}-\BM}{\det\BGL}.
\]
This implies that $\BM^{*}=(\BL_{11}^{*})^{-1}(\BL^{*}_{12}+\BR_{\perp})$
depends only on $\BM(\Bx)=\BL_{11}(\Bx)^{-1}(\BL_{12}(\Bx)+\BR_{\perp})$,
which can be computed from $(\BGL\otimes\BL)^{*}=\BGL^{*}\otimes\BL^{*}$ by
expressing $\BM$ in terms of $\BL$ and $\BGL$.

\section{Links for $(W,\bb{R}\BZ_{0})$}
\subsection{$\Pi_{0}=(W,\bb{R}\BZ_{0})$}
\[
\Pi_{0}=\{\SFK:2\BK_{12}=\BR_{\perp}(\BK_{11}-\BK_{22}+(\Trc\BK_{22})\BI_{2}),\
\Trc(\BK_{11})=\Trc(\BK_{22}).\}.
\]
$\dim\Pi_{0}=$co$\dim\Pi_{0}=5$. The idea is to compute the representation of
$\Pi_{0}$ from the condition that corresponding $\SFL$ satisfies
\begin{equation}
  \label{preER19}
  (\SFL^{\BR}-\SFT)(\BJ\otimes\BR_{\perp})(\SFL^{\BR}-\SFT)=0,\quad
\forall\BR\in SO(2),
\end{equation}
where
\[
\BJ=\mat{1}{0}{0}{-1},\qquad
\SFL^{\BR}=(\BR\otimes\BI_{2})\SFL(\BR^{T}\otimes\BI_{2}).
\]
Factoring $\BR\otimes\BI_{2}$ and $\BR^{T}\otimes\BI_{2}$ out and recalling
that $\BR_{\perp}$ is isotropic we obtain
\[
(\SFL-\SFT)(\BR\otimes\BI_{2})(\BJ\otimes\BR_{\perp})\BR^{T}\otimes\BI_{2}(\SFL-\SFT)=0
\]
Thus, writing $\BJ=\psi(1)$ we obtain that $(W,\bb{R}\BZ_{0})$ can be
described by the equation
\[
(\SFL-\SFT)(\psi(z)\otimes\BR_{\perp})(\SFL-\SFT)=0\quad\forall z\in\bb{C}.
\]
In other words 
\[
\begin{cases}
  (\SFL-\SFT)(\psi(1)\otimes\BR_{\perp})(\SFL-\SFT)=0,\\
  (\SFL-\SFT)(\psi(i)\otimes\BR_{\perp})(\SFL-\SFT)=0
\end{cases}
\]
The first equation is written as (\ref{ER20comp}), while the second equation
adds one more scalar condition:
\begin{equation}
  \label{ER19add}
  \Trc(\BL_{22}\cof(\BS))=0.
\end{equation}
Equivalently, this can also be written as
\[
\det(\BL_{22}+\BL_{12})=\det\BL_{22}+\det\BL_{12}.
\]
Let us write a complete system of equations for reference purposes:
\begin{equation}
  \label{ER19}
  \begin{cases}
    \BL_{11}=(\BL_{12}+\BR_{\perp})\BL_{22}^{-1}(\BL_{12}+\BR_{\perp})^{T},\\
    \det\BL_{22}=\det(\BL_{12}+\BR_{\perp}),\\
    \Trc(\BL_{22}\cof(\BL_{12}))=0
  \end{cases}
\end{equation}
The third equation can also be written as
\[
\det(\BL_{22}+\BL_{12})=\det\BL_{22}+\det\BL_{12}.
\]
In order to find a parametrization of $\SFL$ it will be convenient to write
the ER in terms of $\BM=\BL_{11}^{-1}(\BL_{12}+\BR_{\perp})$:
\begin{equation}
  \label{ER19M}
  \BL_{12}=\BL_{11}\BM-\BR_{\perp},\quad\BL_{22}=\BM^{T}\BL_{11}\BM,\quad\det\BM=1,\quad\Trc\BM=0.
\end{equation}
For example we can write
\begin{equation}
  \label{Msqm1}
  \BM=\mat{m_{11}}{m_{12}}{-\frac{m_{11}^{2}+1}{m_{12}}}{-m_{11}}
\end{equation}
Another equivalent formulation of constraints satisfied by $\BM$ is $\BM^{2}=-\BI_{2}$.
The tensor $\SFL=\SFL(\BL_{11},\BM)$ is positive definite \IFF 
\begin{equation}
  \label{WRZpos}
\BL_{11}>0,\quad  \frac{\BL_{11}}{\det\BL_{11}}+2\BR_{\perp}\BM<0.
\end{equation}
Thus we can also write
\begin{equation}
  \label{ER19fin}
\bb{M}=\left\{\mat{\BL}{\BL\BM-\BR_{\perp}}{\BM^{T}\BL+\BR_{\perp}}{\BM^{T}\BL\BM}:\BM^{2}=-\BI_{2},
\ \BL>0,\ \BL+2\BR_{\perp}\BM\det\BL<0\right\}.
\end{equation}
This equation is obtained immediately from the representations (\ref{ER20fin})
and (\ref{ER20t0fin}) for $(W,V_{\infty})$ and $(W,V_{0})$, respectively.

\subsection{Link $(W,\bb{R}\BZ_{0})$, $\Phi(K(X,Y))=K(\Ga X,Y)$}
The strategy is use Maple to compute 
\begin{enumerate}
\item $\SFL=\SFL(\BL_{22},z)$,
\item $\SFK=2(\SFL+\SFI)^{-1}-\SFI$,
\item $\SFK^{\Ga}=\Phi_{\Ga}(\SFK)$,
\item $\SFL^{\Ga}=2(\SFK_{\Ga}+\SFI)^{-1}-\SFI$,
\item $\BL_{22}^{\Ga}$ and $z^{\Ga}$.
\end{enumerate}
If we write $\SFL=\SFL(\BL_{11},\BM)$, where $\BM^{2}=-\BI_{2}$ then for any
$\Gg_{0}\in\bb{R}$ for which the resulting $\SFL$ is positive definite we have
\begin{equation}
  \label{WRZalpha}
  \SFL'=\SFL\left(\frac{\BQ\det\BL_{11}}{\det\BQ},\BM\right),\qquad
\BQ=\Gg_{0}\BL_{22}+(1+\Gg_{0})\BL_{11}+2\Gg_{0}(\det\BL_{11})\BR_{\perp}\BM.
\end{equation}
The parameters $\Ga$ and $\Gg_{0}$ are related by the formula
$\Ga=2\Gg_{0}+1$. The restrictions on $\Gg_{0}$ are
\[
\BQ>0,\qquad\frac{\BQ}{\det\BL_{11}}+2\BR_{\perp}\BM<0.
\]
An even nicer form is obtained if instead of $\BQ$ we use
$\BP=\cof(\BQ)/\det\BL_{11}$, so that
\begin{equation}
  \label{WRZalpha0}
  \SFL'=\SFL\left(\BP^{-1},\BM\right),\qquad
\BP=\Gg_{0}\BL^{-1}_{22}+(1+\Gg_{0})\BL^{-1}_{11}+2\Gg_{0}\BM\BR_{\perp},
\end{equation}
where the parameter $\Gg_{0}$ is constrained by the inequalities
\[
\BP>0,\qquad\BP+2\BM\BR_{\perp}<0,
\]
understood in the sense of quadratic forms.
We note that the map
$\Phi_{\Gg_{0}}$ fails to be bijective for $\Gg_{0}=-1/2$, but it remains a valid link between
$\bb{M}_{19}$ and ER \#18, whose definition includes the additional relation
\begin{equation}
  \label{ER18ad}
  \BM\BL^{-1}-\BL^{-1}\BM^{T}=2\BR_{\perp}
\end{equation}
between parameters $\BL$ and $\BM$.

\subsection{Link $(W,\bb{R}\BZ_{0})=(\bb{C}\BI,\bb{R}\BZ_{0})\oplus(\bb{C}(\tns{\Bz_{0}}),0)$}
\subsubsection{$\Pi_{0}=(\bb{C}\BI,\bb{R}\BZ_{0})$}
This ER is redundant, but $\Pi_{0}\oplus(\tns{\Bz_{0}},0)=(W,\bb{R}\BZ_{0})$,
where $(\tns{\Bz_{0}},0)$ is an ideal in $(W,\bb{R}\BZ_{0})$. Both the ER
$(W,\bb{R}\BZ_{0})$ and the link  corresponding to the above decomposition are
unresolved. The calculation is therefore considered useful for the puproses of
computing the unresolved cases. We have
\[
\Pi_{0}=\{K(\rho\BZ_{0},c\BI):\rho\in\bb{R},\ c\in\bb{C}\}.
\]
We will use the inversion formula
\[
\SFL=2(\SFK+\SFI)^{-1}-\SFI.
\]
\[
\SFK+\SFI=K(\BX,\BY),\qquad \BX=\mat{\rho+1}{i\rho}{-i\rho}{\rho+1},\quad\BY=c\BI_{2}.
\]
We will use formula (\ref{Kinv}):
\[
(\SFK+\SFI)^{-1}=K(S_{X}^{1},S_{Y}^{-1}),\quad
S_{X}=\BX-\BY\bra{\BX}^{-1}\bra{\BY},\quad S_{Y}=\bra{\BY}-\bra{\BX}\BY^{-1}\BX.
\]
We compute
\[
\bra{\BX}^{-1}=\frac{\BX}{2\rho+1},\quad\bra{\BX}\BX=(2\rho+1)\BI_{2}.
\]
Using these formulas we compute
\[
S_{X}=\frac{2\rho+1-|c|^{2}}{2\rho+1}\BX,\quad
S_{Y}=\frac{|c|^{2}-2\rho-1}{c}\BI_{2}.
\]
Hence,
\[
S_{X}^{-1}=\frac{\bra{\BX}}{2\rho+1-|c|^{2}},\quad
S_{Y}^{-1}=\frac{c}{|c|^{2}-2\rho-1}\BI_{2}.
\]
Writing $K(\bra{\BX},0)=(\rho+1)\SFI+\rho\SFT$ we compute
\[
\SFL=\frac{(1+|c|^{2})\SFI+2\rho\SFT-K(0,2c\BI_{2})}{2\rho+1-|c|^{2}}.
\]
Converting to the block matrix form we get
\[
\SFL=\mat{\BL}{-\Gth\BR_{\perp}}{\Gth\BR_{\perp}}{\BL},\quad
\BL=\frac{\Gvf(1+|c|^{2})-\psi(2c)}{2\rho+1-|c|^{2}},\quad
\Gth=\frac{2\rho}{2\rho+1-|c|^{2}}.
\]
Now it is easy to see that $\det\BL=(1-\Gth)^{2}$. Requiring that $\BL>0$ we
get the restriction that $2\rho+1-|c|^{2}>0$. If this is satisfied then
$\SFL>0$ holds \IFF $\Gth^{2}<\det\BL$. The two constraints can be combined
into one: 
\[
2|\rho|<1-|c|^{2},
\]
giving the ER
\begin{equation}
  \label{ER7}
  \bb{M}=\left\{\mat{\BL}{\Gth\BR_{\perp}}{-\Gth\BR_{\perp}}{\BL}: \det\BL=(1+\Gth)^{2},\
\BL>0,\ \Gth>-1/2\right\}
\end{equation}
Of course the isomorphic ER $(\bb{C}\BI,\bb{R}\bra{\BZ_{0}})$ is obtained from
this by replacing $\Gth$ with $-\Gth$ in the matrix, while keeping all other
constraints the same.

\subsubsection{Link calculation}
The strategy is use Maple to compute 
\begin{enumerate}
\item $\SFL=\SFL(\BL_{11},\BM)$, where $\BM$ is given by (\ref{Msqm1})
\item $\SFK=2(\SFL+\SFI)^{-1}-\SFI$,
\item Projection $\SFP$ of $\SFK$ onto $(\bb{C}\BI,\bb{R}\BZ_{0})$
\item $\SFL'=2(\SFP+\SFI)^{-1}-\SFI$,
\item $\Gth$: $\Gth\BI_{2}=\BR_{\perp}^{T}\BL'_{12}$.
\item $\BL'_{11}$
\end{enumerate}
We obtain
\[
\Gth+1=-\frac{2\det\BL_{11}}{\Trc(\BL_{11}\BM\BR_{\perp})}=
-\frac{2}{\Trc(\BL_{11}^{-1}\BR_{\perp}\BM)}
\]
\[
\BL'_{11}=-(\Gth+1)\BR_{\perp}\BM=\frac{2\BR_{\perp}\BM}{\Trc(\BL_{11}^{-1}\BR_{\perp}\BM)}.
\]
If we use the link $(W,\bb{R}\BZ_{0})/{\rm
  Ann}(\bb{C}\bra{\Bz_{0}})\cong(\bb{C}\BI,0)$, which is inherited from $(W,V)$,
then we obtain that $\BGs^{*}=-\BR_{\perp}\BM^{*}$ is the effective
conductivity of the 2D conducting composite with local conductivity
$\BGs(\Bx)=-\BR_{\perp}\BM(\Bx)$, also satisfying $\det\BGs=1$. This link
shows that the scalar $s^{*}=\av{\BL_{11}^{*}/\det\BL^{*}_{11},\BGs^{*}}$ depends
only on $\BGs(\Bx)$ and $s(x)$ via the effective thermoelectricity in the ER (\ref{ER7}).

\subsection{Link $(W,\bb{R}\BZ_{0})=(W,0)\oplus(0,\bb{R}\BZ_{0})$}
\subsubsection{$(W,0)$}
Need to verify that one additional equation that must be added to the system
(\ref{ER19M}) is 
\[
\BL_{11}\BM-\BM^{T}\BL_{11}=2\det\BL_{11}\BR_{\perp}.
\]
If we write $\BL_{12}=\BS+\Gb\BR_{\perp}$, then this additional equation says
that $\Gb+1=\det\BL_{11}$. Recalling the condition of positivity
(\ref{WRZpos}) of $\SFL$ we can also write the extra equation as
\[
\BL_{22}=-(\BL_{11}+2\BR_{\perp}\BM\det\BL_{11})>0.
\]

\subsection{Link calculation}
Same strategy using Maple produces a link
$\Phi(\SFL(\BL_{11},\BM))=\SFL(\BL'_{11},\BM)$, where $\BL'_{11}$ satisfies
the additional relation
\[
\BL'_{11}\BM-\BM^{T}\BL'_{11}=2\det\BL'_{11}\BR_{\perp}.
\]
We compute using Maple:
\[
\SFL'=\SFL\left(\frac{\BQ\det\BL_{11}}{\det\BQ},\BM\right),\qquad
\BQ=\frac{\BL_{11}-\BL_{22}}{2}-(\det\BL_{11})\BR_{\perp}\BM,
\]
which coincides with formula (\ref{WRZalpha}) when $\Gg_{0}=-1/2$. Indeed,
$\Gg_{0}=-1/2$ corresponds to $\Ga=0$, so that the map $\Phi(K(X,Y))=K(\Ga
X,Y)$ is no longer the automorphism, but the link
$(W,\bb{R}\BZ_{0})=(W,0)\oplus(0,\bb{R}\BZ_{0})$ instead.

\section{Summary of the essential exact relations and links}
Here we refer to various exact relations by their number in the list at the
end of Section~\ref{sec:JMA}. In order to streamline our notation it will be
convenient to introduce the function
\begin{equation}
  \label{LMpar}
  \mathfrak{L}(\BL,\BM)=\mat{\BL}{\BL\BM}{\BM^{T}\BL}{\BM^{T}\BL\BM}+\SFT,
\qquad\SFT=\BR_{\perp}\otimes\BR_{\perp},
\end{equation}
since many of the exact relations below can be described in terms of
$\mathfrak{L}(\BL,\BM)$. Here is the list.

%\begin{itemize}
%\item%[\#8] $\Pi_{0}=(\bb{C}\BI_{2},\BGF)$
  \[
    \bb{M}_{8}=\{\BI_{2}\otimes\BL+t\SFT:\BL\in\Sym^{+}(\bb{R}^{2}),\ \det\BL>t^{2}\},
  \]
%\item%[\#13] $\Pi_{0}={\rm Ann}(\bra{\Bz_{0}})$, $\Pi_{0}^{2}=\{0\}$.
\[
  \bb{M}_{13}=\left\{\mathfrak{L}(\BL,\BR_{\perp}):\BL>\hf\BI_{2} \right\},\qquad
\BL_{*}=\av{\BL^{-1}}^{-1}.
\]
%\item%[\#17] $\Pi_{0}=(\ClD,\ClD')$.
\[
  \bb{M}_{17}=\{\SFL>0:\SFL(\BJ\otimes\BR_{\perp})\SFL=\BJ\otimes\BR_{\perp}\},\quad\BJ=\mat{1}{0}{0}{-1}.
\]
In block-components we can rewrite this as
\[
  \bb{M}_{17}=\left\{\SFL>0:\BL_{11}=\frac{\BL_{12}\cof(\BL_{22})\BL_{12}^{T}}{\det\BL_{12}},\
\det\BL_{22}-\det\BL_{12}=1\right\}.
\]
Of course, there is symmetry between indices and we also have
\[
\bb{M}_{17}=\left\{\SFL>0:\BL_{22}=\frac{\BL_{12}^{T}\cof(\BL_{11})\BL_{12}}{\det\BL_{12}},\
\det\BL_{11}-\det\BL_{12}=1\right\}.
\]
%\item
\[
\bb{M}_{19}=\left\{\mathfrak{L}(\BL,\BM):\BM^{2}=-\BI_{2}, \ \BL>0,\ 
\BL^{-1}+2\BM\BR_{\perp}<0\right\}.
\]
This exact relation has two different links, which are not consequences of
other relations or links listed here.
\begin{enumerate}
\item This is an infinite family of links that we describe in terms of the
  function $\mathfrak{L}(\BL,\BM)$, given by (\ref{LMpar}). The family of links
  are the maps $\Phi_{\Gg_{0}}:\bb{M}_{19}\to\bb{M}_{19}$, given by
\[
  \Phi_{\Gg_{0}}(\mathfrak{L}(\BL,\BM))=\mathfrak{L}\left(\BP_{\Gg_{0}}^{-1},\BM\right),
\quad\BP_{\Gg_{0}}=\Gg_{0}\BM\BL^{-1}\BM^{T}+(1+\Gg_{0})\BL^{-1}+2\Gg_{0}\BM\BR_{\perp},
\]
where the parameter $\Gg_{0}$ is constrained by the inequalities
\[
\BP_{\Gg_{0}}>0,\qquad\BP_{\Gg_{0}}+2\BM\BR_{\perp}<0,
\]
understood in the sense of quadratic forms. We note that the map
$\Phi_{\Gg_{0}}$ fails to be bijective for $\Gg_{0}=-1/2$, but it remains a valid link between
$\bb{M}_{19}$ and
\[
\bb{M}_{18}=\{\mathfrak{L}(\BL,\BM):\BM^{2}=-\BI_{2}, \ 
\BM\BL^{-1}-\BL^{-1}\BM^{T}=2\BR_{\perp},\ \BL>0\}.
\]
\item The second link is between $\bb{M}_{19}$ and
  \[
    \bb{M}_{7}=\{\mathfrak{L}(\mu\BGs,\BR_{\perp}\BGs):\det\BGs=1,\ \BGs>0,\ \mu>1/2\}.
  \]
The link is then given by the formulas
\[
  \BGs=-\BR_{\perp}\BM,\qquad\mu=\frac{2}{\Trc(\BL\BGs)},
\]
so that $\BM^{*}=\BR_{\perp}\BGs^{*}$, where $\BGs^{*}$ is the
effective conductivity of the 2D polycrystal with texture
$\BGs(\Bx)=-\BR_{\perp}\BM(\Bx)$, as before, while additionally we have
\[
  \Trc(\BL^{*}\BGs^{*})=\frac{2}{\mu^{*}}.
\]
\end{enumerate}
%\item%[\#20] $\Pi_{0}=(W,V_{\infty})$
\[
  \bb{M}_{20}=\left\{\SFL>0:(\SFL-\SFT)(\BJ\otimes\BR_{\perp})(\SFL-\SFT)=0\right\}.
\]
We can also write this ER in parametric form 
\[
  \bb{M}_{20}=\left\{\mathfrak{L}(\BL,\BM):\det\BM=1,\ \BL>0,\ 
\BL^{-1}<-2\BM\BR_{\perp}\right\},
\]
where inequalities are undestood in the sense of quadratic forms.
%\item%[\#21] $\Pi_{0}=(W,V)$.
\[
  \bb{M}_{21}=\left\{\mathfrak{L}(\BL,\BM):\BL>0,\ 
\BL^{-1}<-2\BM\BR_{\perp}\right\}.
\]
There is also a link associated with this ER. It says that
$\BM^{*}$ does not depend on $\BL(\Bx)$ in the parametrization
(\ref{ER21fin}). The effective tensor $\BM^{*}$ can be computed from the exact
relation described by
\[
  \bb{M}_{9}=\{\SFL=\BGL\otimes\BP>0:
\det\SFL=\det\BGL\det\BP=1\}.
\]
Specifically, $\SFL=\BGL\otimes\BP\in\bb{M}_{9}$ is uniquely determined by a pair of symmetric,
positive definite $2\times 2$ matrices
$\BGL$ and $\BP$, satisfying $\det\BGL\det\BP=1$, provided we fix
$\GL_{11}=1$. We will denote this parametrization by
$\SFL=\SFL_{9}(\BGL,\BP)$. The fact that $\bb{M}_{9}$ is an exact relation
means that $\SFL_{9}(\BGL,\BP)^{*}=\SFL_{9}(\BGL^{*},\BP^{*})$, for some
$\BGL^{*}$ and $\BP^{*}$ that depend on the \mc\ of the composite.
The link between $\bb{M}_{21}$, given by (\ref{ER21fin}) and
$\bb{M}_{9}$ is given by a bijective transformation
$\BM\mapsto(\BGL(\BM),\BP(\BM))$
\[
  \BGL(\BM)=\mat{1}{\Trc\BM/2}{\Trc\BM/2}{\det\BM},\qquad
  \BP(\BM)=-\BR_{\perp}\frac{\BM-(\Trc\BM)\BI_{2}/2}{\det\BGL(\BM)}.
\]
The link says that $\BM^{*}$ is determined via the formula
\[
  \SFL_{9}(\BGL(\BM),\BP(\BM))^{*}=\SFL_{9}(\BGL(\BM^{*}),\BP(\BM^{*})),
\]
so, that
\[
  \BM^{*}=\GL_{12}^{*}\BI_{2}+\BR_{\perp}\BP^{*}\det\BGL^{*}.
\]
%\item%[\#22] $\Pi_{0}=(\Sym(\bb{C}^{2}),0)$
\[
  \bb{M}_{22}=\{\SFL>0:\SFL(\BI_{2}\otimes\BR_{\perp})\SFL=\BI_{2}\otimes\BR_{\perp}\}.
\]
In block-components we can rewrite this as
\[
  \bb{M}_{22}=\left\{\SFL>0:\BL_{11}=-\frac{\BL_{12}\cof(\BL_{22})\BL_{12}^{T}}{\det\BL_{12}},\
\det\BL_{22}+\det\BL_{12}=1\right\}.
\]
Of course, there is symmetry between indices and we also have
\[
\bb{M}_{22}=\left\{\SFL>0:\BL_{22}=-\frac{\BL^{T}_{12}\cof(\BL_{11})\BL_{12}}{\det\BL_{12}},\
\det\BL_{11}+\det\BL_{12}=1\right\}.
\]
%\end{itemize}

\section{An application: Isotropic polycrystals}
The subspace $\Pi_{0}=(\Sym(\bb{C}^{2}),0)$ contains the unique isotropic
tensor $\SFK=0$. Therefore, the corresponding exact relation
(\ref{ER22}) passes through the unique isotropic tensor $\SFL=\SFI$. The
global automorphisms 
\[
\Psi_{\Ga,\BB}(\SFL)=(\BB\otimes\BI_{2})(\SFL-\Ga\SFT)(\BB\otimes\BI_{2}),\quad\BB^{T}=\BB.
\]
Map $\SFL=\SFI$ into $\BB^{2}\otimes\BI_{2}-(\Ga\det\BB)\SFT$. Every isotropic
tensor has the form $\BGL\otimes\BI_{2}+\Gb\SFT$, where $\BGL$ is symmetric
and positive definite (additionally we must also have
$\det\BGL>\Gb^{2}$). Thus, there exists a unique symmetric and positive
definite matrix $\BB$, such that $\BB^{2}=\BGL$ and $\Ga=\Gb/\det\BB$, so that
$\Psi_{\Ga,\BB}(\SFI)=\BGL\otimes\BI_{2}+\Gb\SFT$. We conclude that for every
isotropic, symmetric and positive definite tensor $\SFL$ there exists a unique
symmetric and positive definite real $2\times 2$ matrix $\BB$ and real number
$|\Ga|<1$, such that $\SFL=\Psi_{\Ga,\BB}(\SFI)$. Each such transformation
maps the exact relation (\ref{ER22}) into another exact relation, which are
all disjoint and therefore foliate an open \nbh\ of the space of isotropic
tensors. Suppose $\SFL_{0}$ is an anisotropic tensor. Then for sufficiently
small $\Ge>0$ the tensor $\Ge\SFL_{0}$ will be in that \nbh\ foliated by exact
relations. Hence, there exists a unique exact relation $\bb{M}_{\Ge}$
isomorphic to (\ref{ER22}) that passes through $\Ge\SFL_{0}$. But then
$\Ge^{-1}\bb{M}_{\Ge}$ is the exact relation isomorphic to (\ref{ER22})
that passes through $\SFL_{0}$. Thus, regardless of texture, the effective
tensor of an isotropic polycrystal made of the single cystallite $\SFL_{0}$
will be uniquely determined by $\SFL_{0}$, as in 2D conductivity
$\Gs^{*}=\sqrt{\det\BGs_{0}}$. If its effective tensor is
$\SFL^{*}$, then if $\Psi_{\Ga,\BB}(\SFL^{*})=\SFI$, then
$\SFL'=\Psi_{\Ga,\BB}(\SFL_{0})$ will belong to the exact relation
(\ref{ER22}). Specifically,
\[
2(\SFL'+\SFI)^{-1}-\SFI\in\Pi_{0}=(\Sym(\bb{C}^{2}),0).
\] 
In order to find equations satisfied by $\BB$ and $\Ga$ we will write
$\SFL_{0}=K(\BX,\BY)$. Then
\[
\SFL'=\Psi_{\Ga,\BB}(\SFL_{0})=K(\BB(\BX-i\Ga\BR_{\perp})\BB,\BB\BY\BB).
\]
Thus, we need to find $\BB\in\Sym^{+}(\bb{R}^{2})$ and $\Ga\in(-1,1)$, such that
\[
2K(\BB(\BX-i\Ga\BR_{\perp})\BB+\BI_{2},\BB\BY\BB)^{-1}-K(\BI_{2},0)=K(\Bzr,\cdot).
\]
Using formula (\ref{Kinv}) we obtain $S_{X}=2\BI_{2}$,
where
\[
S_{X}=\BB(\BX-i\Ga\BR_{\perp})\BB+\BI_{2}-\BB\BY\BB(\BB(\bra{\BX}+i\Ga\BR_{\perp})\BB+\BI_{2})^{-1}\BB\bra{\BY}\BB.
\]
Hence, the equation $S_{X}=2\BI_{2}$ becomes
\[
\BX-i\Ga\BR_{\perp}-\BY(\bra{\BX}+i\Ga\BR_{\perp}+\BB^{-2})^{-1}\bra{\BY}=\BB^{-2}
\]
We now observe that $\SFL^{*}=K(\BB^{-2}+i\Ga\BR_{\perp},0)$. Thus, denoting 
\[
\BL^{*}=\BB^{-2}+i\Ga\BR_{\perp}\in\mathfrak{H}(\bb{C}^{2}),
\]
we obtain the equation for $\BL^{*}$:
\begin{equation}
  \label{isopoly}
  \BX-\BL^{*}=\BY(\bra{\BX}+\BL^{*})^{-1}\bra{\BY},\quad\SFL_{0}=K(\BX,\BY),\quad
\SFL^{*}=K(\BL^{*},0).
\end{equation}
% We can also restate this equation as requirement that both Schur complements
% $\BS_{11}$ and $\BS_{22}$ of the block matrix
% \[
% \mat{\BX-\BL^{*}}{\BY}{\bra{\BY}}{\bra{\BX}+\BL^{*}}
% \]
% are zero.
If we make a change of variables $\BZ=\bra{\BX}+\BL^{*}$ then $\BZ$ is still
self-adjoint (and positive definite) and solves
\begin{equation}
  \label{Zeq}
  \BZ+\BY\BZ^{-1}\BY^{H}=\BX+\bra{\BX}.
\end{equation}
We can first solve 4 real linear equations with 4 real unknowns:
\begin{equation}
  \label{Zeqlin}
  \BZ+\Gth\BY\cof(\BZ)^{T}\BY^{H}=\BX+\bra{\BX},
\end{equation}
Obtaining a solution $\Hat{\BZ}(\Gth)$. We then find $\Gth>0$ from the
equation $\Gth\det\Hat{\BZ}(\Gth)=1$. Let us analyze the linear equation
(\ref{Zeqlin}), assuming first that $\BY$ is invertible. Let
$\mathfrak{B}_{\BY}\in\End_{\bb{R}}(\mathfrak{H}(\bb{C}^{2}))$ be defined by
\[
\mathfrak{B}_{\BY}\BZ=\BY\cof(\BZ)^{T}\BY^{H}.
\]
Let $\Gl$ be an eigenvalue of $\mathfrak{B}_{\BY}$.
Then, taking determinants in the equation $\mathfrak{B}_{\BY}\BZ=\Gl\BZ$ we obtain
\[
\Gl^{2}\det\BZ=|\det\BY|^{2}\det\BZ.
\]
Hence, either $\Gl=\pm|\det\BY|$ or $\det\BZ=0$. In the latter case, $\BZ$
is a real multiple of $\Ba\otimes\bra{\Ba}$ for some nonzero vector
$\Ba\in\bb{C}^{2}$. Then
\[
\BY\BR_{\perp}\bra{\Ba}\otimes\bra{\BY}\BR_{\perp}\Ba=\Gl\Ba\otimes\bra{\Ba}.
\]
Taking traces we obtain
\[
\Gl=\frac{|\BY\BR_{\perp}\bra{\Ba}|^{2}}{|\Ba|^{2}}\ge 0.
\]
Thus, the only possible negative eigenvalue of $\mathfrak{B}_{\BY}$ is
$\Gl=-|\det\BY|$. Observing that
$\mathfrak{B}_{e^{i\Ga}\BY}=\mathfrak{B}_{\BY}$ for any $\Ga\in\bb{R}$, we may
assume, \WLOG, that $\det\BY>0$. It is then easy to see that 
\[
\mathfrak{B}_{\BY}\re(\BY)=(\det\BY)\re(\BY),\qquad
\mathfrak{B}_{\BY}\im(\BY)=-(\det\BY)\im(\BY).
\]
If $\BY$ is real and symmetric, then 
\[
\mathfrak{B}_{\BY}\BZ_{\pm}(c)=(\pm\det\BY)\BZ_{\pm}(c),\quad
\BZ_{+}=\phi(c)\BY+\BY\phi(\bra{c}),\quad\BZ_{-}=\psi(c)\BY+\BY\psi(c)
\]
for any $c\in\bb{C}$ for which $\BZ_{\pm}(c)\not=0$.
In fact, the characteristic polynomial of $\mathfrak{B}_{\BY}$ (as computed by Maple) is 
\[
p(x)=(x^{2}-|\det\BY|^{2})(x^{2}+|\det\BY|^{2}+x\av{\BY,\cof(\BY)})),\quad
\av{\BA,\BB}=\Trc(\BA\BB^{H}).
\]
One can check that the roots of $p(x)$, other than $\pm|\det\BY|$, are either
both complex or both positive.

Hence, $\Hat{\BZ}(\Gth)=(\SFI+\Gth\mathfrak{B}_{\BY})^{-1}(\BX+\bra{\BX})$ and
$\Gth$ is found as a positive root of 
\begin{equation}
  \label{rootheta}
  \Gth\det\left[(\SFI+\Gth\mathfrak{B}_{\BY})^{-1}(\BX+\bra{\BX})\right]=1,
\end{equation}
such that $(\SFI+\Gth\mathfrak{B}_{\BY})^{-1}(\BX+\bra{\BX})>\bra{\BX}$.
% Observe that $\mathfrak{B}_{\BY}$ maps positive definite matrices into
% positive definite ones. Therefore, if $\Gth_{1}>\Gth_{2}>0$, then, since $\BX+\bra{\BX}>0$, 
% \[
% (\SFI+\Gth_{1}\mathfrak{B}_{\BY})^{-1}(\BX+\bra{\BX})<(\SFI+\Gth_{2}\mathfrak{B}_{\BY})^{-1}(\BX+\bra{\BX}).
% \]
% Hence, $\Gth$ must be the smallest positive root of (\ref{rootheta}), since we
% know that $\BL^{*}$ satisfying all the requirements exists and is unique.

The case $\det\BY=0$ was not considered, but being the
limiting case of the general one, our conclusion stays the same: $\Gth$ is
a positive root of (\ref{rootheta}), such that
\begin{equation}
  \label{Lst}
  \BL^{*}=(\SFI+\Gth\mathfrak{B}_{\BY})^{-1}(\BX+\bra{\BX})-\bra{\BX}>0.
\end{equation}
We conjecture that $\Gth$ must be the smallest positive root of
(\ref{rootheta}). This is easily verified if $|\BY|$ is sufficiently small.

In the special case when $\BY$ is real (or purely imaginary) we can write reasonably compact
equations:
\[
\BZ=\frac{2\re(\BX)}{1-\Gth\det\BY}-\frac{2\Gth\Trc(\cof(\BY)\re(\BX))\BY}{1-\Gth^{2}\det\BY^{2}}
\]
\[
\frac{\Gth\det\re(\BX)}{(1-\Gth\det\BY)^{2}}-\frac{\Gth^{2}(\Trc(\cof(\BY)\re(\BX)))^{2}}
{(1-\Gth^{2}\det\BY^{2})^{2}}=\nth{4}.
\]
If we change variables $t=\Gth\det\BY$, then we can write the equation for $t$ as
\[
t(1+t)^{2}\det(\BY^{-1}\re(\BX))-t^{2}(\Trc(\BY^{-1}\re(\BX)))^{2}=\nth{4}(1-t^{2})^{2}.
\]
If $s_{1}$ and $s_{2}$ are the eigenvalues of
$(\re\BX)^{1/2}\BY^{-1}(\re\BX)^{1/2}$, then (assuming, \WLOG, that
$\det\BY>0$) we obtain $|s_{j}|>1$ as a consequence of positive definteness of
$\SFL_{0}$, while the equation for $t$ reads
\[
p(t)=t(1+t)^{2}|s_{1}||s_{2}|-t^{2}(|s_{1}|+|s_{2}|)^{2}-\nth{4}(1-t^{2})^{2}=0.
\]
The case $s_{1}=s_{2}=s$ can be solved explicitly. The roots of $p(t)$ are $t=1$
of multiplicity 2 and $t_{\pm}(s)=2s^{2}-1\pm2s\sqrt{s^{2}-1}$. It is obvious
that $p(t)<0$, when $t\le 0$, so all the real roots of $p(t)$ have to be
positive. The product $t_{+}(s)t_{-}(s)=1$, so one of the roots is in $(0,1)$,
while the other is in $(1,+\infty)$. The discriminant of $p(t)$ is
\[
\GD[p]=(s_{1}^{2}-1)^{2}(s_{2}^{2}-1)^{2}(s_{1}^{2}-s_{2}^{2})^{2}.
\]
Thus, for all $(s_{1},s_{2})\in D=\{(s_{1},s_{2}):s_{1}>s_{2}>1\}$ the number
of real roots is the same. We also compute $p(0)=-1/2$ and
$p(1)=-(s_{1}-s_{2})^{2}$, which show that the number of roots in $(0,1)$ and
in $(1,+\infty)$ remains the same. It is easy to check that there are 4 real
roots when $s_{1}=s_{2}+\Ge$: two in $(0,1)$ and two in $(1,+\infty)$.

\section{An application: two phase composites with isotropic phases}
\label{sec:twophase}
The results in this section consitute Master's thesis of Sarah Childs (MS
2020, Temple University).

\subsection{Analysis}
\begin{figure}[t]
  \centering
  \includegraphics[scale=0.39]{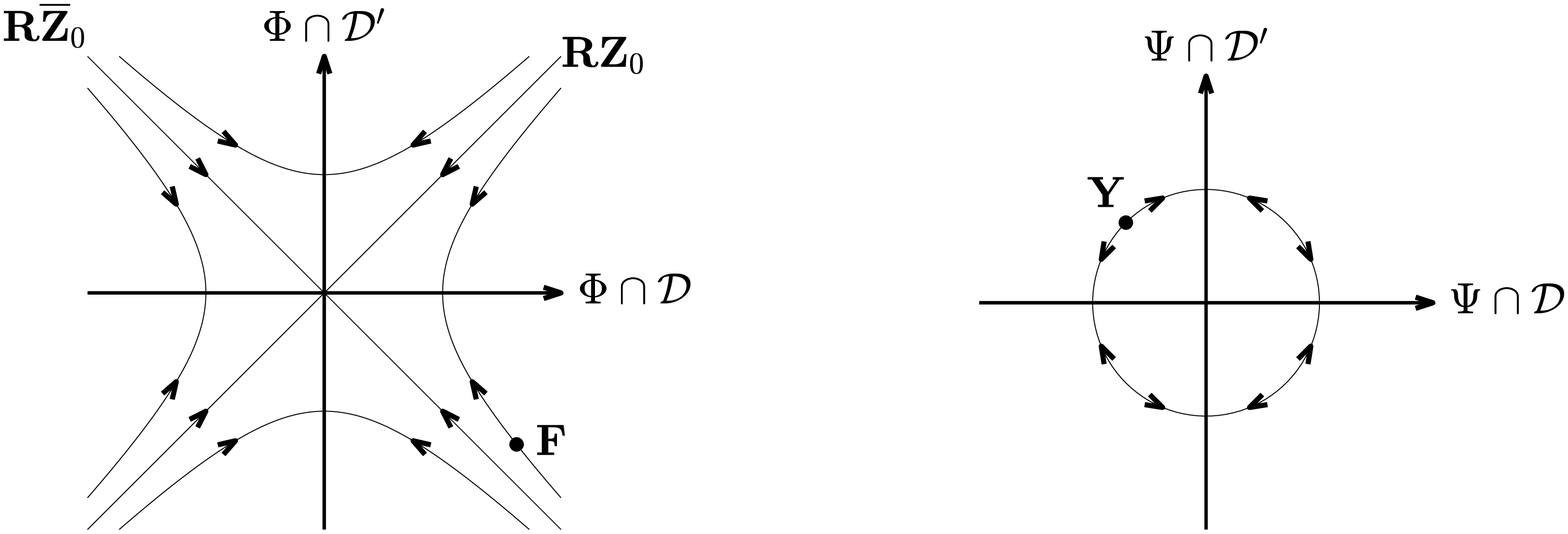}
  \caption{Action of the global automorphism (\ref{Jautoex}) on $\BX_{0}=\BF+\BY$.}
  \label{fig:bincomp}
\end{figure}

If we have two given isotropic tensors $\SFL_{1}$, $\SFL_{2}$ we first use the
global link to map $\SFL_{1}$ to $\SFI$, while $\SFL_{2}$ will be mapped to
some other isotropic tensor $\SFL_{0}$. Next we compute
$\SFK_{0}=2(\SFL_{0}+\SFI)^{-1}-\SFI$ and write is as $\SFK_{0}=K(\BX_{0},0)$ for
some $\BX_{0}\in\mathfrak{H}(\bb{C}^{2})$. We then apply the global automorphism 
(\ref{Jautoex}) and map $\BX_{0}$ to $\BC\BX_{0}\BC^{H}$ using $\BC\in O(2,\bb{C})$. Hence,
we need to understand the action of $G=O(2,\bb{C})$ on
$\mathfrak{H}(\bb{C}^{2})$. The group $G$ has two connected components
$G_{+}=SO(2,\bb{C})$ and $G_{-}=\psi(1)G_{+}$. Hence, it is enough to
understand the action of the subgroup $G_{+}$ on $\mathfrak{H}(\bb{C}^{2})$.
We do this by first identifying invariant subspaces
of $G_{+}$. The calculations have already been done. These invariant subspaces are
$\bb{R}\BZ_{0}$, $\bb{R}\bra{\BZ_{0}}$, and $\BGY$, where 
\[
\mathfrak{H}(\bb{C}^{2})=\bb{R}\BZ_{0}\oplus\bb{R}\bra{\BZ_{0}}\oplus\BGY.
\]
We recall that 
$\BGY=\{\psi(z):z\in\bb{C}\}$. We then have for $\BC_{+}(c)\in G_{+}$,
given by (\ref{Cpmdef})
\[
\BC_{+}(c)\psi(z)\BC_{+}(c)^{H}=\psi(e^{-2i\re(c)}z),
\]
and
\[
\BC_{+}(c)\BZ_{0}\BC_{+}(c)^{H}=e^{2\im(c)}\BZ_{0},\qquad
\BC_{+}(c)\bra{\BZ_{0}}\BC_{+}(c)^{H}=e^{-2\im(c)}\bra{\BZ_{0}}.
\]
These formulas show that $G_{+}$ contains two 1-parameter subgroups 
\[
H_{+}=\{\BC_{+}(c): c\in\bb{R}\}=SO(2,\bb{R}),\qquad
H_{-}=\{\BC_{+}(c): \re(c)=0\}.
\]
All points in the subspace $\BGY$ are fixed by $H_{-}$, while all points
in the subspace $\BGF=\Span_{\bb{R}}\{\BZ_{0},\bra{\BZ_{0}}\}$ are fixed by $H_{+}$. At
the same time $H_{+}$ acts by rotations on $\BGY$, while $H_{-}$ acts by
hyperbolic rotations on $\BGF$. Thus, in order to understand how we can
transform $\BX_{0}$ by $G_{+}$ we first split $\BX_{0}$ into its $\BGF$ and
$\BGY$ components: 
\[
\BX_{0}=\BF+\BY,\qquad\BF\in\BGF,\quad \BY\in\BGY.
\] 
We first apply $H_{+}$ to transform $\BY$ as desired, while $\BF$ is
unchanged. We then apply $H_{-}$ to transform $\BF$ as desired, while the
transformed matrix $\BY$ is unchanged. What can be accomplished is shown in
Fig.~\ref{fig:bincomp}. 
It depends on
whether we are in a generic situation, where $\BY\not=0$ and $\BF$ is neither
a real multiple of $\BZ_{0}$ nor of $\bra{\BZ_{0}}$, or in one of the special
ones. In the generic case we can rotate $\BF$ to a matrix
$f\BI_{2}\in\ClD\cap\BGF$ or $if\BR_{\perp}\in\ClD'\cap\BGF$, depending on
whether $|\Trc\BF|$ is greater or smaller than $2|\im(F_{12})|$,
respectively. While we can rotate $\BY$ to $\psi(y)\in\ClD\cap\BGY$ or in
$\psi(iy)\in\ClD'\cap\BGY$, $y>0$, as desired, i.e. according to which space the
component $\BF$ can be rotated to. In that case either the physically trivial
exact relation $(\ClD,\ClD)$ can be used or an exact relation $(\ClD,\ClD')$
is applicable, unless $\det(\BX_{0})=0$, in which case
$\Ann(\bb{C}\Be_{2})\subset(\ClD,\ClD)$ exact relation would apply. There are the
following special cases
\begin{enumerate}
\item $\BY\not=0$, $\BF\in\bb{R}\BZ_{0}$ or $\BF\in\bb{R}\bra{\BZ_{0}}$ and
  $\BF\not=0$. Exact relation $(W,V_{\infty})\sim(\bra{W},\bra{V_{\infty}})$
  is applicable.
\item $\BY\not=0$, $\BF=0$. Exact relation
  $(\bb{C}\BI,\bb{R}\BY)\sim(\bb{C}\BI,\bb{R}\psi(i))=(\bb{C}\BI,\Psi)\cap(\ClD,\ClD')$
  is applicable.
\item $\BY=0$, $\BF\not\in\bb{R}\BZ_{0}$,
  $\BF\not\in\bb{R}\bra{\BZ_{0}}$. Exact relation $(\bb{C}\BI,\bb{R}\BF)$ is
  applicable. If we are in a strongly coupled case then
  $(\bb{C}\BI,\bb{R}\BF)\sim(\bb{C}\BI,i\BR_{\perp})=(\bb{C}\BI,\Psi)\cap(\ClD,\ClD')$,
  otherwise,
  $(\bb{C}\BI,\bb{R}\BF)\sim(\bb{C}\BI,\bb{R}\BI)=(\bb{C}\BI,\Psi)\cap(\ClD,\ClD)$. 
\item $\BY=0$ and $\BF\in\bb{R}\BZ_{0}$ or $\BF\in\bb{R}\bra{\BZ_{0}}$ and $\BF\not=0$. Exact
  relation $(0,\bb{R}\BZ_{0})\sim(0,\bb{R}\bra{\BZ_{0}})$ is applicable.
\end{enumerate}
Thus, if both components of $\BX_{0}$ are non-zero then there are 4 exact
relations between components of $\SFL^{*}$. What is interesting is that the
form of these relations depends very much on specific values of the
components.  If one of the components of $\BX_{0}$ happens to be 0, then there
would be at least 7 exact relations, rising to 9 in a very special case
$\BX_{0}=x_{0}\BZ_{0}$ or $\BX_{0}=x_{0}\bra{\BZ_{0}}$, in this very special case
\emph{all} components of $\SFL^{*}$ will be uniquely determined if the volume
fractions of the components are known.

If we write $\SFL_{j}=\BGs_{j}\otimes\BI_{2}+r_{j}\SFT$. Then the
transformation
\begin{equation}
  \label{Psi0}
  \Psi(\SFL)=(\BGs_{1}^{-1/2}\otimes\BI_{2})(\SFL-r_{1}\SFT)(\BGs_{1}^{-1/2}\otimes\BI_{2})
\end{equation}
maps $\SFL_{1}$ to $\SFI$, while
\[
\SFL'=\Psi(\SFL_{2})=K\left(\BGs_{1}^{-1/2}\BGs_{2}\BGs_{1}^{-1/2}+\frac{r_{2}-r_{1}}{\sqrt{\det{\BGs_{1}}}}i\BR_{\perp},0\right).
\]
Next we apply transformation $\Psi_{\BA,\BI_{2}}$, where 
\[
 \BA=\mat{a_{0}}{1}{1}{a_{0}}
\]
These transformations have the property that
$\Psi(\SFI)=\SFI$. Then denoting
\begin{equation}
  \label{sigmarho}
  \BGs=\BGs_{1}^{-1/2}\BGs_{2}\BGs_{1}^{-1/2},\qquad\rho=\frac{r_{2}-r_{1}}{\sqrt{\det{\BGs_{1}}}},
\end{equation}
we obtain
\[
\Psi_{\BA,\BI_{2}}(K(\BGs+i\rho\BR_{\perp},0))=K(\BL,0),\quad
\BL=i\BR_{\perp}(\BGs+i(\rho+a_{0})\BR_{\perp})^{-1}(a_{0}\BGs+i(a_{0}\rho+1)\BR_{\perp}).
\]
Since all the matrices involved in the formula for $\BL$ are $2\times 2$ we
have
\[
(\BGs+i(\rho+a_{0})\BR_{\perp})^{-1}=
\frac{\BR_{\perp}^{T}(\BGs-i(\rho+a_{0})\BR_{\perp})\BR_{\perp}}{\det\BGs-(\rho+a_{0})^{2}}.
\]
Thus,
\begin{equation}
  \label{Psi2L}
  \BL=\frac{(1-a_{0}^{2})\BGs+i(a_{0}\det\BGs-(\rho+a_{0})(a_{0}\rho+1))\BR_{\perp}}
{\det\BGs-(\rho+a_{0})^{2}}.
\end{equation}
Let us now choose $a_{0}$ so that the matrix $\BL$ is real and symmetric. This
means that $a_{0}$ must be a root of
\begin{equation}
  \label{diaga0}
  (a_{0}^{2}+1)\rho=a_{0}(\det\BGs-(\rho^{2}+1)).
\end{equation}
It is easy to check that the discriminant of the quadratic equation
(\ref{diaga0}) is nonnegative \IFF 
\begin{equation}
  \label{wkcpld}
  |r_{1}-r_{2}|\le\left|\sqrt{\det\BGs_{1}}-\sqrt{\det\BGs_{2}}\right|.
\end{equation}
We will call this case ``weakly coupled'' because there is a choice of $a_{0}$
(a root of (\ref{diaga0})) that eliminates the thermoelectric coupling from
both materials. The exceptional cases are $\det\BGs=(\rho\pm 1)^{2}$, in
which case the inverse in the definition of $\BL$ does not exist. We note that
there are several possibilities.
\begin{enumerate}
\item We could have mapped $\SFL_{2}$ to $\SFI$, instead of $\SFL_{1}$. That
  means that instead of a pair of numbers $(\rho,\det\BGs)$ we will get a pair
  of numbers
\[
(\rho',\det\BGs')=\left(-\frac{\rho}{\sqrt{\det\BGs}},\nth{\det\BGs}\right).
\]
It is easy to check that the sign of the discriminant of (\ref{diaga0}) does
not depend on the choice of which $\SFL_{j}$ gets mapped to $\SFI$.
\item In each case there are two real roots $a_{0}$. If one root is $a_{0}$,
  the other root is $1/a_{0}$. 
\end{enumerate}
We want $\BL>0$.
This means that in each case we must have the inequality
\[
\frac{1-a_{0}^{2}}{\det\BGs-(\rho+a_{0})^{2}}>0.
\]
In the weakly coupled case (\ref{diaga0}) we have
\[
\frac{1-a_{0}^{2}}{\det\BGs-(\rho+a_{0})^{2}}=\frac{a_{0}}{\rho+a_{0}}=\frac{a_{0}\rho+1}{\det\BGs}.
\]
If $a_{1}$ and $a_{2}$ are the two roots of (\ref{diaga0}), then
\[
\frac{a_{1}}{\rho+a_{1}}\cdot\frac{a_{2}}{\rho+a_{2}}=\nth{\det\BGs}>0,
\]
\[
\frac{a_{1}}{\rho+a_{1}}+\frac{a_{2}}{\rho+a_{2}}=\frac{\det\BGs-\rho^{2}+1}{\det\BGs}.
\]
The inequality $\det\BGs-\rho^{2}+1>0$
always holds in the weakly coupled case.

In the strongly coupled case, i.e. when inequality (\ref{wkcpld}) is violated,
we first apply the transformation
\[
\Psi_{0}(\SFL)=(\BGs_{1}^{-1/2}\otimes\BI_{2})\SFL(\BGs_{1}^{-1/2}\otimes\BI_{2})
\] 
which mapps $\BL_{1}$ into
$\BL'_{1}=\BI_{2}+i\rho_{1}\BR_{\perp}$, and $\BL_{2}$ into
$\BL'_{2}=\BGs+i\rho_{2}\BR_{\perp}$, where
\[
\rho_{j}=\frac{r_{j}}{\sqrt{\det{\BGs_{1}}}}.
\]
In this case we look for a simpler transformation
\[
\Psi_{3}(\SFL)=a\SFL+b\SFT,
\]
which maps $\SFL'_{j}$ onto the ER $(\ClD,\ClD')$ if we work in the frame in
which $\BGs_{11}=\BGs_{22}$. The coefficients $a>0$ and $b\in\bb{R}$ need to
be chosen so that
\[
\det(a\BI_{2}+i(a\rho_{1}+b)\BR_{\perp})=\det(a\BGs+i(a\rho_{2}+b)\BR_{\perp})=1
\]
This means that
\[
a^{2}=(a\rho_{1}+b)^{2}+1,\qquad a^{2}\det\BGs=(a\rho_{2}+b)^{2}+1.
\]
Subtracting the two equations we find
\[
b=a\frac{\det\BGs-1+\rho_{1}^{2}-\rho_{2}^{2}}{2\rho},
\]
where $\rho$ is the same as before. We then find that
\[
a^{2}\frac{((\rho+1)^{2}-\det\BGs)(\det\BGs-(\rho-1)^{2})}{4\rho^{2}}=1.
\]
We note that since $|\rho_{1}|<1$ and $|\rho_{2}|<\sqrt{\det\BGs}$, then
\[
|\rho|=|\rho_{2}-\rho_{1}|\le|\rho_{2}|+|\rho_{1}|<\sqrt{\det\BGs}+1.
\]
This is equivalent to
\[
|r_{2}-r_{1}|<\sqrt{\det\BGs_{1}}+\sqrt{\det\BGs_{2}}.
\]
Using our original parameters we can write
\[
a^{2}=\frac{4\GD r^{2}\det\BGs_{1}}{
(\GD r^{2}-(\sqrt{\det\BGs_{1}}-\sqrt{\det\BGs_{2}})^{2})
((\sqrt{\det\BGs_{1}}+\sqrt{\det\BGs_{2}})^{2}-\GD r^{2})},\qquad\GD r=r_{2}-r_{1}
\]
This shows that in the strongly coupled regime we can always choose $a>0$.
We will therefore write 
\[
a=2a_{0}\sqrt{\det\BGs_{1}},\quad
a_{0}=\frac{|\GD r|}{\sqrt{(\GD r^{2}-(\sqrt{\det\BGs_{1}}-\sqrt{\det\BGs_{2}})^{2})
((\sqrt{\det\BGs_{1}}+\sqrt{\det\BGs_{2}})^{2}-\GD r^{2})}}
\]
\[
b=a_{0}\frac{\det\BGs_{2}+r_{1}^{2}-r_{2}^{2}}{r_{2}-r_{1}}.
\]

Let us find the conditions for each case, including special cases explicitly
in terms of $\BGs_{j}$, $r_{j}$, $j=1,2$. We compute
\[
\BX_{0}=2\left(\BGs_{1}^{-1/2}\BGs_{2}\BGs_{1}^{-1/2}+\BI_{2}+\frac{r_{2}-r_{1}}{\sqrt{\det{\BGs_{1}}}}i\BR_{\perp}\right)^{-1}-\BI_{2}.
\]
It fairly easy to show that $\BY\not=0$ \IFF $\BGs_{1}$ and $\BGs_{1}$ are not scalar
multiples of one another. it is also easy to show that $\BF=0$ \IFF $r_{1}=r_{2}$ and
$\det\BGs_{1}=\det\BGs_{2}$. In order to split $\BX_{0}$ into the $\BGF$ and
$\BGY$ parts we can use the formula
\[
(\Gvf(\Ga)+\psi(a)+ir\BR_{\perp})^{-1}=\frac{\Gvf(\Ga)-\psi(a)-ir\BR_{\perp}}{\Ga^{2}-|a|^{2}-r^{2}}.
\]
So, if $\BF=\Gvf(\Ga)+ir\BR_{\perp}$ and $\BY=\psi(a)$, then we have the
equation
\[
\BGs_{1}^{-1/2}\BGs_{2}\BGs_{1}^{-1/2}+\BI_{2}+\frac{r_{2}-r_{1}}{\sqrt{\det{\BGs_{1}}}}i\BR_{\perp}=2\frac{\Gvf(\Ga+1)-\psi(a)-ir\BR_{\perp}}{(\Ga+1)^{2}-|a|^{2}-r^{2}}.
\]
Hence we have the following system
\[
\begin{cases}
  \frac{r_{2}-r_{1}}{\sqrt{\det{\BGs_{1}}}}=-\frac{2r}{(\Ga+1)^{2}-|a|^{2}-r^{2}},\\
\Trc(\BGs_{2}\BGs_{1}^{-1})+2=\frac{4(\Ga+1)}{(\Ga+1)^{2}-|a|^{2}-r^{2}},\\
\frac{\det\BGs_{2}}{\det\BGs_{1}}=\det\left(2\frac{\Gvf(\Ga+1)-\psi(a)}{(\Ga+1)^{2}-|a|^{2}-r^{2}}-\Gvf(1)\right).
\end{cases}
\]
The condition that $\BF\in\bb{R}\BZ_{0}$ or $\BF\in\bb{R}\bra{\BZ_{0}}$ is
equivalent to $r^{2}=\Ga^{2}$.
\[
\begin{cases}
  |a|^{2}=\frac{\Trc(\BGs_{1}\cof(\BGs_{2}))^{2}-4\det(\BGs_{1}\BGs_{2})}{(\det(\BGs_{1}+\BGs_{2})-(r_{1}-r_{2})^{2})^{2}},\\
r^{2}=\frac{4(r_{1}-r_{2})^{2}\det(\BGs_{1})}{(\det(\BGs_{1}+\BGs_{2})-(r_{1}-r_{2})^{2})^{2}},\\
\Ga=\frac{(r_{1}-r_{2})^{2}+\det(\BGs_{1})-\det(\BGs_{2})}{\det(\BGs_{1}+\BGs_{2})-(r_{1}-r_{2})^{2}}
\end{cases}
\]
The condition that $\BF\in\bb{R}\BZ_{0}$ or $\BF\in\bb{R}\bra{\BZ_{0}}$ is
equivalent to
\begin{equation}
  \label{RZ0}
  |r_{1}-r_{2}|=\left|\sqrt{\det\BGs_{1}}-\sqrt{\det\BGs_{2}}\right|.
\end{equation}
The composite made with the pair of isotropic materials $\SFL_{j}$ will be
called weakly thermoelectrically heterogenerous if
\begin{equation}
  \label{weak}
  |r_{1}-r_{2}|<\left|\sqrt{\det\BGs_{1}}-\sqrt{\det\BGs_{2}}\right|,
\end{equation}
and strongly coupled otherwise. (We note that requirement of positive
definteness of $\SFL_{j}$ implies that $r_{j}^{2}<\det\BGs_{j}$.) The
thermoelectric interactions in a weakly thermoelectrically heterogenerous
composite can be decoupled. Such a decoupling is impossible in a strongly
thermoelectrically heterogenerous composite. Here is the summary.
\begin{enumerate}
\item $\BGs_{1}\not=\Gth\BGs_{2}$ (generic case)
  \begin{enumerate}
  \item $|r_{1}-r_{2}|<\left|\sqrt{\det\BGs_{1}}-\sqrt{\det\BGs_{2}}\right|$
    (weakly coupled case)
    \begin{enumerate}
    \item $|r_{1}-r_{2}|^{2}\not=\det(\BGs_{1}-\BGs_{2})$: ER $(\ClD,\ClD)$.
      This implies that the 10 components of $\SFL^{*}$ depend only on 6
      microstructure-dependent parameters. Moreover, the link
      $(\ClD,\ClD)/{\rm Ann}(\bb{C}\Be_{2})\cong{\rm Ann}(\bb{C}\Be_{2})$ that relates
      $\SFL^{*}$ to the effective tensor of a 2D conducting composite applies.
    \item $|r_{1}-r_{2}|^{2}=\det(\BGs_{1}-\BGs_{2})$: ER Ann$(\bb{C}\Be_{2})$.
      This implies that the 10 components of $\SFL^{*}$ depend only on 3
      microstructure-dependent parameters. Moreover, they can be expressed in
      terms of the effective tensor of a 2D conducting composite.
    \end{enumerate}
  \item $|r_{1}-r_{2}|>\left|\sqrt{\det\BGs_{1}}-\sqrt{\det\BGs_{2}}\right|$:
    ER $(\ClD,\ClD')$ (strongly coupled case)
    This implies that the 10 components of $\SFL^{*}$ depend
    only on 6 microstructure-dependent parameters.  
  \item $|r_{1}-r_{2}|=\left|\sqrt{\det\BGs_{1}}-\sqrt{\det\BGs_{2}}\right|$
  (borderline strongly coupled case)
    \begin{enumerate}
    \item $r_{1}\not=r_{2}$: ER $(W,V_{\infty})$. This implies that the 10
      components of $\SFL^{*}$ depend only on 6 microstructure-dependent
      parameters. Moreover, the link $(W,V_{\infty})/{\rm
        Ann}(\bb{C}\bra{\Bz_{0}})\cong(\bb{C}\BI,\bb{R}\psi(i))$ is applicable. This
      means that there is a link between this case and the $r_{1}=r_{2}$ case,
      since ERs $(\bb{C}\BI,\bb{R}\psi(i))$ and $(\bb{C}\BI,\bb{R}\psi(1))$
      are isomorphic.
    \item $r_{1}=r_{2}$: ER $(\bb{C}\BI,\bb{R}\psi(1))$. This implies that the
      10 components of $\SFL^{*}$ depend only on 3 microstructure-dependent
      parameters, expressible in terms of 2D conductivity.
    \end{enumerate}
\end{enumerate}
\item $\BGs_{1}=\Gth_{1}\BGs$, $\BGs_{2}=\Gth_{2}\BGs$ (special nongeneric case)
  \begin{enumerate}
  \item $|r_{1}-r_{2}|<|\Gth_{1}-\Gth_{2}|\sqrt{\det\BGs}$: 
    ER $(\bb{C}\BI,\bb{R}\BI)$. This implies that the 10 components of $\SFL^{*}$ depend
    only on 3 microstructure-dependent parameters, which are expressible in
    terms of the effective tensor of a 2D conducting composite.
  \item $|r_{1}-r_{2}|>|\Gth_{1}-\Gth_{2}|\sqrt{\det\BGs}$: 
    ER $(\bb{C}\BI,i\BR_{\perp})$. This implies that the 10 components of $\SFL^{*}$ depend
    only on 3 microstructure-dependent parameters.
  \item $|r_{1}-r_{2}|=|\Gth_{1}-\Gth_{2}|\sqrt{\det\BGs}$:
    ER $(0,\bb{R}\BZ_{0})$. This implies $\SFL^{*}$ is completely determined, regardless
    of microstructure, if the volue fractions are known.
  \end{enumerate}
\end{enumerate}

An example is a binary composite made with isotropic thermoelectric materials
in which the Seebeck coefficient $\BS$ is a scalar. In this case we have
$r_{1}=r_{2}=0$. Let $\GS^{*}(h)$ be the effective tensor of an isotropic
conducting composite made with two isotropic materials, whose conductivities
are 1 and $h$. Then $\BGS^{*}$ can be applied to symmetric matrices according
to the rule 
\[
\GS^{*}(\BS)=\BR\mat{\GS^{*}(s_{1})}{0}{0}{\GS^{*}(s_{2})}\BR^{T},\qquad
\BS=\BR\mat{s_{1}}{0}{0}{s_{2}}\BR^{T}.
\]
Then the effective thermoelectric tensor of such a composite will be isotropic
and also have scalar Seebeck coefficient
$r^{*}=0$ and
\[
\SFL^{*}=\BGs^{*}\otimes\BI_{2},\qquad
\BGs^{*}=\BGs_{1}^{1/2}\GS^{*}(\BGs_{1}^{-1/2}\BGs_{2}\BGs_{1}^{-1/2})\BGs_{1}^{1/2}.
\]
\subsection{Results}
In many cases the results will be formulated in terms of the \mc-dependent function
$\BGS(h)$, representing the effictive conductivity of a two-phase composite
with isotropic constituent conductivities $\BI_{2}$ and $h\BI_{2}$, replacing
materials $\SFL_{1}$, $\SFL_{2}$ in the original composite. Another convenient
notation will be matrices $\BS_{1}$ and $\BS_{2}$ definied by
\[
\BS_{1}=\frac{\BGs_{2}-\Gl_{1}\BGs_{1}}{\Gl_{2}-\Gl_{1}},\qquad
\BS_{2}=\frac{\BGs_{2}-\Gl_{2}\BGs_{1}}{\Gl_{1}-\Gl_{2}},
\]
where $\Gl_{1}$ and $\Gl_{2}$ are the two roots of the quadratic equation
$\det(\BGs_{2}-\Gl\BGs_{1})=0$. The ordering of the roots is unimportant. This
notation will not be used when $\BGs_{1}$ and $\BGs_{2}$ are scalar multiples
of one another, since in this case $\Gl_{1}=\Gl_{2}$ and the
matrices $S_{j}$ are undefined.

$\bullet$ (2(c)) $\BGs_{1}=\Gth_{1}\BGs_{0}$, $\BGs_{2}=\Gth_{2}\BGs_{0}$,
$|r_{1}-r_{2}|=|\Gth_{1}-\Gth_{2}|\sqrt{\det\BGs_{0}}$. We apply the global
automorphism
\[
\Psi(\SFL)=\nth{\Gth_{1}}(\BGs_{0}^{-1/2}\otimes\BI_{2})(\SFL-r_{1}\SFT)(\BGs_{0}^{-1/2}\otimes\BI_{2}),
\]
which maps $\SFL_{1}$ into $\SFI$ and $\SFL_{2}$ into
\[
\SFL_{0}=\mat{\frac{\Gth_{2}}{\Gth_{1}}\BI_{2}}{-\frac{r_{2}-r_{1}}{\Gth_{1}\sqrt{\det\BGs_{0}}}\BR_{\perp}}{\frac{r_{2}-r_{1}}{\Gth_{1}\sqrt{\det\BGs_{0}}}\BR_{\perp}}{\frac{\Gth_{2}}{\Gth_{1}}\BI_{2}},
\]
corresponding to the ER $(0,\bb{R}\BZ_{0})$ or $(0,\bb{R}\bra{\BZ_{0}})$.
We choose indexing in such a way that $\Gth_{2}\ge \Gth_{1}$, so that
$\SFL_{0}>0$. If we denote
\[
\Gl=\frac{\Gth_{2}}{\Gth_{1}}\ge 1,
\]
then either
\[
\frac{r_{2}-r_{1}}{\Gth_{1}\sqrt{\det\BGs_{0}}}=\Gl-1,\text{ or }
\frac{r_{2}-r_{1}}{\Gth_{1}\sqrt{\det\BGs_{0}}}=-(\Gl-1),
\]
depending on whether
$r_{1}-r_{2}=(\Gth_{1}-\Gth_{2})\sqrt{\det\BGs_{0}}$ or
$r_{1}-r_{2}=-(\Gth_{1}-\Gth_{2})\sqrt{\det\BGs_{0}}$, respectively.
In either case the effective tensor of a composite made with $\SFI$ and
$\SFL_{0}$ will have the form
\[
\SFL^{*}_{0}=\mat{\Gl^{*}\BI_{2}}{\mp(\Gl^{*}-1)\BR_{\perp}}{\pm(\Gl^{*}-1)\BR_{\perp}}{\Gl^{*}\BI_{2}},
\qquad\nth{\Gl^{*}}=\av{\Gl(x)^{-1}}=f_{1}+\frac{f_{2}}{\Gl}.
\]
Applying $\Psi^{-1}$ transformation we obtain the formula for
$\SFL^{*}=K(\BGs^{*}+ir^{*}\BR_{\perp},0)$: 
\[
  \BGs^{*}=(f_{1}\BGs_{1}^{-1}+f_{2}\BGs_{2}^{-1})^{-1},\qquad
r^{*}=r_{1}+\frac{(r_{2}-r_{1})f_{2}\Gth_{1}}{f_{1}\Gth_{2}+f_{2}\Gth_{1}}=
\frac{r_{1}f_{1}\Gth_{1}^{-1}+r_{2}f_{2}\Gth_{2}^{-1}}{f_{1}\Gth_{1}^{-1}+f_{2}\Gth_{2}^{-1}}.
\]
In other words,
\begin{equation}
  \label{2c}
\BGs^{*}=\av{\BGs(\Bx)^{-1}}^{-1},\qquad r^{*}=\frac{\av{r(x)\Gth(x)^{-1}}}{\av{\Gth(x)^{-1}}}.
\end{equation}

  \textbf{In every case we will work in the frame in which
    $\BGs_{1}^{-1/2}\BGs_{2}\BGs_{1}^{-1/2}$ is diagonal! This is not a
    physical frame. Rather it is a mathematical one, where two different linear
    combinations of the original curl-free and divergence-free fields are
    chosen as a pair of intensity and flux fields.}

$\bullet$ (1(c)ii) $\BL_{1}=\BGs_{1}+ir_{0}\BR_{\perp}$,
  $\BL_{2}=\BGs_{2}+ir_{0}\BR_{\perp}$, moreover, $\det\BGs_{1}=\det\BGs_{2}$.
The global link (\ref{Psi0}) maps $\BL_{1}$ to $\BI_{2}$ and $\BL_{2}$ to
$\BGs=\BGs_{1}^{-1/2}\BGs_{2}\BGs_{1}^{-1/2}$, which is always assumed to be
diagonal. Since in this case we have
$\det\BGs=1$ we can write
\[
\BGs=\mat{\Gl}{0}{0}{\Gl^{-1}},
\]
where the eigenvalues $\Gl$ and $1/\Gl$ solve
the quadratic equation 
\begin{equation}
  \label{geneig12}
  \det(\BGs_{2}-\Gl\BGs_{1})=0.
\end{equation}
The effective tensor of the resulting composite is
\[
\SFL_{0}^{*}=\mat{\BGs^{*}}{0}{0}{\frac{\BGs^{*}}{\det\BGs^{*}}}=
\mat{1}{0}{0}{\nth{\det\BGs^{*}}}\otimes
\BGs^{*},\qquad\BGs^{*}=\BGS(\Gl).
\]
where $\BGs^{*}$ is the effective conductivity of the composite made
with 2 isotropic materials of conductivities 1 and $\Gl$ and the same \mc\ as
the original composite. % Let
% \[
% \Gd(\Gl)=\sqrt{\det\BGs^{*}(\Gl)}.
% \]
% Then, $\Gd(\Gl^{-1})=\Gd(\Gl)^{-1}$, since 
% \[
% \BGs^{*}\left(\nth{\Gl}\right)=\frac{\BGs^{*}(\Gl)}{\det\BGs^{*}(\Gl)}.
% \]
We conclude that
\begin{equation}
  \label{1cii}
  \SFL^{*}=\BGs_{1}^{1/2}\mat{1}{0}{0}{\nth{\det\BGs^{*}}}\BGs_{1}^{1/2}
  \otimes\BGs^{*}+r_{0}\SFT,
\end{equation}
where we have used the the form in which terms on both sides of the $\otimes$
sign are invariant with respect to the $\Gl\mapsto\Gl^{-1}$ permutation.
We can express the answer without using $\BGs_{1}^{1/2}$ by observing that
\[
\BGs_{1}^{1/2}\BI_{2}\BGs_{1}^{1/2}=\BGs_{1},\qquad
\BGs_{1}^{1/2}\BGs\BGs_{1}^{1/2}=\BGs_{2}.
\]
In the frame in which $\BGs$ is diagonal we obtain
\[
\SFL^{*}=r_{0}\SFT+((\det\BGs^{*}-\Gl^{2})\BGs_{1}+\Gl(1-\det\BGs^{*})\BGs_{2})
\otimes\frac{\BGs^{*}}{(1-\Gl^{2})\det\BGs^{*}},
\]
where $\Gl>0$ solves (\ref{geneig12}) and $\BGs^{*}=\GS_{\rm 2D cond}(1,\Gl)$,
and the result is independent of the choice of the root in (\ref{geneig12}).
We note that in the notation $\BA\otimes\BB$, the frame of the operator $\BB$
is the physical one, while the frame of $\BA$ is mathematical. In the formula
for $\SFL^{*}$ above the first factor of the tensor product would look the
same in the original frame, since the tensors $\BGs_{j}$ are transforming
together with the mathematical frame. Using our notation $\BS_{1}$ and
$\BS_{2}$ are can write the final answer as
\[
\SFL^{*}=r_{0}\SFT+\left(\frac{\BS_{1}}{\det\BGs^{*}}+\BS_{2}\right)\otimes\BGs^{*}.
\]

$\bullet$ (1(a)ii) $|r_{1}-r_{2}|^{2}=\det(\BGs_{1}-\BGs_{2})$. In this case the
quadratic equation (\ref{diaga0}) has two solutions $(\Gl_{j}-1)/\rho$,
$j=1,2$, where the eigenvalues $\Gl_{1}$ and $\Gl_{2}$ of $\BGs$ solve
(\ref{geneig12}). If we choose $a_{0}=(\Gl_{1}-1)/\rho$, then
\begin{equation}
  \label{L0st1aii}
  \BS=\mat{\frac{\Gl_{1}}{\Gl_{2}}}{0}{0}{1},\qquad
\SFL_{0}^{*}=\mat{\BGs^{*}}{0}{0}{\BI_{2}},\qquad\BGs^{*}=\BGS\left(\frac{\Gl_{1}}{\Gl_{2}}\right).
\end{equation}

Inverting our global link that mappled $\SFL_{1}$ to $\tns{\BI_{2}}$ and
$\SFL_{2}$ to $\BS\otimes\BI_{2}$ we obtain
\[
\SFL^{*}=r_{1}\SFT+(\BGs_{1}^{1/2}\otimes\BI_{2})(\SFT-a_{0}\SFL_{0}^{*})(\SFL_{0}^{*}-a_{0}\SFT)^{-1}\SFT(\BGs_{1}^{1/2}\otimes\BI_{2}).
\]
Observe that
\[
\SFP^{*}=(\SFT-a_{0}\SFL_{0}^{*})(\SFL_{0}^{*}-a_{0}\SFT)^{-1}\SFT=
(1-a_{0}^{2})\SFT(\SFL_{0}^{*}-a_{0}\SFT)^{-1}\SFT-a_{0}\SFT.
\]
We then compute (verified by Maple)
\[
\SFP^{*}_{0}=(1-a_{0}^{2})\SFT(\SFL_{0}^{*}-a_{0}\SFT)^{-1}\SFT=\frac{1-a_{0}^{2}}{\det(\BGs^{*}-a_{0}^{2})}
\mat{\det\BGs^{*}-a_{0}^{2}\BGs^{*}}{a_{0}\BR_{\perp}^{T}(\BGs^{*}-a_{0}^{2})}
{a_{0}(\BGs^{*}-a_{0}^{2})\BR_{\perp}}{\BGs^{*}-a_{0}^{2}}.
\]
Thus,
\begin{equation}
  \label{Lst1aii}
  \SFL^{*}=(r_{1}-a_{0}\sqrt{\det\BGs_{1}})\SFT+(\BGs_{1}^{1/2}\otimes\BI_{2})\SFP^{*}_{0}
(\BGs_{1}^{1/2}\otimes\BI_{2}).
\end{equation}

In order to write $\SFL^{*}$ without the explicit reference to $\BGs_{1}^{1/2}$ we 
write $\SFP_{0}^{*}$ as a sum of tensor products. In order to accomplish this
we first write $\BGs^{*}=\phi(x^{*})+\psi(y^{*})$ and then observe that 
\[
\BR_{\perp}^{T}\psi(y^{*})=-\psi(iy^{*})=\psi(y^{*})\BR_{\perp}.
\]
Thus, we get
\[
\SFP^{*}_{0}=\frac{1-a_{0}^{2}}{\det(\BGs^{*}-a_{0}^{2})}\left(
\mat{\det\BGs^{*}}{0}{0}{-a_{0}^{2}}\otimes\BI_{2}+\mat{-a_{0}^{2}}{0}{0}{1}\otimes\BGs^{*}
+a_{0}(x^{*}-a_{0}^{2})\SFT-a_{0}\psi(i)\otimes\psi(iy^{*})
\right)
\]
We will then be able to compute $\SFL^{*}$ in the covariant form, if we can
express $\BGs_{1}^{1/2}\psi(i)\BGs_{1}^{1/2}$, in the covariant form.
To do this we first write $\psi(i)=\phi(i)\psi(1)$, and then
\[
\BGs_{1}^{1/2}\psi(i)\BGs_{1}^{1/2}=
\BGs_{1}^{1/2}\phi(i)\BGs_{1}^{1/2}\BGs_{1}^{-1}\BGs_{1}^{1/2}\psi(1)\BGs_{1}^{1/2}=
\sqrt{\det\BGs_{1}}\BR_{\perp}\BGs_{1}^{-1}(\BGs_{1}^{1/2}\psi(1)\BGs_{1}^{1/2}).
\]
Thus, we compute
\[
\BGs_{1}^{1/2}\psi(1)\BGs_{1}^{1/2}=\frac{(\Gl_{1}+\Gl_{2})\BGs_{1}-2\BGs_{2}}{\Gl_{2}-\Gl_{1}}
=\BS_{2}-\BS_{1}.
\]
Therefore,
\[
\BGs_{1}^{1/2}\psi(i)\BGs_{1}^{1/2}=\sqrt{\det\BGs_{1}}\BR_{\perp}\BGs_{1}^{-1}(\BS_{2}-\BS_{1}),
\]
which is now in a frame-covariant form. Putting everything together we obtain
\begin{multline}
  \label{1aiimonster}
\SFL^{*}=\left(r_{1}+a_{0}\left(\frac{(1-a_{0}^{2})(x^{*}-a_{0}^{2})}
{\det(\BGs^{*}-a_{0}^{2})}-1\right)\sqrt{\det\BGs_{1}}\right)\SFT+
\frac{1-a_{0}^{2}}{\det(\BGs^{*}-a_{0}^{2})}\times\\
\{\BS_{1}\otimes(\BGs^{*}-a_{0}^{2})+\BS_{2}\otimes(\det\BGs^{*}-a_{0}^{2}\BGs^{*})+
a_{0}\sqrt{\det\BGs_{1}}\SFT[\BGs_{1}^{-1}(\BS_{1}-\BS_{2})\otimes(\BGs^{*}-x^{*})]\}
\end{multline}
This expression is invariant with respect to the choice of $\Gl_{1}$ and
$\Gl_{2}$. This means that if we interchange $\Gl_{1}$ and $\Gl_{2}$,
$\BS_{1}$ and $\BS_{2}$,  repalce $a_{0}$ with
$1/a_{0}$, and $\BGs^{*}$ with $\BGs^{*}/\det\BGs^{*}$, the value of $\SFL^{*}$
will not change.  

We also want to see if there is a symmetry wrt to material index
interchange ``1''$\mapsto$''2'' (the result should not depend on naming conventions). Hence we
need to make the following replacements:
\[
\det\BGs_{1}\mapsto\Gl_{1}\Gl_{2}\det\BGs_{1},\ 
a_{0}\mapsto a_{0}\sqrt{\frac{\Gl_{2}}{\Gl_{1}}},\
\BGs^{*}\mapsto\frac{\Gl_{2}}{\Gl_{1}}\BGs^{*},\ \Gl_{j}\mapsto\nth{\Gl_{j}},\
\BS_{j}\mapsto\Gl_{1}\Gl_{2}\frac{\BS_{j}}{\Gl_{j}}.
\]
We can verify that the second term is invariant. This is due to the
observation that, since
\[
a_{0}=\frac{\Gl_{1}-1}{\rho}=\frac{\rho}{\Gl_{2}-1},
\]
we can also write
\[
a_{0}^{2}=\frac{\Gl_{1}-1}{\Gl_{2}-1}\quad\thus\quad
1-a_{0}^{2}=\frac{\Gl_{2}-\Gl_{1}}{\Gl_{2}-1}.
\]
This shows that under the index-interchange we also have
$1-a_{0}^{2}\mapsto\nth{\Gl_{1}}(1-a_{0}^{2})$. The first term was also
verified to be invariant, since we can write
\[
r_{1}=\frac{r_{1}+r_{2}}{2}-\hf\rho\sqrt{\det\BGs_{1}}=\frac{r_{1}+r_{2}}{2}
-\hf a_{0}(\Gl_{2}-1)\sqrt{\det\BGs_{1}}.
\]

If $\BGs^{*}=x^{*}\BI_{2}$, the result simplifies:
\begin{equation}
  \label{1aiiiso}
  \SFL^{*}=\left(r_{1}+\frac{(1-x^{*})a_{0}}
{x^{*}-a_{0}^{2}}\sqrt{\det\BGs_{1}}\right)\SFT+
\frac{1-a_{0}^{2}}{x^{*}-a_{0}^{2}}(\BS_{1}+x^{*}\BS_{2})\otimes\BI_{2}.
\end{equation}
In order to verify (\ref{1aiimonster}) with Maple we can parametrize this case
by $\Bs=\BGs_{1}^{1/2}$, $\Gl_{1}$, $\rho$, $r_{1}$,
so that
\[
\BGs_{1}=\Bs^{2},\quad r_{2}=\rho\det\Bs+r_{1},\quad\Gl_{2}=\frac{\rho^{2}}{\Gl_{1}-1}+1,\quad
\BGs_{2}=\Bs\mat{\Gl_{1}}{0}{0}{\Gl_{2}}\Bs,\quad a_{0}=\frac{\Gl_{1}-1}{\rho}.
\]
Substituting these and $\SFL_{0}^{*}$ given by (\ref{L0st1aii}) into
(\ref{Lst1aii}) we obtain (\ref{1aiimonster}), as verified by Maple. Formula
(\ref{1aiiiso}) has also been verified.

Finally we point out that the case $a_{0}=\infty$ corresponding to
$r_{1}=r_{2}=r_{0}$ and $\det(\BGs_{1}-\BGs_{2})=0$ is also included by taking
a limit as  $a_{0}\to\infty$ in (\ref{1aiimonster}):
\[
\SFL^{*}=r_{0}\SFT+\BS_{1}\otimes\BI_{2}+\BS_{2}\otimes\BGS^{*}(\Gl_{1}),\quad\Gl_{1}\not=1=\Gl_{2}.
\]

$\bullet$ (1(c)i):
$|r_{1}-r_{2}|=\left|\sqrt{\det\BGs_{1}}-\sqrt{\det\BGs_{2}}\right|$. The
relevant ER is $(W,V_{0})$, described as 
\begin{equation}
  \label{WV0}
  \left\{\mat{\BL}{\pm(\BL\BM-\BR_{\perp})}{\pm(\BM^{T}\BL+\BR_{\perp})}{\BM^{T}\BL\BM}:\Trc\BM=0,\
\BL>0,\ \BL+2\BR_{\perp}\BM\det\BL<0\right\}.
\end{equation}
The relevant link is as follows. $\BM^{*}=\BR_{\perp}\BGs^{*}$, where $\BGs^{*}$
is the 2D effective condictivity tensor of the composite with local
conductivity $\Bc(\Bx)=-\BR_{\perp}\BM(\Bx)$, which is symmetric and positive
definite for $\BM(\Bx)$ satisfying the constraints in (\ref{WV0}).

Any isotropic tensor $\SFL$ in this ER must have the form
\[
\SFL=\mat{\Gl_{1}\BI_{2}}{\pm(\sqrt{\Gl_{1}\Gl_{2}}-1)\BR_{\perp}}
{\mp(\sqrt{\Gl_{1}\Gl_{2}}-1)\BR_{\perp}}{\Gl_{2}\BI_{2}},
\qquad\Gl_{1}>0,\ \Gl_{1}\Gl_{2}>\nth{4},
\]
corresponding to $\BM=\sqrt{\Gl_{2}/\Gl_{1}}\BR_{\perp}$.
 Thus,
\begin{equation}
  \label{L01ci}
  \SFL_{0}^{*}=\mat{\BL^{*}}{\pm\BR_{\perp}\cof(\BL^{*})\BGs^{*}}{\pm\BGs^{*}\cof(\BL^{*})\BR_{\perp}^{T}}{\BGs^{*}\cof(\BL^{*})\BGs^{*}}\pm\SFT,
\end{equation}
where $\BGs^{*}=\BGS^{*}(\sqrt{\Gl_{2}/\Gl_{1}})$, where $\Gl_{1,2}$ are the
two roots of (\ref{geneig12}), and $\BL^{*}$ is \mc-dependent tensor that
depends on material moduli only through $\Gl_{1}$ and $\Gl_{2}$. It is not
expressible in terms of the effective conductivity. To compute $\SFL^{*}$ we
take $\Psi^{-1}$ and obtain as in the case 1(a)ii
\[
\SFL^{*}=\BS_{1}\otimes\BGs^{*}\cof(\BL^{*})\BGs^{*}+\BS_{2}\otimes\BL^{*}
+\Ga\SFT\left[\BGs_{1}^{-1}(\BS_{1}-\BS_{2})\otimes\BA^{*}\right]+\Gb\SFT,
\]
where
\[
\BA^{*}=\cof(\BL^{*})\BGs^{*}-a^{*}\BI_{2},\quad
\BGs^{*}=\BGS^{*}\left(\sqrt{\frac{\Gl_{2}}{\Gl_{1}}}\right),\quad 
a^{*}=\hf\Trc(\cof(\BL^{*})\BGs^{*}),
\]
\[
\Ga=\frac{\sqrt{\det\BGs_{1}}-\sqrt{\det\BGs_{2}}}{r_{1}-r_{2}}\sqrt{\det\BGs_{1}},\qquad
\Gb=r_{1}+\Ga(a^{*}-1).
\]
\noindent$\bullet$ (2a) $\BGs_{1}=\Gth_{1}\BGs$, $\BGs_{2}=\Gth_{2}\BGs$,
$|r_{1}-r_{2}|<|\Gth_{1}-\Gth_{2}|\sqrt{\det\BGs}$. The relevant ER is
\[
\left\{\mat{\BL}{0}{0}{\BL}:\BL>0\right\}.
\]
and the relevant link is that $\BL^{*}$ is the effective conductivity tensor
of the conducting composite with local conductivity $\BL(\Bx)$.

We first apply (\ref{Psi0}) mappling $\SFL_{1}$ to $\SFI$ and $\SFL_{2}$ to 
\[
\SFL'_{2}=\frac{\Gth_{2}}{\Gth_{1}}\tns{\BI_{2}}+\rho\SFT.
\]
Let $a_{0}$ be a root of (\ref{diaga0}), where
$\BGs=(\Gth_{2}/\Gth_{1})\BI_{2}$. Then apply
\begin{equation}
  \label{Psi2}
\Psi_{2}(\SFL)=(a_{0}\SFL+\SFT)(\SFL+a_{0}\SFT)^{-1}\SFT=
\SFT(\SFL+a_{0}\SFT)^{-1}(a_{0}\SFL+\SFT).
\end{equation}
This gives $\Psi_{2}(\SFI)=\SFI$ and
\[
\Psi_{2}(\SFL'_{2})=\frac{(1-a_{0}^{2})\Gth_{1}\Gth_{2}}{\Gth_{2}^{2}-(\Gth_{1}\rho+\Gth_{1}a_{0})^{2}}\SFI.
\]
Then defining
\[
\BGs^{*}=\BGS^{*}\left(\frac{(1-a_{0}^{2})\Gth_{1}\Gth_{2}}{\Gth_{2}^{2}-(\Gth_{1}\rho+\Gth_{1}a_{0})^{2}}\right)=\BGS^{*}\left(\frac{a_{0}}{\rho+a_{0}}\right)=
\BGS^{*}\left(\frac{\Gth_{1}}{\Gth_{2}}(\rho a_{0}+1)\right)
\]
We obtain
\[
\SFL^{*}=\Psi_{1}^{-1}(\Psi_{2}^{-1}(\BI_{2}\otimes\BGs^{*})).
\]
Recalling that $\Psi_{2}^{-1}$ is $\Psi_{2}$ with $a_{0}$ replaced by $-a_{0}$
and that $\mathfrak{A}:\BA\otimes\BB\mapsto\BB\otimes\BA$ is the algebra
automorphism we can write the expression for
$\SFL'_{*}=\Psi_{2}^{-1}(\BGs^{*}\otimes\BI_{2})$ immediately from (\ref{Psi2L}):
\[
\SFL'_{*}=\frac{(1-a_{0}^{2})(\BI_{2}\otimes\BGs^{*})}{\det\BGs^{*}-a_{0}^{2}}
+\frac{a_{0}(1-\det\BGs^{*})}{\det\BGs^{*}-a_{0}^{2}}\SFT.
\]
Thus,
\[
\SFL^{*}=\frac{(1-a_{0}^{2})(\BGs_{1}\otimes\BGs^{*})}{\det\BGs^{*}-a_{0}^{2}}
+\left(r_{1}+\frac{a_{0}(1-\det\BGs^{*})\sqrt{\det\BGs_{1}}}{\det\BGs^{*}-a_{0}^{2}}\right)\SFT.
\]

\noindent$\bullet$ (1b) $|r_{1}-r_{2}|>\left|\sqrt{\det\BGs_{1}}-\sqrt{\det\BGs_{2}}\right|$.
The relevant ER is
\[
\bb{M}_{17}=\{\SFL>0:\SFL(\BJ\otimes\BR_{\perp})\SFL=\BJ\otimes\BR_{\perp}\},\quad\BJ=\mat{1}{0}{0}{-1}
\]
In this case 
\[
\SFL^{*}_{0}=a((\BGs_{1}^{-1/2}\otimes\BI_{2})\SFL^{*}(\BGs_{1}^{-1/2}\otimes\BI_{2}))+b\SFT.
\]
Then
\[
(a\SFL^{*}+b\sqrt{\det\BGs_{1}}\SFT)(\BGs_{1}^{-1/2}\BJ\BGs_{1}^{-1/2}\otimes\BR_{\perp})
(a\SFL^{*}+b\sqrt{\det\BGs_{1}}\SFT)=\BGs_{1}^{1/2}\BJ\BGs_{1}^{1/2}\otimes\BR_{\perp}
\]
Observe that
\[
\BGs_{1}^{-1/2}\BJ\BGs_{1}^{-1/2}=-\frac{\cof(\BGs_{1}^{1/2}\BJ\BGs_{1}^{1/2})}{\det\BGs_{1}}
=-\frac{\cof(\BM)}{\det\BGs_{1}},\qquad\BM=\BGs_{1}^{1/2}\BJ\BGs_{1}^{1/2}
\]
Substituting the values of $a$ and $b$
we obtain
\[
-(2\GD r\SFL^{*}+b_{0}\SFT)(\cof(\BM)\otimes\BR_{\perp})(2\GD r\SFL^{*}+b_{0}\SFT)
=c_{0}\BM\otimes\BR_{\perp},
\]
where
\[
b_{0}=\det\BGs_{2}-\det\BGs_{1}+r_{1}^{2}-r_{2}^{2},\quad
c_{0}=(\GD r^{2}-(\sqrt{\det\BGs_{1}}-\sqrt{\det\BGs_{2}})^{2})
((\sqrt{\det\BGs_{1}}+\sqrt{\det\BGs_{2}})^{2}-\GD r^{2})
\]

In order to compute matrix $\BM$ we recall that we now work in the frame where
\[
\BGs=\mat{s_{1}}{s_{2}}{s_{2}}{s_{1}}.
\]
The numbers $s_{1}$ and $s_{2}$ are related to the eigenvalues $\Gl_{1}$,
$\Gs_{2}$ of $\BGs$ via
\[
\Gl_{1}=s_{1}+s_{2},\qquad \Gl_{2}=s_{1}-s_{2}.
\]
Thus, from
\[
\BGs=\BGs_{1}^{-1/2}\BGs_{2}\BGs_{1}^{-1/2}=s_{1}\BI_{2}+s_{2}\psi(i)
\]
we obtain
\[
\BGs_{1}^{1/2}\psi(i)\BGs_{1}^{1/2}=\frac{\BGs_{2}-s_{1}\BGs_{1}}{s_{2}}=\BS_{2}-\BS_{1}=\GD\BS.
\]
Thus,
\[
\BM=\BGs_{1}^{1/2}\psi(i)\phi(i)\BGs_{1}^{1/2}=
(\BS_{2}-\BS_{1})\BGs_{1}^{-1}\sqrt{\det\BGs_{1}}\BR_{\perp}.
\]
Since $\BM$ is a symmetric matrix we obtain
\[
\BM=s(\BGs_{2}\BGs_{1}^{-1}\BR_{\perp})_{\rm sym}
\]
for some scalar $s$.
So, we have
\begin{equation}
  \label{2a}
(\SFL^{*}+\Gb_{0}\SFT)(\cof(\BA(\BGs_{1},\BGs_{2}))\otimes\BR_{\perp})(\SFL^{*}+\Gb_{0}\SFT)
=-\Gg_{0}\BA(\BGs_{1},\BGs_{2})\otimes\BR_{\perp},
\end{equation}
\[
\BA(\BGs_{1},\BGs_{2})=(\BGs_{2}\BGs_{1}^{-1}\BR_{\perp})_{\rm sym},\quad
\Gb_{0}=\frac{b_{0}}{2\GD r},\quad\Gg_{0}=\frac{c_{0}}{4(\GD r)^{2}}.
\]

We have derived the equation for $\SFL^{*}$:
\[
(a\SFS\SFL^{*}\SFS+b\SFT)\SFJ(a\SFS\SFL^{*}\SFS+b\SFT)=\SFJ,
\]
written to highlight the structure. Now factor out the $\SFS$s:
\[
\SFS(a\SFL^{*}+b\SFS^{-1}\SFT\SFS^{-1})\SFS\SFJ\SFS(a\SFL^{*}+b\SFS^{-1}\SFT\SFS^{-1})\SFS=\SFJ,
\]
Now multiply by $\SFS^{-1}$ on both sides:
\[
(a\SFL^{*}+b\SFS^{-1}\SFT\SFS^{-1})\SFS\SFJ\SFS(a\SFL^{*}+b\SFS^{-1}\SFT\SFS^{-1})=\SFS^{-1}\SFJ\SFS^{-1}.
\]
Observe now that all the tensors aside from $\SFL^{*}$ are nice tensor
products. Multiply them:
\[
\SFS^{-1}\SFT\SFS^{-1}=\sqrt{\det\BGs_{1}}\SFT,\quad
\SFS\SFJ\SFS=\BGs_{1}^{-1/2}\BJ\BGs_{1}^{-1/2}\otimes\BR_{\perp},\quad
\SFS^{-1}\SFJ\SFS^{-1}=\BGs_{1}^{1/2}\BJ\BGs_{1}^{1/2}\otimes\BR_{\perp}.
\]
We also have (using $\BJ^{-1}=\BJ$)
\[
\BGs_{1}^{-1/2}\BJ\BGs_{1}^{-1/2}=\frac{\cof(\BGs_{1}^{1/2}\BJ\BGs_{1}^{1/2})}{\det(\BGs_{1}^{1/2}\BJ\BGs_{1}^{1/2})}=-\frac{\cof(\BGs_{1}^{1/2}\BJ\BGs_{1}^{1/2})}{\det\BGs_{1}}.
\]
It remains to figure out $\BGs_{1}^{1/2}\BJ\BGs_{1}^{1/2}$. Before we do it,
divide the equation by $a^{2}$ and use constants $A$ and $B$ from case 2b.
\[
(\SFL^{*}+A\SFT)\SFS\SFJ\SFS(\SFL^{*}+A\SFT)=B\frac{\SFS^{-1}\SFJ\SFS^{-1}}{\det\BGs_{1}}.
\]
Denoting $\BZ=\BGs_{1}^{1/2}\BJ\BGs_{1}^{1/2}$ we obtain
\[
-(\SFL^{*}+A\SFT)\SFT(\BZ\otimes\BR_{\perp})\SFT(\SFL^{*}+A\SFT)=B\BZ\otimes\BR_{\perp}.
\]

We work in the frame where 
\[
\BGs=\mat{\mu}{\nu}{\nu}{\mu}.\]
It is easy to see that the eigenvalues of $\BGs$ are $\mu+\nu$ and
$\mu-\nu$. These eigenvalues have been denoted $\Gl_{1}$ and $\Gl_{2}$ before,
so
\[
\Gl_{1}=\mu+\nu,\qquad\Gl_{2}=\mu-\nu.
\]
Solving for $\mu$ and $\nu$ we obtain
\[
\mu=\frac{\Gl_{1}+\Gl_{2}}{2},\qquad\nu=\frac{\Gl_{1}-\Gl_{2}}{2}.
\]
Using relations $\BGs_{1}^{1/2}\BGs\BGs_{1}^{1/2}=\BGs_{2}$ and $\BGs_{1}^{1/2}\BI_{2}\BGs_{1}^{1/2}=\BGs_{1}$ we obtain
\[
\BGs_{2}=\BGs_{1}^{1/2}\mat{\mu}{\nu}{\nu}{\mu}\BGs_{1}^{1/2}=\mu\BGs_{1}+\nu\BX,
\]
where
$
\BX=\BGs_{1}^{1/2}\mat{0}{1}{1}{0}\BGs_{1}^{1/2}
$.
Solving for $\BX$ we obtain
\[
\BX=\frac{\BGs_{2}-\mu\BGs_{1}}{\nu}=
\frac{2\BGs_{2}-(\Gl_{1}+\Gl_{2})\BGs_{1}}{\Gl_{1}-\Gl_{2}}=\BS_{2}-\BS_{1}.
\]
Now using $\BJ=\mat{0}{1}{1}{0}\BR_{\perp}$ we get
\begin{multline*}
\BZ=  \BGs_{1}^{1/2}\BJ\BGs_{1}^{1/2}=\BX\BGs_{1}^{-1}\BGs_{1}^{1/2}\BR_{\perp}\BGs_{1}^{1/2}=
\BX\frac{\BR_{\perp}\BGs_{1}\BR_{\perp}^{T}}{\det\BGs_{1}}\sqrt{\det\BGs_{1}}\BR_{\perp}=
\frac{\BX\BR_{\perp}\BGs_{1}}{\sqrt{\det\BGs_{1}}}=\\
\frac{(\BS_{2}-\BS_{1})\BR_{\perp}(\BS_{2}+\BS_{1})}{\sqrt{\det\BGs_{1}}}=
\frac{\BS_{2}\BR_{\perp}\BS_{1}-\BS_{1}\BR_{\perp}\BS_{2}}{\sqrt{\det\BGs_{1}}}
\end{multline*}
The final answer is the equation
\[
(\SFL^{*}+A\SFT)\SFT(\BZ_{0}\otimes\BR_{\perp})\SFT(\SFL^{*}+A\SFT)+B\BZ_{0}\otimes\BR_{\perp}=0,
\qquad\BZ_{0}=\BS_{2}\BR_{\perp}\BS_{1}-\BS_{1}\BR_{\perp}\BS_{2}.
\]

$\bullet$ (2b)
$\BGs_{1}=\Gth_{1}\BGs_{0}$, $\BGs_{2}=\Gth_{2}\BGs_{0}$,
$|r_{1}-r_{2}|>|\Gth_{1}-\Gth_{2}|\sqrt{\det\BGs_{0}}$. The formula for the
effective tensor is
\[
\SFL^{*}=\BGs_{1}\otimes\BL^{*}+t^{*}\SFT,\qquad
\det\BGs_{1}\det\BL^{*}=\left(t^{*}+A\right)^{2}+B,
\]
where $A$ and $B$ are given by (\ref{Asc}) and (\ref{Bsc}), respectively.

You have derived two formulas
\[
\det\BL_{0}^{*}=\left(\frac{a}{\sqrt{\det\BGs_{1}}}t^{*}+b\right)^{2}+1,
\]
and
\[
\det\BL^{*}=\nth{a^{2}}\det\BL_{0}^{*}.
\]
Now combine the two formulas to eliminate $\BL_{0}^{*}$:
\[
\det\BL^{*}=\nth{a^{2}}\left(\frac{a}{\sqrt{\det\BGs_{1}}}t^{*}+b\right)^{2}+\nth{a^{2}}.
\]
Use the formula $x^{2}y^{2}=(xy)^{2}$ in the first terms and leave the second
as it is:
\[
\det\BL^{*}=\left(\frac{t^{*}}{\sqrt{\det\BGs_{1}}}+\frac{b}{a}\right)^{2}+\nth{a^{2}}.
\]
Now multiply both sides by $\det\BGs_{1}$ and use the same formula by writing
$\det\BGs_{1}=(\sqrt{\det\BGs_{1}})^{2}$:
\[
\det\BGs_{1}\det\BL^{*}=\left(t^{*}+\frac{b\sqrt{\det\BGs_{1}}}{a}\right)^{2}+\frac{\det\BGs_{1}}{a^{2}}.
\]
Now let
\[
A=\frac{b\sqrt{\det\BGs_{1}}}{a},\qquad B=\frac{\det\BGs_{1}}{a^{2}}.
\]
Let us use the formulas for $a$ and $b$ to derive explicit formulas for $A$
and $B$:
\[
b=a\frac{\det\BGs-1+\rho_{1}^{2}-\rho_{2}^{2}}{2\rho},\quad
a^{2}=\frac{4\rho^{2}}{((\rho+1)^{2}-\det\BGs)(\det\BGs-(\rho-1)^{2})}.
\]
For compactness of notation let us denote $\GD r=r_{2}-r_{1}$, so that
\[
\rho=\frac{\GD r}{\sqrt{\det\BGs_{1}}}.
\]
We obtain
\[
A=\sqrt{\det\BGs_{1}}\frac{\frac{\det\BGs_{2}}{\det\BGs_{1}}-1+\frac{r_{1}^{2}}{\det\BGs_{1}}-\frac{r_{2}^{2}}{\det\BGs_{1}}}{2\frac{\GD
    r}{\sqrt{\det\BGs_{1}}}}.
\]
Converting to a normal fraction we obtain
\[
A=\sqrt{\det\BGs_{1}}\frac{\sqrt{\det\BGs_{1}}}{\det\BGs_{1}}\frac{\det\BGs_{2}-\det\BGs_{1}+r_{1}^{2}-r_{2}^{2}}{2\GD r}.
\]
The factor in front simplifies to 1 and we obtain
\begin{equation}
  \label{Asc}
  A=\frac{\det\BGs_{2}-\det\BGs_{1}+r_{1}^{2}-r_{2}^{2}}{2\GD r}.
\end{equation}
We perform the same simplification for $B$:
\[
B=\frac{((\GD r+\sqrt{\det\BGs_{1}})^{2}-\det\BGs_{2})(\det\BGs_{2}-(\GD
  r-\sqrt{\det\BGs_{1}})^{2})}{4(\GD r)^{2}}.
\]
Let us make the numerator $nB$ of $B$ a bit more symmetric:
\begin{multline*}
nB=(\GD r+\sqrt{\det\BGs_{1}}-\sqrt{\det\BGs_{2}})
(\GD r+\sqrt{\det\BGs_{1}}+\sqrt{\det\BGs_{2}})\times\\
(\sqrt{\det\BGs_{2}}+\GD r-\sqrt{\det\BGs_{1}})
(\sqrt{\det\BGs_{2}}-\GD r+\sqrt{\det\BGs_{1}}).
\end{multline*}
We now combine first and third terms together and also the other two together:
\[
nB=((\GD r)^{2}-(\sqrt{\det\BGs_{1}}-\sqrt{\det\BGs_{2}})^{2})
((\sqrt{\det\BGs_{1}}+\sqrt{\det\BGs_{2}})^{2}-(\GD r)^{2}).
\]
This gives
\begin{equation}
  \label{Bsc}
  B=\frac{((\GD r)^{2}-(\sqrt{\det\BGs_{1}}-\sqrt{\det\BGs_{2}})^{2})
((\sqrt{\det\BGs_{1}}+\sqrt{\det\BGs_{2}})^{2}-(\GD r)^{2})}{4(\GD r)^{2}}.
\end{equation}

\end{document}